\documentclass[11pt]{article}
\usepackage{amsmath,amsfonts,amssymb,amsthm,graphicx}
\usepackage[arrow, matrix, tips]{xy}

\newtheorem*{theorem*}{Theorem}
\newtheorem{theorem}[subsubsection]{Theorem}
\newtheorem{cor}[subsubsection]{Corollary}
\newtheorem{lem}[subsubsection]{Lemma}
\newtheorem{prop}[subsubsection]{Proposition}

\newtheorem*{statement*}{Statement}
\newtheorem{conj}{Conjecture}
\newtheorem{thesis}{Thesis}

\theoremstyle{definition}
\newtheorem{defi}[subsection]{Definition}
\newtheorem{example}[subsection]{Example}
\newtheorem{remark}[subsection]{Remark}
\newtheorem{notation}[subsubsection]{Notation}
\newtheorem{terminology}[subsubsection]{Terminology}

\numberwithin{equation}{section}

\newcommand{\cob}{{\mathsf{Cob}_{d,(k,N)}^\mathcal{F}}}
\newcommand{\cp}{{\mathbb{C}P}^n}

\title{Geometric cobordism categories}

\author{David Ayala}

\begin{document}

\maketitle

\begin{abstract}
In this paper we study cobordism categories consisting of manifolds which are endowed with geometric structure.  Examples of such geometric structures include symplectic structures, flat connections on principal bundles, and complex structures along with a holomorphic map to a target complex manifold. A general notion of ``geometric structure" is defined using sheaf  
theoretic constructions.  Our main theorem is the identification of the homotopy type of such cobordism categories in terms of certain Thom spectra.  This extends work of Galatius-Madsen-Tillmann-Weiss who identify the homotopy type of cobordism categories of manifolds with fiberwise structures on their tangent bundles.  Interpretations of the main theorem are discussed which have relevance to topological field theories, moduli spaces of geometric structures, and h-principles.  Applications of the main theorem to various examples of interest in geometry, particularly holomorphic curves, are elaborated upon.  
\end{abstract}

\section*{Acknowledgments}
    
I would like to thank my advisor Ralph Cohen for his support and well-tuned efforts.  I am indebted to S{\o}ren Galatius for his guidance and for generously sharing his ideas and intuition.  The inspiration for the ideas of this paper came largely from a course S{\o}ren taught at Stanford University in 2006.  I would like to recognize that some of the results in this paper overlap with those of Oscar Randal-Williams.  Lastly, I am indebted to Randy Montana for at least three instances of incredible navigation and clarity of mind in severe wilderness situations.

\section{Introduction}

\subsection{Context}

Recently developed theory of Madsen and Weiss~\cite{madsen-weiss:mumford} has provided a proof of Mumford's conjecture (\cite{mumford:conjecture}): a conjecture on the cohomology of the moduli space of Riemann surfaces of high genus. In the end, their statement was homotopy theoretic: $\mathbb{Z}\times B\Gamma^+_\infty\simeq \Omega^\infty BU(1)^{-\gamma_2}$. This innovative proof of Madsen and Weiss' introduced several new techniques, making use of singularity theory, h-principles, and categorical methods in homotopy theory. The proof was further distilled by the later work of the four authors Galatius, Madsen, Tillmann, and Weiss in~\cite{galatius-madsen-tillmann-weiss:cobordism} where they define a cobordism category of $d$-manifolds; Mumford's conjecture being a sort of group completion argument (a la~\cite{mcduff-segal:group-completion},\cite{tillmann:mapping-class-group},~\cite{dwyer:centralizer}) on this cobordism category with $d=2$.  The development of this cobordism category is still new and its role in the homotopy theory of manifolds, be it surgery theory, $K$-theory, or stable diffeomorphism groups, is not yet well understood.

This new categorical framework of cobordism theory developed in~\cite{galatius-madsen-tillmann-weiss:cobordism} informs us some on the relationship between Pontryagin-Thom theory and the `stable diffeomorphism group' of $d$-manifolds. Indeed, the four authors show, using Philips' h-principle for submersions (\cite{phillips:submersions}), that it is possible to integrate the tangential data, coming from Pontryagin-Thom theory, of a map to a certain Thom spectrum to the data of a smooth family of $d$-manifolds. In this way, they show that a certain Thom spectrum and the cobordism category of $d$-manifolds classify (via~\cite{weiss:classify} and~\cite{segal:classifying-spaces}) the same objects up to concordance and are thus weakly equivalent.  In fact, as with the spirit of classical Pontryagin-Thom theory, in~\cite{galatius-madsen-tillmann-weiss:cobordism} the four authors identify the homotopy type of the cobordism category of $d$-manifolds endowed with prescribed \textit{tangential} structures.

An insight of Galatius in~\cite{galatius:graphs} was to replace the use of Pontryagin-Thom theory in~\cite{galatius-madsen-tillmann-weiss:cobordism} with scanning arguments (as in~\cite{mcduff:configurations-scanning},\cite{segal:scanning}) applied to graphs in order to study stable automorphisms of free groups, thus avoiding reference to tangent spaces. Proving a conjecture of Hatcher and Vogtmann (\cite{hatcher-vogtmann:splitting-conjecture}), his statement in the end was $\mathbb{Z}\times BAut(F_\infty)^+\simeq QS^0$. Using these methods of scanning suggests that these techniques should hold for spaces which, unlike ordinary manifolds, are not described by tangent data. ÊThe content of this paper concerns manifolds with \textit{geometric} structure and not just tangential structure.   Examples of geometric structures which are not tangential include symplectic structures, flat-connections, the data of a complex structure and a holomorphic map to a target manifold, and many more.

\subsection{Statements of results}

The goal of this paper is to identify the homotopy type of the classifying space of certain cobordism categories.  These categories are founded on $d$-manifolds endowed with certain types of geometric structure.  This geometric structure is encoded by a sheaf $\mathcal{F}$, and we refer to such a cobordism category as $\mathsf{Cob}^\mathcal{F}_d$.  Recall that the four authors in~\cite{galatius-madsen-tillmann-weiss:cobordism} studied cobordism categories of manifolds with fiberwise structures on their tangent bundles.  The results of this paper expand theirs in that the structures at hand are geometric and not just tangential, thus requiring an alternative to the Pontryagin-Thom theory.  We begin by describing what is meant by geometric structure.

Let $\mathsf{Emb}_d$ be a topological category whose objects are $d$-dimensional submanifolds of $\mathbb{R}^\infty$ and whose morphisms are embeddings of $d$-manifolds.  Consider a sheaf
\[
\mathcal{F}:{\mathsf{Emb}_d}^{op}\rightarrow \mathsf{Top}
\]
which is continuous in an appropriate sense.  There results an action of ${Diff}(W)$ on $\mathcal{F}(W)$ for each $W^d\in \mathsf{Emb}_d$.  Think of $g\in\mathcal{F}(W)$ as an $\mathcal{F}$-structure, or geometric structure, on $W$.  The terminology comes from natural examples which include $\mathcal{F}(W)$ given by the following spaces of structures:
\begin{itemize}
\item $\{$Riemannian metrics on $W$, possibly with conditions on the curvature$\}$,

\item  given a fixed complex manifold $Y$, $\{$complex structures $i$ on $W$ along with a holomorphic map $(W,i)\rightarrow Y\}$, 

\item  $\{$symplectic structures on $W\}$,

\item $\{$principal $G$-bundles over $W$ along with a flat connection$\}$, 

\item $\{$immersions from $W$ into a fixed manifold $Y\}$,

\item $\{$unordered configurations of points in $W\}$.  
\end{itemize}
It is suggested that the reader keep such examples in mind.

A rough description of the cobordism category $\mathsf{Cob}^\mathcal{F}_{d,N}$ is given as follows.  Declare
\[
ob\text{ }\mathsf{Cob}^\mathcal{F}_{d,N} =
\{(M^{d-1},g)\}
\]
where $M^{d-1}$ is a closed $(d-1)$-manifold embedded in $\mathbb{R}^{N-1}$ and $g$ is a geometric structure on $(-\epsilon,\epsilon)\times M$ for some $\epsilon >0$ (i.e., $g\in \mathcal{F}( (-\epsilon,\epsilon)\times M^{d-1})$).  Declare
\[
mor\text{ }\mathsf{Cob}^\mathcal{F}_{d,N} = \{identities\}\amalg \{(W^d,g)\}
\]
where $W^d$ is a compact collared $d$-manifold collared-embedded in $[0,1]\times\mathbb{R}^{N-1}$  and $g\in\mathcal{F}(W)$.  The source and target maps are given by intersecting with $\{0\}\times \mathbb{R}^{N-1}$ and $\{1\}\times\mathbb{R}^{N-1}$ respectively.  Composition is given by concatenating cobordisms in $\mathbb{R}^N$ and gluing $\mathcal{F}$-structures via the sheaf-gluing property of the sheaf $\mathcal{F}$.  These categories $\mathsf{Cob}^\mathcal{F}_{d,N}$ are natural in the variable $\mathcal{F}$.  In the case $N=\infty$ the subscript is dropped and we simply write $\mathsf{Cob}^\mathcal{F}_d$.

The category $\mathsf{Cob}^\mathcal{F}_{d,N}$ is naturally a topological category.  For instance, when $N=\infty$ its set of morphisms should classify bundles of collared $\mathcal{F}$-manifolds.   We verify that indeed $mor\text{ }\mathsf{Cob}^\mathcal{F}_d$ has a topology with which its connected components are indexed by diffeomorphism classes of collared $d$-manifolds with $[W^d]$-component homotopy equivalent to the moduli space
\[
\mathcal{F}(W)//{Diff}(W):= E{Diff}(W)\times_{{Diff}(W)} \mathcal{F}(W).
\]
of $\mathcal{F}$-structures on $W$.   
The terminology is justified because this moduli space indeed classifies, up to isomorphism, smooth families of $\mathcal{F}$-manifolds which are diffeomorphic to $W$.  In the above notation, the double forward slash ``//'' denotes the homotopy orbit space.  

\begin{remark}
One could equally well work with topological stacks rather than these homotopy orbit spaces.  Moreover, one could proceed using the language of $\infty$-categories (complete Segal spaces~\cite{rezk:segal-spaces}) rather than topological categories.  See~\cite{hopkins-lurie:tfts} for such a treatment when the geometric structure $\mathcal{F}$ is that of a tangential structure in the sense of~\cite{galatius-madsen-tillmann-weiss:cobordism}. 
\end{remark}

Write $Fr_{d,N}$ for the orthogonal frame bundle over the Grassmann $Gr_{d,N}$ of $d$-planes in $\mathbb{R}^N$.  Write $\gamma_{d,N}$ for the universal rank $d$ vector bundle over $Gr_{d,N}$ and $\gamma^\perp_{d,N}$ as its orthogonal compliment.  There is a projection 
\[
p_\mathcal{F}: Fr_{d,N}\times_{O(d)} \mathcal{F}(\mathbb{R}^d)  \longrightarrow Gr_{d,N}.
\]
Let $Th$ denote the Thom space functor.
We can now state the main theorem.

\begin{theorem}[Main Theorem]\label{thm: main theorem}

There is a weak homotopy equivalence
\[
B\mathsf{Cob}^\mathcal{F}_{d,N} \simeq \Omega^{N-1}Th(p^*_\mathcal{F}\gamma^\perp_{d,N}).
\]

\end{theorem}

Recently, Randal-Williams (\cite{randal-williams:embedding}) proved this result as well in the case of no $\mathcal{F}$-structure.

The classical map $\epsilon^1\oplus \gamma^\perp_{d,N}\rightarrow \gamma^\perp_{d,N+1}$ over the inclusion $Gr_{d,N}\rightarrow Gr_{d,N+1}$ induces structure maps for a spectrum having $Th(p^*_\mathcal{F}\gamma^\perp_{d,N})$ as its $N^{th}$ space.  We call this spectrum $MT\mathcal{F}(d)$.  Using superscript notation to denote the Thom spectrum functor, a quick way of writing this spectrum is as
\[
MT\mathcal{F}(d) := (\mathcal{F}(\mathbb{R}^d)//O(d))^{-p^*_\mathcal{F}\gamma_d}.
\]
To simplify notation, the symbol $p^*_\mathcal{F}$ in the superscript is often dropped.

\begin{notation}
When the sheaf $\mathcal{F}$ is that associated to an understood representation $G(d)\xrightarrow{\rho} O(d)$ in the sense of \textsection\ref{subsec: fiberwise equivariant sheaves}, the spectrum $MT\mathcal{F}(d)$ is written $MTG(d)$.  
\end{notation}

\begin{theorem}[Main Theorem, $N=\infty$]\label{thm: main theorem infinity}

There is a weak homotopy equivalence
\[
B\mathsf{Cob}^\mathcal{F}_d\simeq \Omega^{\infty-1}MT\mathcal{F}(d).
\]
\end{theorem}

\begin{remark}
Though Theorem~\ref{thm: main theorem infinity} is a generalization of the main result of~\cite{galatius-madsen-tillmann-weiss:cobordism}, the techniques in this paper are distinct from those in \cite{galatius-madsen-tillmann-weiss:cobordism} and more closely follow those introduced in
\cite{galatius:graphs}, which provides methods for the case of graphs and an outline  
for the case of smooth manifolds.  Morally speaking, the latter avoids reference to tangent data and Pontryagin-Thom theory, thus allowing the theory to hold for manifolds with non-tangential geometric structure.
\end{remark}

\begin{remark}
As will be made precise, for any such sheaf $\mathcal{F}$, there is associated a \textit{fiberwise} sheaf $\tau\mathcal{F}$ and an h-principle comparison morphism $\mathcal{F}\rightarrow \tau\mathcal{F}$.  A $\tau\mathcal{F}$-structure is an example of a \textit{tangential} structure in the sense of~\cite{galatius-madsen-tillmann-weiss:cobordism} (see \textsection\ref{subsec: fiberwise equivariant sheaves}).  A different phrasing of the main theorem states that this h-principle comparison morphism induces an equivalence on the classifying spaces of the resulting cobordism categories.  
\end{remark}

\begin{remark}
There is an obvious action of the little $(N-1)$-cubes operad on $B\mathsf{Cob}^\mathcal{F}_{d,N}$ obtained upon choosing a homeomorphism $\mathbb{R}^{N-1}\cong int(I^{N-1})$ of the ambient euclidean space of the objects and morphisms of $\mathsf{Cob}^\mathcal{F}_{d,N}$.  The  weak homotopy equivalence in Theorem~\ref{thm: main theorem} above is equivariant with respect to this little $(N-1)$-cubes action.  It follows that the weak homotopy equivalence in Theorem~\ref{thm: main theorem infinity} is an equivalence as $E_\infty$-spaces.  

Better still, as mentioned in \textsection\ref{subsec:delooping-categories} it is possible to define a topological $k$-category $\mathsf{Cob}^\mathcal{F}_{d,(k,N)}$ whose objects are closed $(d-k)$-manifolds in $I^{N-k}$, whose $1$-morphisms are $(d-k+1)$-manifolds in $I^{N-k+1}$ regarded as cobordisms, ... , and whose $k$-morphisms are $d$-manifolds in $I^N$ regarded as cobordisms.  This $k$-category has an action of the little $(N-k)$-cubes operad similar to the situation in the previous paragraph.  This $k$-category is the $(k-1)$-fold delooping of $\mathsf{Cob}^\mathcal{F}_d$ (see Corollary~\ref{cor:delooping}).  Of particular interest is the case $k=d$ which, in the limit $N\rightarrow \infty$, recovers the topological categorical version of the $(\infty,d)$-category of~\cite{hopkins-lurie:tfts}.  
\end{remark}

\subsection{Some examples, motivation, and interpretation}

The inspiration for such generality in the sheaf $\mathcal{F}$ comes from the large list of examples, three of which will be highlighted here.  Each example illustrates a more general phenomenon as will be observed.  See \textsection\ref{subsec: interpretations} for a more complete discussion.

\subsection*{Example 1: stable moduli spaces and holomorphic curves}

It can be expected that the topology of the classifying space $\mathsf{Cob}^\mathcal{F}_d$ contains information about the moduli spaces of $\mathcal{F}$-structures on $d$-manifolds.  After all, recall that each component of the morphism space of $\mathsf{Cob}^\mathcal{F}_d$ is homotopy equivalent to the moduli space $\mathcal{F}(W)//{Diff}(W)$ for some collared $W^d$.

There is a good sense in which the classifying space $B\mathsf{Cob}^\mathcal{F}_d$ can be regarded as the \textit{stable moduli space of $\mathcal{F}$-structures on $d$-manifolds}.
Up to homotopy, there are distinguished maps 
\[
\mathcal{F}(W)//{Diff}(W)\rightarrow \Omega B\mathsf{Cob}^\mathcal{F}_d
\]
such that the operation of gluing $\mathcal{F}$-manifolds along common boundary is compatible with the loop product.  In fact, the classifying space $B\mathsf{Cob}^\mathcal{F}_d$ is initial with respect to this property and as such it can be regarded as encoding the \textit{stable} topology of the moduli spaces of $\mathcal{F}$-structures.  We regard the cohomology of $B\mathsf{Cob}^\mathcal{F}_d$ as parametrizing `stable' characteristic classes of
bundles of $\mathcal{F}$-manifolds.  When appropriate stability exists among the above gluing operations, be it homological stability as in~\cite{harer:stability} or else, the classifying space $B\mathsf{Cob}^\mathcal{F}_d$ \textit{is} the stabilization of the moduli spaces with respect to the gluing operation.  Examples include $B\Gamma^+_\infty$ (\cite{madsen-weiss:mumford}),  the moduli space $\mathcal{M}^d_g(\cp)$ of degree $d$ genus $g$ holomorphic curves in $\cp$ (see \textsection\ref{sec:stability}), and  $B\Sigma^+_\infty$ (\cite{barratt-priddy:symmetric-group},\cite{segal:gamma-spaces}).  See~\cite{hatcher-wahl:stability} for a somewhat general presentation of such homological stability.

Fix a complex manifold $Y$.  Motivated by Gromov-Witten theory, take $d=2$ and $\mathcal{F}(W): = \{(j,h)\}$ where $j$ is a complex structure on $W$ and $(W,j)\xrightarrow{h}Y$ is holomorphic.  The homotopy orbit space $\mathcal{F}(W)//{Diff}(W)$ is the familiar moduli space $\mathcal{M}(Y)$ of holomorphic curves in $Y$.  An application of Theorem~\ref{thm: main theorem infinity} shows that the category of holomorphic curves in $Y$ has the same homotopy type as the category of surfaces mapping \textit{continuously}  to $Y$, namely $\Omega^{\infty-1}(MTSO(2)\wedge Y_+)$.  The cohomology of this $\Omega^\infty$-space is well understood, particularly over the rationals where it is freely generated by the cohomology of $Y$ and the descendant classes $\kappa_i$ each of degree $2i$.  It follows among other things that the `stable' cohomology of the moduli spaces of holomorphic curves in $Y$ does not detect the complex structure of $Y$.

It is of interest to know the connectivity of the maps 
\[
\mathcal{M}_g(Y)\rightarrow \Omega^\infty MTSO(2)\wedge Y_+.
\]
Let $\mathcal{M}^\alpha_g(Y)\subset\mathcal{M}_g(Y)$ be the subspace consisting of holomorphic maps $h$ having homological class $h_*[F]=\alpha\in H_2(Y)$; this subspace is a union of connected components.  In \textsection\ref{sec:stability} we show for $Y={\mathbb{C}P}^n$ and for $\alpha = d>>g>>1$, 
\[
H_k(\mathcal{M}^d_g(\cp);\mathbb{Q})\cong H_k(\Omega^\infty MTSO(2)\wedge {\mathbb{C}P}^n_+;\mathbb{Q}),
\]
the left hand side being well understood.  In general, one expects  $\Omega^\infty MTSO(2)\wedge Y_+$ to approximate the moduli spaces $\mathcal{M}^\alpha_g(Y)$ for $\alpha$ and $g$ large in some sense.  The proof of this theorem follows an argument of G. Segal (\cite{segal:scanning}).

\subsection*{Example 2: h-principles and gauge theory}

Roughly, a smooth map $E^{d+r}\rightarrow X^r$ is a \textit{stable family of $d$-manifolds} if the topology of the $d$-dimensional fibers is allowed to change by way of cobordism.  The space $B\mathsf{Cob}^\mathcal{F}_d$ classifies \textit{stable} families of $d$-dimensional $\mathcal{F}$-manifolds, and as so it parametrizes characteristic classes of stable families.  
Theorem~\ref{thm: main theorem infinity} then states that, up to concordance, all $h$-principles hold in such families.  That is, any stable family of $d$-manifolds which are endowed with \textit{fiberwise} $\mathcal{F}$-structures on their tangent bundles is concordant to another stable family whose fiberwise $\mathcal{F}$-structures come from honest $\mathcal{F}$-structures.

Fix a Lie group $G$ and take $\mathcal{F}$ to be given by $\mathcal{F}(W) = \{(\pi,\omega)\}$ where $\pi$ is a principle $G$-bundle over $W$ and $\omega$ is a flat connection on $\pi$.  Theorem~\ref{thm: main theorem infinity} shows that any family of $d$-manifolds with a connection is concordant to a \textit{stable} bundle of $d$-manifolds with a flat connection.  This offers an alternative proof of results in~\cite{cohen-galatius-kitchloo:flat-connections}.  A  case of particular relevance to four-manifold geometry occurs for a similar sheaf $\mathcal{F}= \mathcal{SD}$ of self-dual connections where the classifying space $B\mathsf{Cob}^\mathcal{SD}_4$ reflects the \textit{stable} topology of the (universal) moduli spaces of self-dual connections.  A similar statement holds when $\mathcal{F}$ is the sheaf of Seiberg-Witten structures.

\subsection*{Example 3: field theories}

A typical field theory is a continuous functor $F: \mathsf{Cob}^\mathcal{F}_d \rightarrow \mathsf{A}$ landing in some, typically algebraic, category $\mathsf{A}$.  Theorem~\ref{thm: main theorem infinity} informs us that no \textit{geometric} invariants are retained upon passing to classifying spaces.  One could phrase this in the context of Hopkins and Lurie's proof of the cobordism hypothesis of Baez and Dolan (\cite{hopkins-lurie:tfts},\cite{baez-dolan:cobordism-hypothesis}) and regard each topological category $\mathsf{Cob}^\mathcal{F}_d$ and possibly $\mathsf{A}$ as an $\infty$-category, thought of as a simplicial space (see \textsection\ref{subsubsec: group completion}), and $F$ as a map of simplicial spaces.  In this case, if $\mathsf{A}$ is an $\infty$-groupoid (i.e., a space), then any such field theory $F$ uniquely factors through $B\mathsf{Cob}^\mathcal{F}_d$; so maps $B\mathsf{Cob}^\mathcal{F}_d\rightarrow \mathsf{A}$ classify field theories $F$ with target the $\infty$-groupoid $\mathsf{A}$.

Of particular relevance to symplectic field theory is the example when $\mathcal{F}(W)$ is the space of symplectic structures on $W$.  Theorem~\ref{thm: main theorem infinity} implies that this category has the same classifying space as the \textit{almost-complex} cobordism category.  The same is conjectured true for the domain category of symplectic field theory.  Another familiar example of a field theory is that of a conformal field theory which is demonstrated when $d=2$ and $\mathcal{F}$ is given by $\mathcal{F}(W^2) = \{\text{complex structures $j$ on $W$}\}$.  Another still, of relevance to Gromov-Witten theory, is demonstrated by $\mathcal{F}(W) = \{$complex structures $j$ on $W$ and holomorphic maps $(W,j)\rightarrow Y\}$.

\subsection{Strategy and organization}\label{subsec: strategy and organization}

There are ten sections in this paper beyond the introduction.
The paper begins with the unfortunately technical task of stating precisely what sort of object the sheaf $\mathcal{F}$ should be.  This, along with useful constructions among such objects, is the content of \textsection\ref{sec: sheaves}.  The difficulty in providing a simple definition comes from the fact that $\mathsf{Top}$-valued sheaves do not admit the ``right" limits or colimits and hence stalks.  As stalks seem to be essential to the ideas involved, a larger target category $\mathsf{QTop}$ for our sheaves is defined, whence the technical discussions.  The light reader is encouraged to skip sections~\ref{subsec: quasi-topological spaces}-~\ref{subsec: continuous sheaves}, referring to them when necessary, and move on to the introduction of equivariant sheaves in \textsection\ref{subsec: equivariant sheaves}.

\textsection\ref{sec: the geometric cobordism category} defines the cobordism category, describes its topology, and  finishes with a precise statement of the main theorem.

\textsection\ref{sec: proof of main theorem} is devoted to the proof of the main theorem whose techniques come largely from the method of proof in~\cite{galatius:graphs} where the objects under study are graphs.  In fact, at the end of~\cite{galatius:graphs} can be found an outline of the proof of the main theorem of this paper when there is no geometric structure $\mathcal{F}$.  An expansion of this outline can be found online from a course Galatius taught at Stanford University in 2006.  We will sketch the idea presently.  
In \textsection\ref{subsec: the sheaf}
we define a space $\Psi^\mathcal{F}_d(\mathbb{R}^N)$ whose points are $d$-dimensional $\mathcal{F}$-submanifolds of $\mathbb{R}^N$ which, as subsets, are closed.  A `scanning'-type argument in the sense of~\cite{segal:scanning} is applied in \textsection\ref{subsec: homotopy type of psi} to identify the weak homotopy type of $\Psi^\mathcal{F}_d(\mathbb{R}^d)$ as the Thom space $Th(p^*_\mathcal{F} \gamma^\perp_{d,N})$.  Group completion techniques in the sense of~\cite{dwyer:centralizer} (also~\cite{mcduff-segal:group-completion} and ~\cite{tillmann:mapping-class-group}) are invoked in \textsection\ref{subsubsec: comparing D and omega B} to identify the iterated loop space $\Omega^{N-1}\Psi^\mathcal{F}_d(\mathbb{R}^N)$ as the  subspace $D^\mathcal{F}_{d,N}$ of $\Psi^\mathcal{F}_d(\mathbb{R}^N)$ consisting of manifolds which are bounded in the last $(N-1)$-directions of $\mathbb{R}^N$.  Applying this group completion is the most difficult part of the proof and occupies~\S\ref{subsec: proof of lemma}.  Easy transversality arguments are applied in \textsection\ref{subsec: bounded manifolds} and \textsection\ref{subsubsec: comparing D and Cob} to identify the weak homotopy type of $D^\mathcal{F}_{d,N}$ as the classifying space $B\mathsf{Cob}^\mathcal{F}_{d,N}$.

\textsection\ref{subsec: interpretations} contains two interpretations of $B\mathsf{Cob}^\mathcal{F}_d$.  The first interprets this classifying space as one which classifies certain families of $d$-dimensional $\mathcal{F}$-manifolds up to concordance.  The second interprets loops on this classifying space as the `stable' moduli space of $\mathcal{F}$-structures in the sense mentioned above.  

\textsection\ref{subsec: fiberwise equivariant sheaves}-\ref{subsec: configurations} are devoted to a few examples of the main theorem, of which there are many.  \textsection\ref{subsec: interpretations}-\ref{subsec: configurations} are intended to use the statement of the main theorem and its setup as a black box without needing knowledge of the proof of the main theorem.  The reader interested only in applications could thus skip straight to \textsection\ref{subsec: interpretations}.

\subsection{Conventions}

Throughout this paper, unless stated otherwise, it will always be meant that a manifold is a smooth manifold and that maps among manifolds are smooth maps.  Spaces consisting of smooth maps are topologized with the $C^\infty$ Whitney topology.  All maps among topological spaces are assumed continuous.  Spaces consisting of continuous maps (between compactly generated spaces) are topologized with the compact-open topology.  All diagrams are commutative.  Each vector space will be assumed a finite dimensional vector space over $\mathbb{R}$ unless stated otherwise or clear from the context.  Write $\mathsf{Set}$ for the category of sets and point-set maps.
Write $\mathsf{Top}$ for the category of topological spaces and continuous maps.  Write $I$ for the standard interval $[0,1]$.  A topological category will be an internal category in $\mathsf{Top}$.  The notion of spectra used here is in the sense of~\cite{adams:stable-homotopy}, also known as prespectra.

Write $\Delta$ for the category of finite linearly ordered sets and order-preserving maps.  There is a tautological functor $\Delta \rightarrow \mathsf{Top}$ given by $\{0,...,n\}\mapsto \Delta^n = \{\underline{t}\in I^{n+1}\mid t_0 + ...+ t_n = 1\}$.  A simplicial space will be a functor $X_\bullet :\Delta^{op}\rightarrow \mathsf{Top}$.  The \textit{classifying space of $X_\bullet$}, or \textit{realization of $X_\bullet$}, is the quotient space 
\[
BX_\bullet = \lVert X_\bullet \rVert : = (\coprod X_l\times \Delta^l)/\sim
\]
where the equivalence relation $\sim$ is generated (only) by the morphisms in $\Delta$ which are injective maps of finite sets.  This notion of classifying space is sometimes referred to as the \textit{thick geometric realization} of $X_\bullet$.  For $\mathsf{C}$ a topological category, define its \textit{simplicial nerve} $N_\bullet \mathsf{C}$ as the simplicial space whose space of $k$-simplicies are $k$-long chains of composable morphisms in $\mathsf{C}$.  That is $N_k\mathsf{C}:= mor \times_{ob}...\times_{ob} mor$.  The \textit{classifying space of $\mathsf{C}$}, written  $B\mathsf{C}$, is the classifying space of the simplicial nerve $N_\bullet \mathsf{C}$.  The \textit{homotopy colimit} of a functor $\mathsf{C}^{op}\xrightarrow{F} \mathsf{Top}$ is the classifying space $B(\mathsf{C}\wr F)$ of the wreath category which has $\{(c,x)\mid c\in ob\text{ }\mathsf{C},\text{ } x\in F(c)\}$ as the space of objects and $mor_{\mathsf{C}\wr F}((c,x),(c^\prime,x^\prime)) := \{c\xrightarrow{m} c^\prime \in mor\text{ }\mathsf{C} \mid F(m)(x^\prime) = x\}$ as the space of morphisms between two objects, appropriately topologized.  Regard a topological group $G$ as a category $\mathsf{G}$ with $ob\text{ }\mathsf{G} = *$ and $mor\text{ } \mathsf{G} = G$ with composition given by the multiplication in $G$.  A functor $\underline{X}:\mathsf{G}\rightarrow \mathsf{Top}$ is the data of a space $X = \underline{X}(*)$ together with a continuous action of $G$.  We denote the homotopy colimit of $\underline{X}$ by $X//G$ and refer to it as the \textit{homotopy orbit space}.

\section{Sheaves}\label{sec: sheaves}

In this section we develop a notion of sheaves sufficient to include most naturally occurring examples of geometric structure while at the same time having a good notion of stalk.  The light reader is encouraged to skip~\textsection\ref{subsec: quasi-topological spaces}-\ref{subsec: continuous sheaves} and simply think ``$\mathsf{Top}$'' whenever ``$\mathsf{QTop}$'' is read.

\subsection{Traditional sheaves}\label{subsec: traditional sheaves}

The familiar, or at least traditional, notion of a sheaf is recalled here.

\subsubsection*{Definition}

For $X\in \mathsf{Top}$ let $\mathcal{O}(X)$ be the  poset of open sets of $X$ ordered by inclusion.  Let $\mathcal{J}(X)$ be the poset of open covers of $M$ ordered by refinement.  Such an open cover $\{U_\alpha\}_{\alpha\in A}\in\mathcal{J}(X)$ determines a diagram (subset) in the poset $\mathcal{O}(X)$ whose entries (elements) are $\cap_{r\in R} U_r$ where $\emptyset \neq R\subset A$ is finite.

\begin{defi}\label{def: traditional sheaf}
A  \textit{traditional} ($\mathsf{Set}$-valued) presheaf on a space $X$ is a functor
\[
\mathcal{T}: {\mathcal{O}(X)}^{op}\rightarrow \mathsf{Set}.
\]
The sheafification of the traditional presheaf $\mathcal{T}$ is the presheaf $\overline{\mathcal{T}}$
given by 
\[
\overline{\mathcal{T}}(U):=colim_\mathcal{J} \text{ } lim \mathcal{T}(\cap_{r\in R} U_r)
\]
where the first colimit is over $\{U_\alpha\}_{\alpha\in A}\in\mathcal{J}(U)$ and the second limit is over finite non-empty subsets $R\subset A$.  We say $\mathcal{T}$ is a \textit{traditional sheaf} if the universal morphism 
\[
\mathcal{T}(U)\rightarrow \overline{\mathcal{T}}(U) 
\]
is an equivalence for each $U\in\mathcal{O}(X)$.
\end{defi}

\begin{example} The prototypical example of a traditional sheaf comes in the situation of a fiber bundle $\alpha= (E\rightarrow X)$ as the assignment $U\mapsto \Gamma(\alpha_{\mid_U})$, the set (space) of sections of $\alpha_{\mid_U}$.  It is straight forward to verify that this is indeed a traditional sheaf on $X$.  
\end{example}

Extend a traditional ($\mathsf{Set}$-valued) sheaf $\mathcal{T}$ on $X$ to arbitrary subsets $A\subset X$ by 
\[
\mathcal{T}(A) := colim_{A\subset U} \mathcal{T}(U) 
\]
where the colimit is over $A\subset U\in \mathcal{O} (X)$.  But more, giving $A\subset X$ the subspace topology, create a new traditional sheaf on $A$,
\[
\mathcal{T}_{\mid_A}:{\mathcal{O}(A)}^{op}\rightarrow \mathsf{Set},
\hspace{.7cm} \text{by} \hspace{.7cm}  
\mathcal{T}_{\mid_A}(V):= colim \mathcal{T}(U)
\]
where the colimit is over the poset $\{U\in \mathcal{O}(X)\mid U\cap A = V\}$.  An important example occurs when $A=\{x\}\subset X$; call the set $\mathcal{T}_{\mid_{\{x\}}}(\{x\})$ the \textit{stalk} of $\mathcal{T}$ at $x\in X$ and write it as $Stalk_x\mathcal{T}$.  There is a universal map $\mathcal{T}(U)\rightarrow Stalk_x\mathcal{T}$ for any $x\in U\in \mathcal{O}(X)$; denote the value of this map at $g\in \mathcal{T}(U)$ by $g_{\mid_x}$ and refer to it as the \textit{germ} of $g$ at $x$.  Note that by definition, for any $x\in A\subset X$, $Stalk_x\mathcal{T} = Stalk_x(\mathcal{T}_{\mid_A})$.

\begin{defi}

A traditional $\mathsf{Top}$-valued sheaf on $X$ is a functor 
\[
\mathcal{T}: \mathcal{O}(X)^{op}\rightarrow \mathsf{Top}
\]
which satisfies a sheaf-gluing condition similar to that in Definition~\ref{def: traditional sheaf}.

\end{defi}

\begin{remark}
The naturally occurring examples of sheaves in this paper tend to be $\mathsf{Top}$-valued sheaves.  However, the category of $\mathsf{Top}$-valued sheaves does not allow for familiar constructions we would like to use in this paper.  We are inclined to adopt the weaker notion of a \textit{continuous sheaf} very similar to that introduced by Gromov~\cite{gromov:h-principle}.   It is for this reason that we involve quasi-topological spaces.
\end{remark}

\subsection{Quasi-topological spaces}\label{subsec: quasi-topological spaces}

Fix a distinguished one-point space $\{*\}$.  For $X\in\mathsf{Top}$, each $x\in X$ determines an inclusion $\iota_x:\{*\}\rightarrow X$ by $*\mapsto x\in X$.  Given a functor 
\[
\mathcal{A}:\mathsf{Top}^{op}\rightarrow \mathsf{Set}
\]
with $A:=\mathcal{A}(\{*\})$, there is a canonical point-set map $\mathcal{A}(\{*\})\rightarrow \{X\rightarrow A\}$ described by $a\mapsto(x\mapsto \iota^*_x a)$.

\begin{defi}
A Quasi-topological space, Q-space for short, is a functor 
\[
\mathcal{A}:\mathsf{Top}^{op}\rightarrow \mathsf{Set}
\]
such that 
\begin{itemize}
\item{for each $X\in \mathsf{Top}$ the restriction $\mathcal{O}(X)^{op}\xrightarrow{\mathcal{A}_{\mid_X}}\mathsf{Set}$ is a traditional ($\mathsf{Set}$-valued) sheaf on $X$},
\item{the natural point-set map $\mathcal{A}(\{*\})\rightarrow \{X\rightarrow A\}$ is injective,}
\item{$\mathcal{A}$ takes filtered colimits to filtered limits.}
\end{itemize}
\end{defi}

A morphism of quasi-topological spaces is simply a natural transformation of functors.

\begin{remark}
It should be noted that this is not the same as the definition of a quasi-topological space due to~\cite{gromov:h-principle}, the difference being that we leave off an extra axiom which requires that for each space $X$, $\mathcal{A}_{\mid_X}$ satisfies the sheaf condition for \textit{finite closed} covers of $X$.
\end{remark}

\subsubsection*{Notation and terminology}

Observe the evaluation functor $ev_{\{*\}}:\mathsf{QTop}\rightarrow \mathsf{Set}$ given as $\mathcal{A}\mapsto A:= \mathcal{A}(\{*\})$.
Refer to the point-set maps in 
\[
\mathcal{A}(X)\subset\{X\rightarrow A\}
\]
as \textit{``continuous"}; the idea being that the set $A$ does not have a topology but does have a notion of continuous maps into it.  For this reason, the set $A$ will often be referred to as the Q-space and a functor $\mathcal{A}$ with $\mathcal{A}(\{*\}) = A$ as a \textit{quasi-topology} on $A$.

To say $\mathcal{A}:\mathsf{Top}^{op}\rightarrow \mathsf{Set}$ is a functor is to say in particular that composition of a continuous function with a ``continuous" function is ``continuous"; to say $\mathcal{A}_{\mid_X}$ is a sheaf is to say that a point-set map $X\rightarrow A$ is ``continuous" exactly when it is ``continuous" and compatible on each $U_\alpha\subset X$ of an open cover of $X$.  In most situations, the quotations will be dropped  and we simply refer to a map $(X\rightarrow A)\in\mathcal{A}(X)$ as continuous, the appropriate notion being understood.  An isomorphism of quasi-topological spaces is called a \textit{homeomorphism}.

The above terminology is inspired by the situation when a set $A$ is given a topology and 
\[
\mathcal{A}(X)=Map(X,A),
\]
the set (space) of continuous maps from $X$ to $A$; in this situation we say $\mathcal{A}$ is \textit{represented}.  Indeed, there is a Yoneda embedding $\mathsf{Top}\hookrightarrow \mathsf{QTop}$ given by sending s space $A$ to the representing functor $Map(-,A)$.

\subsubsection*{Subobjects and mapping objects}

Given a Q-space $A$ and a subset $A^\prime\subset A$, there is a naturally induced quasi-topology on $A^\prime$ as follows.  A point-set map $X\rightarrow A^\prime$ is declared ``continuous" exactly when the composite point-set map $X\rightarrow A^\prime\hookrightarrow A$ is ``continuous" as a map to $A$.  The inclusion $A^\prime \hookrightarrow A$ is then a continuous map of Q-spaces.

Given a Q-space $\mathcal{A}$ and a space $Z$, the set $\mathcal{A}(Z)$ is also a Q-space seen as follows.  Declare a point-set map $X\rightarrow \mathcal{A}(Z)$ to be ``continuous" if the adjoint map (as sets) $X\times Z\rightarrow A$ is ``continuous".  Clearly then, a continuous map of spaces $W\rightarrow Z$ induces a morphism, or continuous map, of quasi-topological spaces $\mathcal{A}(Z)\rightarrow\mathcal{A}(W)$.  In a similar way, given a Q-space $A$, a space $Z$, and a continuous map $A\xrightarrow{p} Z$, the set of sections 
\[
\Gamma(A\xrightarrow{p} Z):=\{Z\xrightarrow{f} A\mid f\text{ is continuous and } p\circ f = id_Z\}
\]
is naturally a quasi-topological space.  Moreover, the restriction maps 
\[
\Gamma(A\rightarrow Z)\rightarrow \Gamma(A_{\mid_U}\rightarrow U)
\]
 are continuous for each $U\in\mathcal{O}(Z)$.

\subsubsection*{Limits and colimits in $\mathsf{QTop}$}\label{subsec: limits and colimits}

Here we point out the driving feature of Q-spaces as admitting well-behaved limits and colimits.

Let $\mathsf{D}$ be any small category and 
\[
F:\mathsf{D}\rightarrow \mathsf{QTop}
\]
a functor.  Write $(co)lim(F)$ for a choice of the categorical (co)limit of the composition of $F$ with the functor $ev_{\{*\}}$ to $\mathsf{Set}$.  For $X$ a space, declare $X\rightarrow lim(F)$ to be ``continuous" if for each $d\in\mathsf{D}$ the composite map
\[
X\rightarrow lim(F)\rightarrow F(d)
\]
 is ``continuous".  So $lim(F)(X):=lim(ev_X\circ F)$ where $\mathsf{QTop}\xrightarrow{ev_X}\mathsf{Set}$, given by $A\mapsto A(X)$, is evaluation on $X$.  It is apparent that this is the weakest quasi-topology on $lim(F)$ such that the universal maps $lim(F)\rightarrow F(d)$ are continuous.  So in fact, $lim(F)$, with this quasi-topology, satisfies the required universal property as a categorical limit in $\mathsf{QTop}$.

Similarly, we endow $colim(F)$ with the weakest quasi-topology so that the universal maps $ F(d)\rightarrow colim(F)$ are continuous.  From this the required universal property for $colim(F)$ in $\mathsf{QTop}$ will follow.  Declare $X\rightarrow colim(F)$ to be ``continuous" if the pull back 
\[
X\times_{colim(F)} F(d)\rightarrow F(d),
\]
is continuous for each $d\in\mathsf{D}$ where the source is a quasi-topological space when regarded as a subset of the categorical limit $X\times F(d)$ in $\mathsf{QTop}$.  It is worth checking that this is consistent with the categorical (co)limit in $\mathsf{Top}$ under the Yoneda-embedding $\mathsf{Top}\hookrightarrow \mathsf{QTop}$.

It can be verified that the Yoneda embedding $\mathsf{Top}\hookrightarrow \mathsf{QTop}$ is continuous (preserves limits) and cocontinuous (preserves colimits).

\subsubsection*{Homotopy theoretic notions in $\mathsf{QTop}$}

Here we establish that many constructions from classical homotopy theory carry over to the category of Q-spaces.  Namely, the notion of homotopy, homotopy equivalence, homotopy groups, and weak homotopy equivalence.

Although not every quasi-topological space is representable  by a topological space on the nose, it is up to homotopy.  To be more precise, given a Q-space $\mathcal{A}$, the assignment $[n]\mapsto \mathcal{A}(\Delta^n)$ is a simplicial-set denoted $\mathcal{A}_\bullet$.  The realization of this simplicial set yields a topological space.  This process is sufficiently natural to describe a functor $\mathsf{QTop}\rightarrow \mathsf{Top}$ so that the composition 
\[
\mathsf{Top}\hookrightarrow\mathsf{QTop}\rightarrow\mathsf{Top}
\]
is level-wise a weak homotopy equivalence; amounting to the classical weak homotopy equivalence $|S_\bullet X | \simeq X$ for $S_\bullet X$ the simplicial set of singular simplicies of $X$.  We say a morphism $\mathcal{A}^\prime\rightarrow\mathcal{A}$ is a \textit{homotopy equivalence} if the induced map of simplicial sets $\mathcal{A}^\prime_\bullet\rightarrow \mathcal{A}_\bullet$ realizes to a homotopy equivalence.

For $\mathcal{A}\in\mathsf{QTop}$, two elements $X\xrightarrow{f_0,f_1} A\in \mathcal{A}(X)$ are said to be \textit{homotopic} if there is an element $(I\times X\xrightarrow{H} A) \in\mathcal{A}(I\times X)$ such that $H_{\mid_{\{\nu\}\times X}} = f_\nu$ ($\nu=0,1$).  This notion of homotopy is an equivalence relation on $\mathcal{A}(X)$ and is natural in the argument $\mathcal{A}$.  Define $\mathcal{A}[X] :=\mathcal{A}(X)/homotopy$ to be the set of homotopy classes of ``continuous" maps $X\rightarrow A$.  We could alternatively say that a morphism of Q-spaces $\mathcal{A}^\prime\rightarrow\mathcal{A}$ is a homotopy equivalence, respectively weak homotopy equivalence, if $\mathcal{A}^\prime[X]\rightarrow \mathcal{A}[X]$ is a bijection for every space $X$, respectively every finite CW-complex $X$.  The equivalence relation generated by both is written as $\simeq$.

Define the homotopy groups of a Q-space as $\pi_q(\mathcal{A}):=\mathcal{A}[S^q]$ which are indeed groups for standard reasons.  A sequence of maps of Q-spaces 
\[
\mathcal{A}^\prime\rightarrow \mathcal{A}\rightarrow\mathcal{A}^{\prime\prime}
\]
is said to be a fibration sequence or simply a fibration, respectively a quasi-fibration, if the homotopy lifting property holds, respectively if it induces a long exact sequence in homotopy groups as usual.  In both cases, $\mathcal{A}^\prime$ is referred to as the fiber of the (quasi-)fibration.

With limits, colimits, and this notion of homotopy, $\mathsf{QTop}$ admits homotopy limits and homotopy colimits in the usual way.  Briefly, one can regard the homotopy (co)limit of a functor as the (co)limit of an associated functor built using iterated homotopies.

Though such language is unnecessary, one could alternatively describe a model structure on $\mathsf{QTop}$ whose weak equivalences are pulled back from the standard weak equivalences in $\mathsf{Set}^{\Delta^{op}}$.  As such, it can be verified that the Yoneda embedding $\mathsf{Top}\rightarrow \mathsf{QTop}$ induces an isomorphism on homotopy categories.

\subsection{Continuous sheaves}\label{subsec: continuous sheaves}

Here we introduce ideas from Gromov~\cite{gromov:h-principle} on \textit{continuous} sheaves.  Think of the category of such as an enlargement of the category of traditional $\mathsf{Top}$-valued sheaves which admits well-behaved colimits.

\begin{defi}
A traditional continuous presheaf on $X$ is a functor 
\[
\mathcal{T}:\mathcal{O}(X)^{op}\rightarrow \mathsf{QTop}.
\]
The sheafification of a presheaf $\mathcal{T}$ is the traditional continuous presheaf $\overline{T}$ determined by
\[
\overline{\mathcal{T}}(U):=colim_\mathcal{J} \text{ } lim \mathcal{T}(\cap_{r\in R} U_r)
\]
where the colimit limit is as before (note that this colimit limit exists in $\mathsf{QTop}$).  A presheaf $\mathcal{F}$ is said to be a sheaf if the universal morphism $\mathcal{T}(U)\rightarrow \overline{\mathcal{T}}(U)$ is a homeomorphism for each $U\in\mathcal{O}(X)$.
\end{defi}

A morphism $\mathcal{T}^\prime\rightarrow\mathcal{T}$ of traditional continuous sheaves on $X$ is a natural transformation of functors.  Such a morphism is called a \textit{homeomorphism} if it is level-wise a homeomorphism, that is $\mathcal{T}^\prime(U)\rightarrow\mathcal{T}(U)$ is a homeomorphism for each $U\in\mathcal{O}(X)$.  Similarly, a morphism is called a \textit{homotopy equivalence}, or \textit{weak homotopy equivalence}, if $\mathcal{T}^\prime(U)\rightarrow \mathcal{T}(U)$ is a homotopy equivalence, or weak homotopy equivalence respectively, for each $U\in\mathcal{O}(X)$.

The Yoneda-embedding $\mathsf{Top}\hookrightarrow \mathsf{QTop}$ results in  an inclusion functor from the category of traditional $\mathsf{Top}$-valued sheaves on $X$ to the category of traditional continuous sheaves on $X$.

\begin{example}
A topological group $G$ is said to \textit{act ``continuously"} on a Q-space $A$ if the point-set action $G\times A\rightarrow A$ is a ``continuous" map of Q-spaces.   The quotient $A/G$, regarded as a colimit, is naturally a Q-space.

Given such a Q-space $A$ with an action of a topological group $G$ and given a principal $G$-bundle $(P\rightarrow B)$ over a topological space $B$, one can form the quasi-topological space $P\times_G A$ and the \textit{quasi-fiber bundle} 
\[
\xymatrix{
P\times_G A  \ar[d]
\\
B.
}
\]
For $\alpha=(E\rightarrow B)$ such a quasi-fiber bundle, the assignment 
\[
(U\subset B)\mapsto \Gamma(\alpha_{\mid_U})
\]
is a traditional continuous sheaf on $B$.  Call this sheaf \textit{the sheaf of sections of $\alpha$} and write it as $\Gamma(\alpha_{\mid})$.  This is the prototypical example of a traditional continuous sheaf on $X$.
\end{example}

\begin{remark}
Because the target category of a continuous sheaf $\mathcal{T}$ on $X$ is $\mathsf{QTop}$, we have a good notion of Stalk.  That is, for each $x\in X$, $colim_{\{x\in U\}} \mathcal{T}(U)$ is a well-behaved Q-space exactly because (good) colimits exist in $\mathsf{QTop}$.  
\end{remark}

\subsection{Equivariant sheaves}\label{subsec: equivariant sheaves}

Here we finally land at the notion of sheaves used in the remainder of the paper.  Think of the category of equivariant sheaves as an enlargement of the category of $\mathsf{Top}$-valued sheaves which admits well-behaved stalks and whose evaluation on an object has a continuous action of the automorphism group of that object.

Let $\mathsf{Emb}_{d,N}$ be the category whose objects are embedded $d$-manifolds $W^d\subset\mathbb{R}^N$ without boundary and whose morphisms are embeddings of manifolds.  Particularly, 
\[
ob\text{ }\mathsf{Emb}_{d,N} = \coprod_{[W]} Emb(W,\mathbb{R}^N)/{Diff}(W)
\]
where the disjoint union is over all diffeomorphism classes of boundaryless $d$-manifolds.    Similarly,
\[
mor\text{ }\mathsf{Emb}_{d,N}=\coprod_{[W_0],[W_1]} Emb(W_0,\mathbb{R}^N)\times_{{Diff}(W_0)}Emb(W_0,W_1)\times_{{Diff}(W_1)} Emb(W_1,\mathbb{R}^N)
\]
where the disjoint union is over all pairs of diffeomorphism classes of boundaryless $d$-manifolds, and the action of each diffeomorphism groups is by pre- and post-composition.  For reasons discussed in \textsection\ref{subsubsec: quasi-topology of cob}, $\mathsf{Emb}_{d,N}$ is naturally a topological category.  When $N=\infty$ drop the subscript and topologize with the colimit topology.  Regard $\mathsf{Emb}_d^{op}$ as a category enriched in $\mathsf{QTop}$.

\begin{defi}
An equivariant presheaf is an enriched functor of categories enriched in $\mathsf{QTop}$
\[
\mathcal{F}:{\mathsf{Emb}_d}^{op}\rightarrow \mathsf{QTop}.
\]
A morphism of equivariant continuous presheaves is simply a natural transformation of functors.
A morphism of equivariant sheaves is a homeomorphism, or (weak) homotopy equivalence, if it is a level-wise homeomorphism, respectively (weak) homotopy equivalence, of quasi-topological spaces.  
A global section $g\in\mathcal{F}(W)$ is an $\mathcal{F}$-structure on $W$.  A pair $(W,g)$, $g\in\mathcal{F}(W)$, is called an $\mathcal{F}$-manifold.  
\end{defi}

For each $W\in ob\text{ }\mathsf{Emb}_d$, there is the obvious inclusion functor
\[
\iota_W: \mathcal{O}(W)\rightarrow \mathsf{Emb}_d.
\]
Given an equivariant presheaf, the composition 
\[
\iota_M^*\mathcal{F}:\mathcal{O}(W)^{op}\xrightarrow{\iota_W}{\mathsf{Emb}_d}^{op}\xrightarrow{\mathcal{F}} \mathsf{QTop}
\]
is a traditional continuous presheaf on $W$, sometimes written $\mathcal{F}_{\mid_W}$.

\begin{defi}
The sheafification of an equivariant presheaf $\mathcal{F}$ is the equivariant presheaf $
\overline{\mathcal{F}}$ determined by $\iota^*_W\overline{\mathcal{F}}:=\overline{\iota^*_W\mathcal{F}}$ for each $W\in \mathsf{Emb}_d$.
We say an equivariant presheaf $\mathcal{F}$ is an equivariant sheaf if the universal morphism
\[
\mathcal{F}\rightarrow \overline{\mathcal{F}}
\]
is a level-wise homeomorphism. 
\end{defi}

Note that for $\mathcal{F}$ an equivariant sheaf, ${Diff}(W)$ acts on the space $\mathcal{F}(W)$ with the action map ${Diff}(W)\times \mathcal{F}(W)\rightarrow \mathcal{F}(W)$ ``continuous" (the group ${Diff}(W)$ is a \textit{topological} group as explained in \textsection\ref{subsubsec: quasi-topology of cob}
).

\subsection{Fiberwise sheaves}\label{subsec: fiberwise sheaves}

Here we single out those equivariant sheaves which are given from fiberwise structure on tangent bundles.  It is this class of equivariant sheaves to which all others compare via an h-principle comparison morphism.  It is this class of equivariant sheaves for which the main theorem was proved in~\cite{galatius-madsen-tillmann-weiss:cobordism}.

\subsubsection*{Traditional sheaves of sections}\label{subsec: trad sheaves of sections}

Let $\alpha=(E\rightarrow B)$ be a $d$-dimensional vector bundle with associated principal $GL(d)$-bundle written as $P_\alpha\rightarrow B$.  Recall that ${Diff}(\mathbb{R}^d)$ acts on $\mathcal{F}(\mathbb{R}^d)$ and thus so does $GL(d)$.  Consider the resulting bundle 
\[
\xymatrix{
P_{\alpha}\times_{GL(d)} \mathcal{F}(\mathbb{R}^d)  \ar[d]
\\
B
}
\]
written as $\mathcal{F}(\alpha)$.  
The assignment $U\mapsto  \Gamma (\mathcal{F}(\alpha)_{\mid_U})$ is a traditional continuous sheaf on $B$ as noted in the example of \textsection\ref{subsec: limits and colimits} where it was written $\Gamma(\mathcal{F}(\alpha)_{\mid})$.

More to the point, because $GL(d)$ acts on the poset $\{0\in U\subset \mathbb{R}^d\}$ of neighborhoods of the origin ordered by inclusion, then it acts on the colimit $Stalk_0(\mathcal{F}_{\mid_{\mathbb{R}^d}})$.  Write the resulting quasi-bundle over $B$, 
\[
\xymatrix{
P_{\alpha}\times_{GL(d)} Stalk_0(\mathcal{F}_{\mid_{\mathbb{R}^d}})  \ar[d]
\\
B,
}
\]
as $Stalk_0(\mathcal{F}(\alpha))$.  Again, the sheaf of sections $U\mapsto\Gamma(Stalk_0(\mathcal{F}(\alpha))_{\mid_U})$, written $\Gamma(Stalk_0(\mathcal{F}(\alpha))_{\mid})$, is a traditional continuous sheaf on $B$.

There is a morphism of traditional continuous sheaves on $B$
\[
\Gamma(\mathcal{F}(\alpha)_\mid)\rightarrow\Gamma(Stalk_0(\mathcal{F}(\alpha))_\mid)
\]
described by sending $g\in\mathcal{F}(E_b)$ to $g_{\mid_0}\in Stalk_0(\mathcal{F}(E_b))$.

\begin{prop}\label{prop: stalk}

For $\alpha=(E\xrightarrow{\pi} B)$ a vector bundle, the morphism above is a weak homotopy equivalence of traditional continuous sheaves on $X$,
\[
\Gamma(\mathcal{F}(\alpha)_\mid)\xrightarrow{\simeq}\Gamma(Stalk_0(\mathcal{F}(\alpha))_\mid)
\]

\end{prop}

\begin{proof}

Choose a metric on $\alpha$ so that its structure group lies in $O(d)$.  The space of such metrics is affine and thus contractible.  Without loss in generality, it only needs to be  verified that the map
\[
\Gamma(\mathcal{F}(\alpha))\rightarrow\Gamma(Stalk_0(\mathcal{F}(\alpha)))
\]
is a weak homotopy equivalence.

Let $K$ be a finite, and therefore compact, CW-complex.  
Let $D_\epsilon\subset\mathbb{R}^d$ be an $\epsilon$-neighborhood of the origin with $D_\infty=\mathbb{R}^d$.  Choose once and for all a \textit{continuous} family of embeddings $ \mathbb{R}^d  \xrightarrow{\phi_\epsilon}    \mathbb{R}^d$, $\epsilon\in (0,\infty]$, with $\phi_\infty = id_{\mathbb{R}^d}$ and such that each $\phi_\epsilon$ has image $D_\epsilon$ and is the identity in a small neighborhood of the origin.

Because open disks form a basis for the topology of $\mathbb{R}^d$, the universal map
\[
colim_\epsilon Map(-,P_\alpha\times_{O(d)}\mathcal{F}(D_\epsilon))\rightarrow Stalk_0(\mathcal{F}_{\mid_{\mathbb{R}^d}})
\]
is a homeomorphism of quasi-topological spaces.  From the definition of stalk as a colimit in $\mathsf{QTop}$, the set of ``continuous" maps $Map(K,\Gamma(Stalk_0(\mathcal{F}(\alpha)))$ is the set 
\[
colim_\epsilon \{K\times B\xrightarrow{f} P_\alpha\times_{O(d)}D_\epsilon\mid \forall x\in K,\text{ }\pi\circ f_{\mid_{\{x\}\times B}} = id_B\}.
\]
For each $f\in Map(K,\Gamma(Stalk_0(\mathcal{F}(\alpha)))$ the compactness of $K$ ensures the existence of $1\geq \epsilon_f>0$ such that $f$ is represented by $f_{\epsilon_f}\in Map(K\times B,P_\alpha\times_{O(d)}D_{\epsilon_f})$.  Pulling-back by the embedding $\phi_{\epsilon_f}$ yields $\tilde{f_{\epsilon_f}}\in Map(K,\Gamma(\mathcal{F}(\alpha)))$ with the property that 
\[
(\tilde{f_{\epsilon_f}})_{\mid_0} = f\in Map(K,\Gamma(Stalk_0(\mathcal{F}(\alpha))).
\]

On the other hand, given $\tilde{f}\in Map(K,\Gamma(\mathcal{F}(\alpha)))$,  performing the above procedure  to $f=\tilde{f}_{\mid_0}\in Map(K,\Gamma(Stalk_0(\mathcal{F}(\alpha)))$ we get $\tilde{f}_\epsilon\in  Map(K,\Gamma(\mathcal{F}(\alpha)))$ for some $1\geq \epsilon>0$.  The family $\phi_{\epsilon^\prime}$, $\epsilon^\prime \in [\epsilon,\infty]$, of embeddings realizes a homotopy from $\tilde{f}_\epsilon$ to $f$.  This demonstrates the desired bijection 
\[
[K,\Gamma(\mathcal{F}(\alpha))]\leftrightarrow [K,\Gamma(Stalk_0(\mathcal{F}(\alpha))].
\]
\end{proof}

\begin{remark}

If $\mathcal{F}$ has values in $\mathsf{Top}$ then the sheaf $\Gamma(\mathcal{F}(\alpha)_\mid)$ also has values in $\mathsf{Top}$.  On the other hand,
the sheaf $\Gamma(Stalk_0(\mathcal{F}(\alpha)_\mid)$ is not in general representable even if $\mathcal{F}$ is.  This is explained by the observation that the quasi-topological space $Stalk_0(\mathcal{F}_{\mid_{\mathbb{R}^d}})$ is not representable even if $\mathcal{F}$ is.  Indeed, for a generic space $X$,
\[
colim_\epsilon Map(X,\mathcal{F}(D_\epsilon)) \neq Map(X,colim_\epsilon \mathcal{F}(D_\epsilon))
\]
where $D_\epsilon\subset\mathbb{R}^d$ is an $\epsilon$-neighborhood of the origin.
\end{remark}

\begin{defi}[First Definition of $\tau\mathcal{F}$]
Applying the above procedure with $\alpha = \tau_W$ for $W\in \mathsf{Emb}_d$, from $\mathcal{F}$ we  obtain a new equivariant sheaf 
\[
\tau\mathcal{F}:{\mathsf{Emb}_d}^{op}\rightarrow \mathsf{QTop}
\]
determined by $\tau\mathcal{F}(W)=\Gamma(Stalk_0(\mathcal{F}(\tau_W)))$ since for any embedding $W^\prime\xrightarrow{e}W$, $\tau_{W^\prime} \xrightarrow{e_*} e^*\tau_W$ is a bundle isomorphism.  
\end{defi}

\subsubsection*{Alternative description of $\tau\mathcal{F}$}\label{subsubsec: alternative description}

Let $Gr_{d,N}$ denote the Grassmann manifold of $d$-planes in $\mathbb{R}^N$.  Write $Fr_{d,N}$ for the Stiefel space of linear orthogonal embeddings $\mathbb{R}^d\hookrightarrow \mathbb{R}^N$.  Let $Fr_{d,N}^\mathcal{F}$ denote the set of pairs 
\[
\{(e,g)\mid e\in Fr_{d,N},\text{   }g\in Stalk_0(\mathcal{F}_{\mid_{e(\mathbb{R}^d)}})\}.
\]
Endow $Fr_{d,N}^\mathcal{F}$ with the quasi-topology according to the bijection 
\[
Fr_{d,N}^\mathcal{F}\leftrightarrow Fr_{d,N}
\times
Stalk_0(\mathcal{F}_{\mid_{\mathbb{R}^d}})
\]
given by $(e,g)\mapsto (e,e^*g)$.
With this quasi-topology there is a continuous free action of $O(d)$.  Give the set $Gr_{d,N}^\mathcal{F}:= Fr_{d,N}^\mathcal{F}/O(d)$ the resulting colimit quasi-topology.  It is clear that the projection $Gr_{d,N}^\mathcal{F}\xrightarrow{p_\mathcal{F}} Gr_{d,N}$ is continuous.  Choose $Gr_{d,\infty}$ as a model for $BO(d)$ using the contractibility of $Fr_{d,\infty}$ along with its free $O(d)$-action.  Write $B\mathcal{F}(d)$ for the colimit 
\[
B\mathcal{F}(d):=Gr^\mathcal{F}_{d,\infty} \cong Stalk_0(\mathcal{F}_{\mid_{\mathbb{R}^d}})//O(d).
\]

\begin{prop}
The projection
\[
p_\mathcal{F}:Gr_{d,N}^\mathcal{F}\rightarrow Gr_{d,N}
\]
is a fibration with fiber $Stalk_0(\mathcal{F}_{\mid_{\mathbb{R}^d}})$
\end{prop}

\begin{proof}
 
For $N=\infty$, the $O(d)$-equivariance of the above bijection establishes the homotopy equivalence
\[
B\mathcal{F}(d)\xrightarrow{\cong}Fr_{d,\infty}\times_{O(d)} Stalk_0(\mathcal{F}_{\mid_{\mathbb{R}^d}}) =  Stalk_0(\mathcal{F}_{\mid_{\mathbb{R}^d}})//O(d)
\]
over $BO(d)$.  
This establishes the fibration sequence
\[
Stalk_0(\mathcal{F}_{\mid_{\mathbb{R}^d}})\rightarrow B\mathcal{F}(d)\rightarrow BO(d).
\]

For finite $N$, there is the pull back square 
\[
\xymatrix{
Gr_{d,N}^\mathcal{F}  \ar[r]  \ar[d]
&
B\mathcal{F}(d)  \ar[d]
\\
Gr_{d,N} \ar[r]
&
BO(d).
}
\]
Thus
\[
Stalk_0(\mathcal{F}_{\mid_{\mathbb{R}^d}})\rightarrow Gr_{d,N}^\mathcal{F}\rightarrow Gr_{d,N}
\]
is a fibration sequence.
\end{proof}

Given $W\in\mathsf{Emb}_d$ there is an explicit map $W\xrightarrow{\tau_W}BO(d)$ classifying $\tau_W$, the tangent bundle of $W$.  The Q-space $\Gamma(Stalk_0(\mathcal{F}(\tau_W))$ is then naturally homeomorphic to the Q-space 
\[
\{l:W\rightarrow B\mathcal{F}(d)\mid p_\mathcal{F}\circ l =\tau_W\}
\]
of lifts 
\[ 
\xymatrix{
&
B\mathcal{F}(d)  \ar[d]^{p_\mathcal{F}}
\\
W  \ar[r]^\tau_W  \ar[ur]^l
&
BO(d).
}
\]
In this way we describe the sheaf $\tau\mathcal{F}$ alternatively by 
\[
\tau\mathcal{F}(W)=\{l:W\rightarrow B\mathcal{F}(d)\mid p_\mathcal{F}\circ l =\tau_W\}.
\]
It is apparent that $\tau\mathcal{F}$ is indeed an equivariant sheaf.

There is a perhaps helpful rephrasing of this as follows.  Recall that the space of objects of the topological category $\mathsf{Emb}_d$ is a disjoint union of the spaces $Emb(W,\mathbb{R}^\infty)/{Diff}(W)$.  There is then a map
\[
ob\text{ }\mathsf{Emb}_d\xrightarrow{\tau} \coprod_{[W]}Map(W,BO(d)),
\]
given by $(W^\prime\subset\mathbb{R}^\infty)\mapsto \tau_{W^\prime}$, yielding the pull back square
\[
\xymatrix{
\tau^*(\coprod Map(W,B\mathcal{F}(d)) ) \ar[r]  \ar[d]^{pr}
&
\coprod Map(W,B\mathcal{F}(d))  \ar[d]^{p_\mathcal{F}}
\\
ob\text{ }\mathsf{Emb}_d   \ar[r]^\tau
&
\coprod Map(W,BO(d))
}
\]
with the vertical maps fibrations.
The equivariant sheaf $\tau\mathcal{F}$ is then given by $\tau\mathcal{F}(W^\prime)=pr^{-1}(W^\prime)$.

Refer to equivariant sheaves obtained in this way from a fibration $B_\theta\xrightarrow{\theta}BO(d)$ as \textit{fiberwise}.  The assignment $\mathcal{F}\mapsto\tau\mathcal{F}$ is then a prescription for how to obtain from a general equivariant sheaf a fiberwise equivariant sheaf.  \textit{Fiberwise} equivariant sheaves are much more homotopy theoretic.

\subsubsection*{The h-principle comparison morphism $\mathcal{F}\rightarrow \tau\mathcal{F}$}

Each $W\in \mathsf{Emb}_{d,N}$ is equipped with a metric coming from the ambient $\mathbb{R}^N$ and there is thus a preferred choice of exponential map $exp:V\rightarrow W$ defined on some neighborhood $V\subset TW$ of the zero-section.  This induces a map 
\[
\mathcal{F}(W)\xrightarrow{exp^*} \tau\mathcal{F}(W)
\]
defined as follows.  For each $p\in W$, there is  an embedding $exp_p:V_p\rightarrow W$.  For $g\in\mathcal{F}(W)$,  we obtain $exp_p^*g\in \mathcal{F}(V_p)$ restricting to $(exp_p^*g)_{\mid_0}\in Stalk_0(\mathcal{F}_{\mid_{\mathbb{R}^d}})$.  In this way, given $g\in\mathcal{F}(W)$ we can produced a section  of $Stalk_0(\mathcal{F}(\tau_W))$ by $p\mapsto (exp_p^*g)_{\mid_0}$.  This procedure describes, for each $W\in\mathsf{Emb}_d$,  a morphism of traditional continuous sheaves on $W$, 
\[
\mathcal{F}_{\mid_W} \rightarrow  \tau\mathcal{F}_{\mid_W}.  
\]
This procedure is sufficiently natural in $M$ to establish a canonical morphism of equivariant sheaves
\[
\mathcal{F}\rightarrow \tau\mathcal{F}.
\]

\begin{remark}
The morphism $\mathcal{F}\rightarrow\tau\mathcal{F}$ is a general way of assigning to an $\mathcal{F}$-structure on a manifold a fiberwise $\mathcal{F}$-structure on a manifold.  This process is familiar in some cases.  For instance, when $\mathcal{F}$ is the sheaf of complex structures, $\tau\mathcal{F}$ is the sheaf of almost-complex structures.   The sheaf $\tau\mathcal{F}$ is an equivariant version of Gromov's traditional sheaf $\mathcal{F}_{\mid_M}^*$ and the morphism $\mathcal{F}\rightarrow\tau\mathcal{F}$ is the equivariant version of the h-principle comparison morphism; that is, it is a weak homotopy equivalence exactly when the parametrized h-principle holds for the sheaf $\mathcal{F}$ (see~\cite{gromov:h-principle} page76).
\end{remark}

The main theorem of this paper can then be phrased as roughly saying the h-principle always holds on the `group completion' of the category of manifolds; or equivalently, each \textit{stable} family of $\tau \mathcal{F}$-manifolds is concordant to one which is simultaneously integrable.  This statement will be made more sensible as the paper progresses.

\begin{remark}
Many examples of equivariant sheaves come from geometric structures which are solutions to differential relations; more generally, geometric structures that are defined locally.  This list of examples is quite large.
\end{remark}

\subsection{Examples of equivariant sheaves}\label{subsec: examples1}

Here we outline some examples of outside interest.  All are equivariant sheaves with values in $\mathsf{Top}$.  In \textsection\ref{subsec: interpretations} - \ref{subsec: configurations} of this paper we will describe some of the examples in much greater detail.

\begin{enumerate}

\item  Take $\mathcal{F}$ to be the sheaf of continuous maps into a fixed target space $Y$.  That is, $\mathcal{F}(W)= Map(W,Y)$.  When $F$ is a smooth manifold and $Y=B{Diff}(F)$, $\mathcal{F}$ becomes the sheaf of bundles with fiber $F$.  

More generally, for $B_\theta\xrightarrow{\theta} BO(d)$ a fibration, take $\mathcal{F}$ to be given by 
\[
\mathcal{F}(W)=\{W\xrightarrow{l} B_\theta\mid \theta\circ l = \tau_W\}.
\] 
Many such examples come from a representation $G\xrightarrow{\rho} O(d)$.   In this event, take $B_\theta= EG\times_\rho Fr_d\xrightarrow{B\rho} BO(d)$ where $Fr_d$ is the infinite Steifel space of orthogonal $d$-frames in $\mathbb{R}^\infty$.  

For $G=SO(d)$, $\mathcal{F}$ is the sheaf of orientations.  For $G=U(d/2)$, $\mathcal{F}$ is the sheaf of almost-complex structures.  

Other examples of this type include the sheaf of $k$-dimensional distributions ($k\leq d$) as seen when $B_\theta$ is the tautological space over $BO(d)$ with fiber over $V^d\in BO(d)$ the  Grassmann $Gr_k(V)$ of $k$-planes in $V$.

\item  Take $\mathcal{F}$ to be the sheaf of immersions into a fixed manifold $Y^n$.  So $\mathcal{F}(W)=Imm(W,Y)$.  Similarly, take $\mathcal{F}$ to be the sheaf of submersions onto $Y$.  More generally, any (not necessarily open) $Diff$-invariant subspace of the jet space $J^k(\mathbb{R}^d,\mathbb{R}^n)$ results in a subspace of the jet space $J^k(W,Y)$.  One could take $\mathcal{F}(W)$ to be the sheaf of holonomic sections landing in this subspace of $J^k(W,Y)$.  In particular, one could prescribe a set $\mathcal{S}$ of allowed singularities and take $\mathcal{F}$ to be the sheaf of $\mathcal{S}$-singular maps to $Y$.

\item  Take $\mathcal{F}$ to be the sheaf of Riemannian metrics, or Riemannian metrics and isometric maps into a fixed Riemannian manifold $Y$.

\item Take $\mathcal{F}$ is the sheaf of complex structures.   Even more, take $\mathcal{F}$ to be the sheaf of complex structures and holomorphic maps into a fixed complex manifold $Y$.  For $d=2$ this amounts to the data of a holomorphic curve.

\item  Similarly, take $\mathcal{F}$ to be the sheaf of symplectic structure, or perhaps the sheaf of contact or Hamiltonian structures.  This has relevance to Symplectic Field Theory and, for $d=4$, to contact geometry.

\item  Take $\mathcal{F}$ to be the sheaf of Lagrangian immersions into a fixed $2d$-dimensional symplectic manifold $Y$.

\item  Take $\mathcal{F}$ to be the sheaf of $k$-dimensional foliations, maybe regarded as integrable $k$-dimensional distributions.

\item   Let $\xi$ be a vector field on the manifold $Y$, for example the gradient vector field of a Morse function on $Y$.  For $d=1$, take $\mathcal{F}$ to be the sheaf of maps into $Y$ which are flow lines of $\xi$.  

\item  Take $\mathcal{F}$ to be the sheaf of connections, or perhaps flat connections.  Of interest to Donaldson theory is when $d=4$ and $\mathcal{F}$ is the sheaf of anti-self-dual connections.  Toward Seiberg-Witten theory, take $\mathcal{F}$ to be the sheaf of Seiberg-Witten structures.

\item  Take $\mathcal{F}$ to be the the sheaf of (unordered) configurations.  Even stronger, for $(Y,Z)$ a pair of spaces, take $\mathcal{F}$ to be the sheaf of unordered configurations along with a map into $Y$ sending the configurations into $Z$.  For $Y$ complex one can further require holomorphicity of the maps into $Y$.  For $d=2$, this has relevance to the moduli space of (unordered) marked holomorphic curves in $Y$ with markings lying in $Z$.  Similar still, one could take $\mathcal{F}$ to the the sheaf of embedded $k$-manifolds, topologized appropriately, along with a map into $Y$ sending the $k$-manifolds into $Z$.  

\end{enumerate}

\section{The Geometric cobordism category}\label{sec: the geometric cobordism category}

We start off this section with some necessary constructions followed by a precise definition of the geometric cobordism category.  Located at the end of this section is a precise statement of the main theorem.

\subsection{Extending equivariant sheaves}\label{subsec: extending equivariant sheaves}

For $l< d$ there is an inclusion functor 
\[
\mathsf{Emb}_{l,N}\rightarrow \mathsf{Emb}_{d,N}
\]
given by $M\mapsto\mathbb{R}^{d-l-1}\times\mathbb{R}\times M$.  Think of this as $(d-l)$-fold suspension.  Extend $\mathcal{F}$ to be defined on $\mathsf{Emb}_{l,N}$, $l< d$, by 
\[
\mathcal{F}(M)=colim_{\epsilon>0} \mathcal{F}( \mathbb{R}^{d-l-1}\times (-\epsilon,\epsilon)\times M).
\]

For $W^d\subset \mathbb{R}^N$ an embedded $d$-manifold, write $pr_W:W\rightarrow \mathbb{R}^k$ for the projection onto the first $k$ coordinates and $pr_k:W\rightarrow \mathbb{R}$ for the projection onto the $k$th coordinate.  On an embedded $d$-manifold  $W^d\subset \mathbb{R}^{k-1}\times[a_0,a_1]\times\mathbb{R}^{N-k}$ such that for  $\partial_\nu W = pr^{-1}_k(a_\nu)$ ($\nu=0,1$) then
\[
\partial W = \partial_0 W\amalg \partial_1 W, 
\]
a \textit{collaring} of $W$ is an equivalence class of the data $(\epsilon>0, W_\epsilon)$ defined as follows.   $W_\epsilon$ is an embedded $d$-manifold 
\[
W_\epsilon\subset \mathbb{R}^{k-1}\times (a_0-\epsilon, a_1+\epsilon)\times\mathbb{R}^{N-k}
\]
satisfying    
\begin{itemize}

\item  $W = pr^{-1}_k([a_0,a_1])$ 

\item  $pr^{-1}_k(a_\nu - \epsilon , a_\nu + \epsilon)  	= (a_\nu - \epsilon , a_\nu + \epsilon) \times \partial_\nu W$ ($\nu=0,1$).  

\end{itemize}
Declare $(\epsilon,W_\epsilon)$ to be equivalent to $(\epsilon^\prime,W^\prime_{\epsilon^\prime})$ if there is a $\delta>0$ such that $pr^{-1}_k(a_0 - \delta, a_1 + \delta) = pr^{-1}_k(a_0-\delta, a_1+\delta)$. 
A \textit{collared embedding} $W\xrightarrow{e}W^\prime$ between collared manifolds is an $\epsilon$-germ of the data of embeddings $W_\epsilon\xrightarrow{e_\epsilon} W^\prime_\epsilon$ which restrict to a product embedding 
\[
(a_\nu - \epsilon, a_\nu +\epsilon)\times \partial_\nu W \xrightarrow{id\times e_\nu}(a_\nu - \epsilon, a_\nu +\epsilon)\times \partial_\nu W^\prime\text{ } (\nu=0,1).
\]
Extend $\mathcal{F}$ to collared-embedded manifolds $W^d\subset \mathbb{R}^{k-1}\times [a_0,a_1]\times\mathbb{R}^{N-k}$ by 
\[
\mathcal{F}(W):=colim_{\epsilon>0} \mathcal{F}(W_\epsilon).
\]

\subsection{Definition of the geometric cobordism category}\label{subsec: definition of cob}

\subsubsection*{Point-set description}\label{subsubsec: point-set description}

Maintain that the letter $\mathcal{F}$ denotes a fixed equivariant sheaf.  Construct the cobordism category with the extra structure $\mathcal{F}$, written $\mathsf{Cob}_{d,N}^\mathcal{F}$, as follows.  Declare
\[
ob\text{ }\mathsf{Cob}_{d,N}^\mathcal{F}=\{(a,M^{d-1},g)\mid g\in \mathcal{F}(M)\}
\]
where $a\in\mathbb{R}$ and $M^{d-1}\subset\mathbb{R}^{N-1}$ is a closed embedded $(d-1)$-manifold.  Declare
\[
mor\text{ }\mathsf{Cob}_{d,N}^\mathcal{F} = \{identities\}\amalg\{(W^d,g)\mid g\in \mathcal{F}(W)\}
\]
where $W\subset [a_0,a_1]\times \mathbb{R}^{N-1}$ is a \textit{compact} collared-embedded $d$-manifold.

Write $pr_1:W\rightarrow \mathbb{R}$ for projection onto the first coordinate.  The source and target maps of $\mathsf{Cob}_{d,N}^\mathcal{F}$ are given by 
\[
s,t:W\subset [a_0,a_1]\times \mathbb{R}^{N-1}\mapsto pr^{-1}_1(\{a_{0,1}\})
\] 
along with restriction of the $\mathcal{F}$-structures via the embeddings 
\[
(-\epsilon,\epsilon)\times s(W)  \hookrightarrow W_\epsilon \hookleftarrow (-\epsilon,\epsilon)\times t(W)
\]
for small enough $\epsilon>0$.
Composition is given by disjoint union in $\mathbb{R}^N$ along with  the gluing property of the equivariant sheaf $\mathcal{F}$.

Think of the role of the equivariant sheaf $\mathcal{F}$ in the definition of $\mathsf{Cob}_{d,N}^\mathcal{F}$ as prescribing \textit{geometric structure}.  When $\mathcal{F}$ is a fiberwise equivariant sheaf (i.e., comes from a fibration over $BO(d)$ as in \textsection\ref{subsec: fiberwise sheaves}), this geometric structure is sometimes referred to as \textit{tangential}; the phrase coming from the observation that in this case the geometric structure on a manifold is fiberwise from the tangent bundle of the manifold.

\subsubsection*{Quasi-topology of the geometric cobordism category}\label{subsubsec: quasi-topology of cob}

For $P$ and $Q$ possibly collared submanifolds of $\mathbb{R}^N$, the $C^\infty$ Whitney topology makes the diffeomorphism group ${Diff}(P)$ into a topological group.  The $C^\infty$ Whitney topology also makes the set $Emb(P,Q)$ of embeddings (collared embeddings when relevant) into a topological space having a free continuous action of ${Diff}(P)$ and ${Diff}(Q)$ by slices given by pre- and post-composition.

Fix $P = M^{d-1}$ closed or $P = W^d$ collared-embedded, without specifying.
For $a_0\leq a_1$, define 
\[
Emb^\mathcal{F}(P,[a_0,a_1]\times\mathbb{R}^{N-1}):=\{(e,g)\mid e\in Emb(P,[a_0,a_1]\times\mathbb{R}^{N-1}),\text{ and }g\in\mathcal{F}(e(P))\}.
\]
There is a bijection 
\[
Emb^\mathcal{F}(P,[a_0,a_1]\times\mathbb{R}^{N-1})\leftrightarrow Emb(P,[a_0,a_1]\times\mathbb{R}^{N-1})\times \mathcal{F}(P)
\]
given by $(e,g)\mapsto (e,e^*g)$.  Endow $Emb^\mathcal{F}(P,[a_0,a_1]\times\mathbb{R}^{N-1}$ with the resulting quasi-topology and note that the (free) diagonal action of ${Diff}(P)$ on $E_N^\mathcal{F}(P)$ is continuous.  Endow 
\[
B^\mathcal{F}_N(P):=Emb^\mathcal{F}(P,[a_0,a_1]\times\mathbb{R}^{N-1})/{Diff}(P)
\]
with the resulting colimit quasi-topology and observe that the projection 
\[
B^\mathcal{F}_N(P)\rightarrow B_N(P)
\]
is continuous.  Explicitly, 
\[
B^\mathcal{F}_N(P)=\{(P^\prime,g)\mid P\cong P^\prime\subset[a_0,a_1]\times\mathbb{R}^{N-1}\text{ and }g\in\mathcal{F}(P^\prime)\}.  
\]
Define the Q-space 
\[
E^\mathcal{F}_N(P):= Emb^\mathcal{F}(P,[a_0,a_1]\times\mathbb{R}^{N-1})\times_{{Diff}(P)} P 
\]
and observe the fiber bundle
\[
E^\mathcal{F}_N(P)\rightarrow B^\mathcal{F}_N(P)
\]
with smooth fiber $P\subset\mathbb{R}^N$.  The inclusion $P\subset\mathbb{R}^N$ yields a canonical embedding 
\[
E^\mathcal{F}_N(P)\hookrightarrow B^\mathcal{F}_N(P)\times \mathbb{R}^N.
\]

\begin{remark}
When $N=\infty$ Whitney's embedding theorem tells us that the spaces $Emb(P,[a_0,a_1]\times\mathbb{R}^{N-1})$ have trivial homotopy groups.  It follows that $Emb(P,[a_0,a_1]\times\mathbb{R}^{N-1})$  and $B_\infty(P)$ are models for $E{Diff}(P)$ and $B{Diff}(M)$ respectively.  Thus $B^\mathcal{F}_\infty(P) = \mathcal{F}(P)//{Diff}(P)$ is the \textit{moduli space of $\mathcal{F}$-structures on $P$}.
\end{remark}

Clearly $B^\mathcal{F}_N(P)$, with this quasi-topology, classifies bundles of pairs $(P^\prime,g)$ such that $P^\prime\subset\mathbb{R}^N$ is an embedded manifold diffeomorphic to $P$ and $g\in\mathcal{F}(P^\prime)$.  In this way there is a bijection
\[
ob\text{ }\mathsf{Cob}_{d,N}^\mathcal{F} \cong \mathbb{R}^\delta \times ( \coprod_{[M^{d-1}]} B^\mathcal{F}_N (M))
\]
where the disjoint union runs over closed $M^{d-1}\subset\mathbb{R}^{N-1}$, one in each diffeomorphism class.  Here, $\mathbb{R}^\delta$ is the underlying set of $\mathbb{R}$ with the discrete topology.  Use this bijection to endow $ob\text{ }\mathsf{Cob}^\mathcal{F}_{d,N}$ with a quasi-topology.

Provide the set of morphisms of $\mathsf{Cob}_{d,N}^\mathcal{F}$ with a quasi-topology in a similar way.  In the end there is a similar bijection
\[
mor\text{ }\mathsf{Cob}_{d,N}^\mathcal{F} \cong ob\text{ }\mathsf{Cob}_{d,N}^\mathcal{F} \amalg (\mathbb{R}^2_+)^\delta \times (\coprod_{[W]}  B_N^\mathcal{F} (W))
\] 
where $(\mathbb{R}^2_+)^\delta$ is the underlying set of $\{(a_0 , a_1) \in \mathbb{R}^2 \mid a_0 < a_1\}$ with the discrete topology, and the disjoint union runs over compact collared-embedded $d$-manifolds $W$, one in each diffeomorphism class.

It is straight forward to see that the structure maps (source, target, composition, and identity) of the category $\mathsf{Cob}_{d,N}^\mathcal{F}$ are indeed continuous.  Therefore, $\mathsf{Cob}_{d,N}^\mathcal{F}$ is a quasi-topological category.  Note that as such, this category is natural in the argument $\mathcal{F}$.

\begin{remark}
Recall that $\mathsf{QTop}$ has limits and colimits and contains $\mathsf{Top}$ via a Yoneda embedding.  Thus, the classifying space $B\mathsf{Cob}^\mathcal{F}_{d,N}$ makes sense as a quasi-topological space.  From \textsection\ref{subsec: quasi-topological spaces}, there is then an associated (weak) homotopy type.  
\end{remark}

\begin{remark}
An alternative to describing $\mathsf{Cob}^\mathcal{F}_d$ as a topological category is to describe it as an $\infty$-category, or complete Segal space (\cite{rezk:segal-spaces}), as in Lurie's survey of topological field theories (\cite{hopkins-lurie:tfts}).  It should be pointed out that the techniques for doing so are not so different.  Indeed, one must still describe a (semi-) simplicial Q-space replacing the nerve $N_\bullet \mathsf{Cob}^\mathcal{F}_d$.  The quasi-topology of each set of $k$-simplicies amounts to that described in this paper (see \textsection\ref{subsec: the sheaf}).  
\end{remark}

\subsubsection*{Smooth families}\label{subsec: smooth families}

In \textsection\ref{subsubsec: quasi-topology of cob}
above, we used that $B^\mathcal{F}_N(P)$ classifies bundles of $\mathcal{F}$-manifolds with underlying manifold diffeomorphic to $P$.  We go further into this idea here.

Recall the fiber bundle
\[
E^\mathcal{F}_N(P)\rightarrow B^\mathcal{F}_N(P)
\]
with smooth fiber $P^d\subset\mathbb [a_0,a_1]\times\mathbb{R}^{N-1}$.  For $X^q$ a smooth manifold, we would like the notion of a `smooth' map $X\rightarrow B^\mathcal{F}_N(P^d)$ to be exactly so that the pull back 
\[
E^{d+q}\rightarrow X
\]
is a \textit{smooth} map between \textit{manifolds}.  Precisely, 
given a map $X\xrightarrow{f} B_N(P^d)$, write $f(x)=(P_x,g_x)$.  Call $f$ a \textit{smooth map} (over $X$) if the subset
\[
E:=\{(x , p)\mid p\in P_x\}\subset X\times [a_0,a_1]\times\mathbb{R}^{N-1}
\]
is a smooth $(q+d)$-manifold.  Write $Map^{sm}$ for the space of \textit{smooth} maps.  Call two smooth maps $f_0,f_1:X\rightarrow B_N(P)$ \textit{concordant} if there is a smooth map $H:\mathbb{R}\times X\rightarrow B_N(P)$ which is constant in a neighborhood of $((-\infty,0]\amalg [1,\infty))\times X$.

Rephrase this notion of smooth as follows.   Say a map $X\xrightarrow{f} B_N(P)$ is \textit{smooth} if the transition functions of the resulting fiber bundle $E\rightarrow X$, written as maps $U\rightarrow {Diff}(P)$ have \textit{smooth} adjoints $U\times P\rightarrow P$.  Phrased as such, because smooth maps approximate continuous maps, the inclusion 
\[
Map^{sm}(X,B_N(P))\hookrightarrow Map(X,B_N(P))
\]
is a weak homotopy equivalence.  Indeed, any map $S^k\rightarrow Map(X,B_N(P))$ is approximated by, and concordant to, a \textit{smooth} map $S^k\times X\rightarrow B_N(P)$ and the same holds for $\mathbb{R}\times S^k$ in place of $S^k$.

\subsection{Comparison of geometric to tangential cobordism categories}\label{subsubsec: comparison}

\subsubsection*{Statement of the main theorem}

The h-principle comparison morphism $\mathcal{F}\rightarrow\tau\mathcal{F}$ results in a functor of geometric cobordism categories
\[
\mathsf{Cob}_{d,N}^\mathcal{F}\rightarrow \mathsf{Cob}_{d,N}^{\tau\mathcal{F}}.
\]
We show on the level of classifying spaces that this functor is a weak homotopy equivalence.  This  is done by computing their weak homotopy types directly and observing that they are identical.

Write $Th$ for the Thom space functor and $\gamma_{d.N}^\perp$ for the tautological perpendicular bundle over $Gr_{d,N}$ of $(N-d)$-planes in $\mathbb{R}^N$.  Recall the fibration $Gr_{d,N}^\mathcal{F}\xrightarrow{p_\mathcal{F}} Gr_{d,N}$ from the beginning of \textsection\ref{subsubsec: alternative description}

\begin{theorem*}[Main Theorem]
The h-principle comparison morphism $\mathcal{F}\rightarrow\tau\mathcal{F}$ induces a weak homotopy equivalence
\[
B\mathsf{Cob}_{d,N}^\mathcal{F}\xrightarrow{\simeq} B\mathsf{Cob}_{d,N}^{\tau\mathcal{F}}
\]
with weak homotopy type 
\[
\Omega^{N-1}Th(p_\mathcal{F}^*\gamma_{d,N}^\perp).
\]
\end{theorem*}

\begin{remark}
It should be noted that although the statement of the main theorem in the introduction is indeed true, what is proved here and is more prepared for applications requires the following modification.  Consider the stalk $Stalk_0(\mathcal{F}_{\mid_{\mathbb{R}^d}})$ of the restriction $\mathcal{F}_{\mid_{\mathbb{R}^d}}$ at $0\in\mathbb{R}^d$.  There is an action of $O(d)$ on $Stalk_0(\mathcal{F}_{\mid_{\mathbb{R}^d}})$.  Redefine the map $p_\mathcal{F}$ as the projection
\[
p_\mathcal{F}:Fr_{d,N}\times_{O(d)} Stalk_0(\mathcal{F}_{\mid_{\mathbb{R}^d}}) \longrightarrow Gr_{d,N}.
\]
This has the benefit of throwing away all information of $\mathcal{F}(\mathbb{R}^d)$ away from the origin and is thus less data to keep track of when applying the main theorem.  The obvious drawback of such a reformulation is the use of the notion of a \textit{stalk} of the $\mathsf{Top}$-valued sheaf.  
\end{remark}


\subsubsection*{Statement of the main theorem for $N=\infty$}

The Thom space functor $Th: \textbf{Vector Bundles}\rightarrow \mathsf{Top}$ extends to virtual vector bundles upon taking values in the category of spectra.  That is, there is a Thom spectrum functor $Th:\mathsf{Top}\text{ } \mathsf{over}\text{ } \mathsf{BO}\rightarrow \mathsf{Spectra}$.  Likewise in the category $\mathsf{QTop}$ there is the Thom (quasi-)spectra functor, $Th:\mathsf{QTop}\text{ } \mathsf{over}\text{ } \mathsf{BO}\rightarrow\mathsf{QSpecra}$.  The space $BO$ is an $H$-space having an inverse $-1:BO\rightarrow BO$ which is seen on the level of finite Grassmann manifolds as $Perp: Gr_{d,N}\rightarrow Gr_{N-d,N}$ given by $V\mapsto V^\perp$.  Write $-\theta$ for the composition 
\[
B_\theta\xrightarrow{\theta} BO(d)\hookrightarrow BO\xrightarrow{-1} BO.
\]
Explicitly, $Th(-\theta)$ is the quasi-spectrum whose $N^{th}$ Q-space is $ Th(\theta^*\gamma^\perp_{d,N})$.

Apply this to $B\mathcal{F}(d)\xrightarrow{p_\mathcal{F}}BO(d)$.
An alternative way to write $B\mathcal{F}(d)$ is as $Stalk_0(\mathcal{F}_{\mid_{\mathbb{R}^d}})//O(d)$ (see \textsection\ref{subsubsec: alternative description}) which is weakly homotopy equivalent to $\mathcal{F}(\mathbb{R}^d)//O(d)$ via Proposition~\ref{prop: stalk} applied to $B=*$.  For notation, write
\[
MT\mathcal{F}(d):= B\mathcal{F}(d)^{-p^*_\mathcal{F}\gamma_d}.
\]

\begin{theorem*}[Main Theorem ($N=\infty$)]

The h-principle comparison morphism $\mathcal{F}\rightarrow \tau\mathcal{F}$ induces a weak homotopy equivalence
\[
B\mathsf{Cob}_d^\mathcal{F}\xrightarrow{\simeq} B \mathsf{Cob}_d^{\tau\mathcal{F}}
\]
with weak homotopy type 
\[
\Omega^{\infty-1}MT\mathcal{F}(d)  
\simeq \Omega^{\infty-1}(\mathcal{F}(\mathbb{R}^d)//O(d))^{-\gamma_d}.
\]
\end{theorem*}

\subsection{Topological manifolds}

The proof of the main theorem to follow does not heavily use the smoothness of the manifolds at hand.  In fact, the author expects that one could consider a cobordism category of topological manifolds.  The notion of embeddings would be replaced by locally \textit{flat} embeddings (\cite{kirby-siebenmann:topological-manifolds}).  A flat embedding of a topological manifold $M^m\xrightarrow{e} N^n$ is an embedding of topological manifolds with the additional property that for each $p\in M$ there are charts $p\in U_\alpha\subset M$ and $e(p)\in V_\beta\subset N$ with the restriction $e:\mathbb{R}^m\cong U_\alpha \rightarrow V_\beta\cong \mathbb{R}^n$ given as an inclusion of vector spaces.  The role of the tautological virtual bundle $-\gamma_d$ would be replaced by the fiberwise Spanier-Whithead dual of the  tautological $\mathbb{R}^d$-bundle also denoted $\gamma_d$ which lies over the `Grassmann' $BTop(d):=Emb^{flat}_0(\mathbb{R}^d,\mathbb{R}^\infty)/Top(d)$.  There is a resulting Thom spectrum $BTop(d)^{-\gamma_d}$.  One could do similarly for PL-manifolds.

\begin{conj}
Let $CAT$ be either PL or Top.  There is a cobordism category $\mathsf{Cob}^{CAT}_d$ of $d$-dimensional $CAT$-manifolds.  This category is a topological category with 
\[
B\mathsf{Cob}^{CAT}_d \simeq \Omega^{\infty-1} BCAT(d)^{-\gamma_d}.
\]
\end{conj}

An outline of a proof of this conjecture could follow from
\begin{enumerate}

\item a tubular neighborhood theorem for flat-embedded topological manifolds,

\item a Whitney embedding theorem for flat embeddings of topological manifolds,

\item transversality for the projection maps $W \rightarrow \mathbb{R}^k$ of flat-embedded topological manifolds $W$.  
\end{enumerate}

\section{Proof of the main theorem}\label{sec: proof of main theorem}

\subsection{Strategy}\label{subsec: strategy}

To prove the main theorem we use ideas from~\cite{galatius:graphs} and~\cite{galatius-madsen-tillmann-weiss:cobordism}.  The strategy is to break the desired weak homotopy equivalence 
\[
B\mathsf{Cob}^\mathcal{F}_{d,N}\simeq \Omega^{N-1}Th(p^*_\mathcal{F}\gamma^\perp_{d,N})
\]
into six stages (we define the relevant objects used below in the subsections to follow; all arrows will be shown to be weak homotopy equivalences):
\[
B\mathsf{Cob}^\mathcal{F}_{d,(k,N)}\xleftarrow{\alpha} B\mathsf{D}^{\mathcal{F},\perp}_{d,(k,N)}\xrightarrow{\beta} B\mathsf{D}^{\mathcal{F},\pitchfork}_{d,(k,N)}\xrightarrow{\gamma}D^\mathcal{F}_{d,(k,N)}\xrightarrow{\delta} \Omega B\mathsf{Cob}^\mathcal{F}_{d,(k-1,N)},
\]
\[
\mathsf{Cob}^\mathcal{F}_{d,N}\cong \mathsf{Cob}^\mathcal{F}_{d,(1,N)}, \hspace{.5cm} \text{and}  \hspace{.5cm}  Th(p^*_\mathcal{F}\gamma^\perp_{d,N})\xrightarrow{\simeq}D^\mathcal{F}_{d,(0,N)}.
\]

The equivalences $\alpha$, $\beta$, and the congruence in the bottom left have little content and are there for book keeping only.  The equivalence $\gamma$ amounts to easy transversality arguments.  Showing the map $\delta$ is an equivalence is a group completion argument in the sense of~\cite{mcduff-segal:group-completion} in the case of the homology of a topological monoid, \cite{dwyer:centralizer} (3.12)  in the case of the homotopy groups of a discrete category, and later of~\cite{tillmann:mapping-class-group} in the case of the homology of a topological category.  The equivalence in the bottom right can be thought of as a sort of scanning map in the sense of~\cite{segal:scanning}.

\subsection{The sheaf $\Psi^\mathcal{F}_d$ of $\mathcal{F}$-manifolds}\label{subsec: the sheaf}

\subsubsection*{$\Psi^\mathcal{F}_d$ as a $\mathsf{Set}$-Valued Sheaf}

Consider the functor $\Psi_d: \mathsf{Emb}^{op}_N\rightarrow \mathsf{Set}$ described by 
\[
\Psi_d(Q)=\{W^d\subset Q\}
\]
where $W$ above is required to be a $d$-dimensional  boundaryless submanifold of $Q$, closed as a subset of $Q$.
Because $e^{-1}(W)\in \Psi_d(Q^\prime)$ for each $e\in Emb(Q^\prime, Q)$, this is indeed a functor.
It is straight forward to verify in fact that ${\Psi_d}_{\mid_Q}$ is a traditional ($\mathsf{Set}$-valued) sheaf for each $Q^N\in\mathsf{Emb}_N$.  It amounts to the defining feature that $d$-manifolds are locally defined.

Define similarly the functor $\Psi^\mathcal{F}_d:\mathsf{Emb}_N^{op}\rightarrow \mathsf{Set}$ described by 
\[
\Psi^\mathcal{F}_d(Q):=\{(W,g)\mid W\in\Psi_d(Q),\text{ and } g\in\mathcal{F}(W)\}.
\]
This is indeed a functor since for $g\in\mathcal{F}(W)$ and $e\in Emb(Q^\prime,Q)$,  $e^*g\in\mathcal{F}(e^{-1}(W))$.
The sheaf-gluing property of $\mathcal{F}$ along with the already established sheaf-property of $\Psi_d$ makes ${\Psi^\mathcal{F}_d}_{\mid_Q}$ into a traditional ($\mathsf{Set}$-valued) sheaf on $Q$ for each $Q\in \mathsf{Emb}_N$.  Note the forgetful morphism $\Psi^\mathcal{F}_d\rightarrow \Psi_d$.

\subsubsection*{$\Psi^\mathcal{F}_d$ as an Equivariant Sheaf}

It will now be shown that $\Psi^\mathcal{F}_d$ is an equivariant sheaf.  For this, declare a point-set map $X\xrightarrow{f} \Psi^\mathcal{F}_d(Q)$ to be ``continuous" if the following holds.  Write $f(x)=(W_x,g_x)$.  For each $x\in X$ and each $q\in Q$ there are open sets $x\in U_x\subset X$ and $q\in V^\prime_q\subset V_q\subset Q$, with $\overline{V^\prime_q}\subset V_q$ compact, and a continuous map
\[
U_x\xrightarrow{(j,h)} Emb(V_q,Q)\times \mathcal{F}(V_q\cap W_x)
\]
with $x\mapsto (\{V_q\subset Q\},\iota^*_{V_q} g_x)$,
such that 
\[
\iota^*_{\overline{V}^\prime_q}f_{\mid_{U_x}}(x^\prime) =\iota^*_{\overline{V}^\prime_q}(j_{x^\prime}^{-1}(W_x),h_{x^\prime}) \in \Psi^\mathcal{F}_d(\overline{V}^\prime_q).
\]

Letting $Map(X,\Psi^\mathcal{F}_d(Q))$ denote the set of such ``continuous" maps, it is evident by inspection that it is natural with respect to continuous maps among the argument $X$; thus $Map(-,\Psi^\mathcal{F}_d(Q))$ is a functor $\mathsf{Top}^{op}\rightarrow \mathsf{Set}$.  But more, from its definition, ``continuity" of such a map is a local condition, so the necessary sheaf condition is satisfied for $\Psi^\mathcal{F}_d(Q)$ to be a quasi-topological space.  Moreover, this quasi-topology is natural with respect to embeddings $Q^\prime\rightarrow Q$.  Also,
it is apparent from this quasi-topology that the restriction map $Emb(P,Q)\times \Psi^\mathcal{F}_d(Q)\rightarrow \Psi^\mathcal{F}_d(P)$ is ``continuous''.  In this way, $\Psi^\mathcal{F}_d$ takes values in $\mathsf{QSpace}$ and becomes an equivariant sheaf.  Note that there was no dependence of $\Psi^\mathcal{F}_d$ on the dimension $N$ of its argument.  So in fact 
\[
\Psi^\mathcal{F}_d:\mathsf{Emb}_N^{op}\rightarrow \mathsf{QTop}
\]
is an equivariant sheaf for any $N$.  As such, it is clear that $\Psi^\mathcal{F}_d$ is natural with respect to morphisms of equivariant sheaves $\mathcal{F}^\prime\rightarrow\mathcal{F}$.  In particular, the forgetful functor $\Psi^\mathcal{F}_d\rightarrow\Psi_d$ is a morphism of equivariant sheaves.

\begin{prop}
If the equivariant sheaf $\mathcal{F}$ is represented, then the equivariant sheaf $\Psi^\mathcal{F}_d$ is also represented.  In particular, $\Psi_d$ is represented.

\end{prop}

\begin{proof}

Recall that to say an equivariant sheaf is represented is to say it takes values in $\mathsf{Top}$ rather than $\mathsf{QTop}$.

Let $V^\prime\subset V\subset Q$ be open sets with $\overline{V}^\prime \subset V$ compact.  Let $U\subset Emb(V, Q)$ be a neighborhood of the inclusion $V\subset Q$ such that for each $e\in U$, $\overline{V}^\prime\subset e(V)$.  Fix $(W,g)\in\Psi^\mathcal{F}_d(Q)$ and let $g\in O_g\subset \mathcal{F}( V\cap W)$ be an open set.  Define 
\[
C(W,V^\prime,V,U,O_g):=\{(W^\prime,g^\prime)\in\Psi^\mathcal{F}_d(Q)\mid \iota^*_{\overline{V^\prime}} (W^\prime,g^\prime) \in \iota^*_{\overline{V^\prime}} (U^* W, O_g)\}
\]
where the last condition holds in $\Psi^\mathcal{F}_d(\overline{V}^\prime)$ and 
\[
U^*W=\{W^\prime \mid \exists j\in U\text{ with } j^{-1} (W^\prime) = W\}.   
\]
The collection $\{C(W,V^\prime, V, U,O_g)\}$ forms a basis for a topology on $\Psi^\mathcal{F}_d(Q)$ as follows.

Let $(W,g)\in C(W_0,V_0^\prime,V_0,U_0,O_{g_0})\cap C(W_1,V_1^\prime,V_1,U_1,O_{g_1})$.  Take $(V^\prime,V) = (V^\prime_0 \cup V^\prime_1, V_0 \cup V_1)$.  Let $\iota^*_V\mathcal{F}(W)\xrightarrow{r} \iota^*_{V_0}\mathcal{F}(W) \times \iota^*_{V_1}\mathcal{F}(W)$ be restriction.  The continuity of $r$ ensures the existence of an open set $g\in O_g \subset r^{-1}(O_{g_0}\times O_{g_1})$.  Similarly, the restriction map $Emb(V,Q)\xrightarrow{R} Emb(V_0,Q)\times Emb(V_1,Q)$ is continuous thus ensuring the existence of a neighborhood $\{V\subset Q\}\in U\subset R^{-1}(U_0\times U_1)$.  It is clear then that 
\[
(W,g)\in C(W,V^\prime,V,U,O_g)\subset C(W_0,V_0^\prime,V_0,U_0,O_{g_0})\cap C(W_1,V_1^\prime,V_1,U_1,O_{g_1}).
\]

It is almost by definition that a point-set map $X\rightarrow \Psi^\mathcal{F}_d(Q)$ is ``continuous" exactly if it is continuous with respect to the topology above.  The rest of the axioms for $\Psi^\mathcal{F}_d$ to be representable follow without trouble.
\end{proof}

\begin{remark}
There are some unusual phenomena with this topology on $\Psi_d$.  For example, a neighborhood of $\emptyset\in \Psi_d(\mathbb{R}^N)$ consists of $W\in\Psi_d(\mathbb{R}^N)$ such that $W\cap K =\emptyset$ for some $K\subset\mathbb{R}^N$ compact.  So for $W$ a compact $d$-manifold, the  assignment $t\mapsto W+1/t$ for $t>0$ and $0\mapsto \emptyset$ is a continuous path $[0,1]\rightarrow \Psi_d(\mathbb{R}^N)$.  In particular, it is \textit{not} true that each member of a continuous family $X\rightarrow\Psi_d(\mathbb{R}^N)$ of $d$-submanifolds of $\mathbb{R}^N$ has the same topological type.  In what follows, notice that this topology on $\Psi_d(\mathbb{R}^N)$ is exactly what is needed in order to `group complete' $\mathsf{Cob}_{d,N}$ as $\Omega^{N-1}\Psi_d$.  
\end{remark}

\subsection{Categories $D^{\mathcal{F},-}_{d,(k,N)}$ of `unbounded' $\mathcal{F}$-manifolds}\label{subsec: bounded manifolds}

\subsubsection*{Comparing $\mathsf{D}^{\mathcal{F},\pitchfork}_{d,(k,N)}$ and $D^\mathcal{F}_{d,(k,N)}$}

Define the quasi-topological subspace 
\[
D^\mathcal{F}_{d,(k,N)}:=\{(W,g)\in \Psi_d^\mathcal{F}(\mathbb{R}^N)\mid W\subset \mathbb{R}^{k}\times \mathring{I}^{N-k}\}
\]
consisting of embedded $\mathcal{F}$-manifolds which are `bounded' in $N-k$ chosen directions.  In the notation above, $\mathring{I}^{N-k}$ denotes the interior $int(I^{N-k})$.

For each $W\in D_{d,(k,N)}$, there is the projection onto the first $k$-coordinates $pr_W :W \rightarrow \mathbb{R}^k$.  Write $(\mathbb{R}^l)^\delta$ for $\mathbb{R}^l$ with the discrete topology.  Let $\mathsf{D}^{\mathcal{F},\pitchfork}_{d,(k,N)}$ be the Q-space 
\[
\{(a,W,g)\mid (W,g)\in {D_{d,(k,N)}^\mathcal{F}}\text{ and } a\in(\mathbb{R}^k)^\delta \text{ is transverse to }pr_W\}.
\]
For $a^0,a^1\in\mathbb{R}^{k}$, write $a^0_k \leq a^1_k$ to mean the $k$th coordinates satisfy $a^0_k\leq a^1_k$ and similarly for strict inequality.  This space is a quasi-topological partially ordered set via the rule $(a^0,W_0,g_0)\prec (a^1,W_1,g_1)$ exactly when $(W_0,g_0)=(W_1,g_1)$ and $a^0\leq a^1$.

Regard the Q-space $D^\mathcal{F}_{d,(k,N)}$ as a quasi-topological category with only identity morphisms.  There is a forgetful functor 
\[
\mathsf{D}^{\mathcal{F},\pitchfork}_{d,(k,N)} \xrightarrow{\gamma} D_{d,(k,N)}^\mathcal{F}
\]
given by $(a,W,g)\mapsto (W,g)$.

\begin{lem}[$\gamma$]\label{lem:gamma}
The functor $\gamma$ induces a weak homotopy equivalence 
\[
B\mathsf{D}^{\mathcal{F},\pitchfork}_{d,(k,N)}\xrightarrow{B\gamma}  D_{d,(k,N)}^\mathcal{F}.
\]

\end{lem}

To prove this Lemma~\ref{lem:gamma} we will use the following lemma from~\cite{galatius:graphs} which is a corollary of a theorem of Segal (\cite{segal:classifying-spaces}, A.1).  It is a sort of Quillen's Theorem A for topological categories.  Recall that a map is \textit{etale} if it is a local homeomorphism and an open map.  

\begin{lem}\label{lem:etale}
Let $X_\bullet$ be a simplicial space and $Y$ a space.  Regard $Y$ as a constant simplicial space.  Let $f_\bullet:X\bullet\rightarrow Y$ be a simplicial map such that each $f_k$ is etale.  If the simplicial set (discrete space) $f^{-1}(y)$ has contractible classifying space for each $y\in Y$, then $Bf:BX_\bullet \rightarrow Y$ is a weak homotopy equivalence.
\end{lem}

\begin{proof}[Proof of Lemma~\ref{lem:gamma}]

For each $l$, the map $N_l \mathsf{D}^\mathcal{F}_{d,(k,N)} \rightarrow D^\mathcal{F}_{d,(k,N)}$ is etale.  The fiber $(N_\bullet \gamma)^{-1}(W,g)$ is the nerve of the (discrete) poset $\mathsf{P}$ of points in $\mathbb{R}^k$ which are transverse to the projection $pr:W\rightarrow \mathbb{R}^k$.  This fiber is apparently non-empty and the partial ordering comes from $a^0\leq a^1$.  For each $a\in\mathsf{P}$, the sub-poset $\mathsf{C}_a := \{a^\prime \mid a^\prime_i = a_i \text{ for } i < k\}$ is linearly ordered and thus $B\mathsf{C}\simeq *$.  Moreover, for each $b\in \mathsf{P}$ there exists an $a^\prime\in \mathsf{C}_a$ such that either $b\leq a^\prime$ or $a^\prime \leq b$.  It follows that homotopy colimits over the discrete category $P$ agree up to weak homotopy equivalence with those over $\mathsf{C}_a$.  In particular $B\mathsf{P}\simeq B\mathsf{C}_a$ as homotopy colimits over the trivial functor.  Lemma~\ref{lem:etale} then applies.

\end{proof}

\subsubsection*{Comparing $\mathsf{D}^{\mathcal{F},\pitchfork}_{d,(k,N)}$ and $\mathsf{D}^{\mathcal{F},\perp}_{d,(k,N)}$}

We say $W\subset \mathbb{R}^{k-1}\times \mathbb{R}\times I^{N-k}$ is \textit{cylindrical} near $a\in\mathbb{R}^k$ if $\{a\}$ is transverse $pr_W:W\rightarrow \mathbb{R}^k$ and there exists $\epsilon > 0$ such that 
\[
pr_W^{-1}(\mathbb{R}^{k-1}\times (a_k-\epsilon,a_k+\epsilon)) = \mathbb{R}^{k-1}\times (a_k-\epsilon,a_k+\epsilon)\times pr_W^{-1}(a)\subset\mathbb{R}^N.
\]
Define $\mathsf{D}^{\mathcal{F},\perp}_{d,(k,N)}\subset \mathsf{D}^{\mathcal{F},\pitchfork}_{d,(k,N)}$ as the full quasi-topological partially ordered set of triples $(a,W,g)$ such that $W$ is \textit{cylindrical} near $a$.  Label this inclusion of quasi-topological partially ordered sets by
\[
\beta: \mathsf{D}^{\mathcal{F},\perp}_{d,(k,N)}\hookrightarrow \mathsf{D}^{\mathcal{F},\pitchfork}_{d,(k,N)}.
\]

\begin{lem}[$B\beta$]
The inclusion $\beta$ induces a homotopy equivalence 
\[
B\beta: B\mathsf{D}^{\mathcal{F},\perp}_{d,(k,N)}\xrightarrow{\simeq} B \mathsf{D}^{\mathcal{F},\pitchfork}_{d,(k,N)}.
\]
\end{lem}

\begin{proof}

We will show that $\beta$ induces a level-wise weak homotopy equivalence of simplicial nerves
\[
N_l\mathsf{D}^{\mathcal{F},\perp}_{d,(k,N)}\simeq N_l \mathsf{D}^{\mathcal{F},\pitchfork}_{d,(k,N)}.
\]
For each $\underline{a}=(a^0\leq ...\leq a^l)\in(\mathbb{R}^k)^{l+1}$ and $\delta:\mathbb{R}^k \rightarrow (0,\infty)$ sufficiently small, let $\lambda = \lambda^\delta_{\underline{a}} :\mathbb{R}^k\rightarrow \mathbb{R}^k$ be a smooth non-decreasing function which takes the constant value $a^i$ near $a^i$ and which is the identity outside a $\delta(a^i)$-neighborhood of each $a^i$ (see Figure~\ref{fig:lambda-1}).


For each $t\in \mathbb{R}$, consider the self-map $\phi_{t,\underline{a}}$ of $\mathbb{R}^k\times I^{N-k}$ given by 
\[
\phi_{t,\underline{a}}(x,v)=((1-s(t))x+s(t)\lambda_{\underline{a}}(x),v)
\]
where $s:\mathbb{R}\rightarrow [0,1]$ is a smooth monotonically increasing function with $s(t) = 0$ for $t\leq 0$ and $s(t) = 1$ for $t\geq 1$.
Let $(\underline{a},W) \in \mathsf{D}^\pitchfork_{d,(k,N)}$.  Although $\phi_{t,\underline{a}}$ is not an embedding, the preimage $\phi^{-1}_{t,\underline{a}}(W)\subset\mathbb{R}^k\times I^{N-k}$ is again a smooth $d$-submanifold.  In fact, $(\underline{a}, \phi^{-1}_{t,\underline{a}}(W))\in \mathsf{D}^\pitchfork_{d,(k,N)}$ for all $t\in \mathbb{R}$ with $(\underline{a},\phi^{-1}_{1,\underline{a}}(W))\in \mathsf{D}^\perp_{d,(k,N)}$ for $t\geq1$.  Define the self-map of $N_l \mathsf{D}^\pitchfork_{d,(k,N)}$ by 
\[
\Phi_t(\underline{a},W) = (\underline{a}, \phi^{-1}_{t,\underline{a}}(W)).
\]
It is routine to verify that $\Phi_t$ realizes a deformation retraction of $N_l \mathsf{D}^\pitchfork_{d,(k,N)}$ onto $N_l\mathsf{D}^\perp_{d,(k,N)}$.  We  must now do similarly with $\mathcal{F}$-structures.

\begin{figure}
\begin{center}
	\includegraphics[width=3.5in,viewport=0 400 600 792,clip=true]{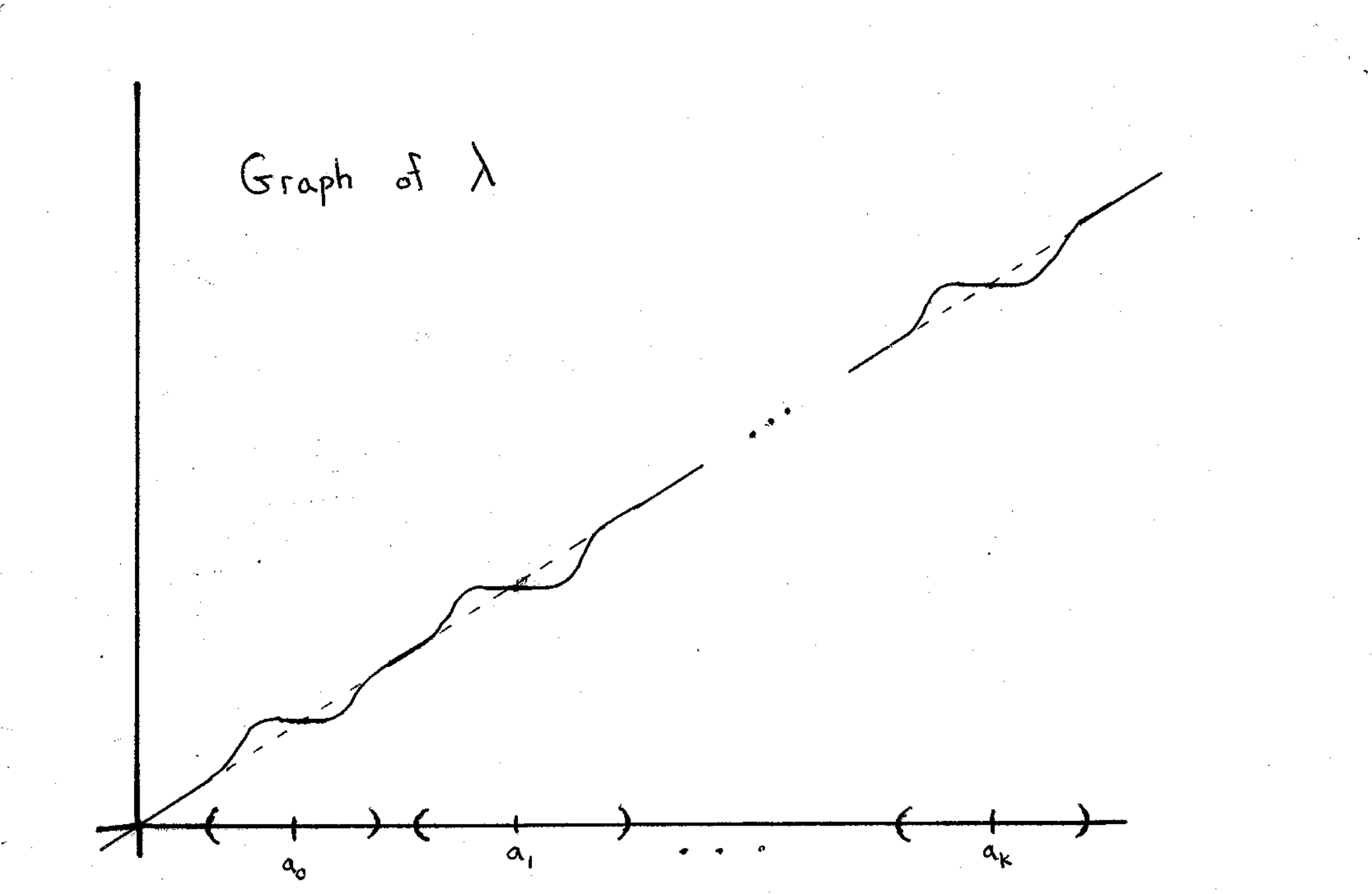}
	\caption{$\lambda$}
	\label{fig:lambda-1}
\end{center}
\end{figure}

Consider a smooth compact family $K\xrightarrow{f} N_l \mathsf{D}^{\mathcal{F},\pitchfork}_{d,(k,N)}$ written $f(x) = (\underline{a}_x,W_x,g_x)$.  There exists a smooth function $\epsilon:\mathbb{R}^k \rightarrow (0,\infty)$ such that for every $x\in K$ an $\epsilon$-neighborhood of $W_x\subset \mathbb{R}^k\times I^{N-k}$ is a tubular neighborhood of $W_x$.  We call this tubular neighborhood $N_\epsilon (W_x)$.  For such an $\epsilon$ there is a small enough $\delta:\mathbb{R}^k\rightarrow (0,\infty)$ as above, such that $\phi_{t,\underline{a}_x}^{-1}(W_x)\subset N_\epsilon (W_x)$ for all $t\in \mathbb{R}$.  Moreover, upon possibly shrinking $\epsilon$ and $\delta$ and choosing once and for all a nice enough function $\lambda$, it can be guaranteed that the projection $N_\epsilon (W_x) \rightarrow W_x$ restricts to a diffeomorphism 
\[
pr_{t,x}: \phi_{t,\underline{a}_x}^{-1}(W_x) \rightarrow W_x
\]
for each $t\in \mathbb{R}$.  This describes a smooth family $\mathbb{R}\times K\xrightarrow{\Phi} N_l\mathsf{D}^{\mathcal{F},\pitchfork}_{d,(k,N)}$ given by 
\[
\Phi_t(x) = (\underline{a}_x,\phi_{t,\underline{a}}^{-1}(W_x), {pr_{t,x}}^* g_x).
\]
Notice that $\Phi_t(x)\in N_l\mathsf{D}^{\mathcal{F},\perp}_{d,(k,N)}$ for $t\geq 1$.  It is not difficult to verify that the homotopy class of $\Phi_t$ is independent of the choices of functions $\epsilon$ and $\delta$.  Evaluating $\Phi_t$ at $t=1$ thus describes an assignment 
\[
[K,N_l\mathsf{D}^{\mathcal{F},\pitchfork}_{d,(k,N)}]\rightarrow [K,N_l\mathsf{D}^{\mathcal{F},\perp}_{d,(k,N)}].
\]
It is routine to verify that the inclusion $\mathsf{D}^{\mathcal{F},\perp}_{d,(k,N)}\rightarrow \mathsf{D}^{\mathcal{F},\pitchfork}_{d,(k,N)}$ demonstrates an inverse to the above assignment.  The result follows.

\end{proof}

\subsection{Some delooping cobordism categories $\mathsf{Cob}^\mathcal{F}_{d,(k,N)}$}\label{subsec:delooping-categories}

We now construct an artificial cobordism category $\mathsf{Cob}_{d,(k,N)}^\mathcal{F}$ motivated by showing 
\[
B\alpha: B\mathsf{D}^{\mathcal{F},\perp}_{d,(k,N)}\xrightarrow{\simeq} B\mathsf{Cob}_{d,(k,N)}^\mathcal{F}\text{\hspace{.5cm} and \hspace{.5cm}}   \delta: D^\mathcal{F}_{d,(k+1,N)}\xrightarrow{\simeq} \Omega B\mathsf{Cob}^\mathcal{F}_{d,(k,N)}.
\]

\subsubsection*{Comparing $\mathsf{D}^{\mathcal{F},\perp}_{d,(k,N)}$ and $\mathsf{Cob}^\mathcal{F}_{d,(k,N)}$}\label{subsubsec: comparing D and Cob}

Declare
\[
ob\text{ }\mathsf{Cob}^\mathcal{F}_{d,(k,N)}=\{(a,M^{d-k},g)\}
\]
where $ a\in\mathbb{R}^k$, $M\subset I^{N-k}$ is a closed $(d-k)$-manifold, and $ g\in \mathcal{F}(M)$. 
Similarly, declare  
\[
mor\text{ }\mathsf{Cob}^\mathcal{F}_{d,(k,N)}=\{identities\}\amalg\{(a^0,a^1,W^d,g)\}
\]
where $a^0 < a^1\in \mathbb{R}^k$ and $W\subset \mathbb{R}^{k-1}\times [a^0_k,a^1_k]\times int(I^{N-k})\subset \mathbb{R}^{k-1}\times [a^0_k,a^1_k]\times \mathbb{R}^{N-k}$ is a collared embedded $d$-manifold, closed as a subset, and a collaring $W_\epsilon$ of $W$ is cylindrical near $\{a^0,a^1\}$, and where $g\in\mathcal{F}(W)$.  Note that such $W$ need not be compact, but the projection $p_W$ onto $\mathbb{R}^k$ must be proper.  The source and target maps are given by $(W,g)\mapsto (M_i:= pr^{-1}_W(a^i),g_{\mid_{M_i}})$ for $i=0$ and $1$ respectively.  Composition in $\mathsf{Cob}^\mathcal{F}_{d,(k,N)}$ is given by union in $\mathbb{R}^N$ while gluing $\mathcal{F}$-structures.  

\begin{notation}
The additional data $(a^0,a^1)$ of a morphism will often not be written for notational ease though it is still there.
\end{notation}

As discrete categories, there is a functor 
\[
\mathsf{D}^{\mathcal{F},\perp}_{d,(k,N)}\xrightarrow{\alpha} \mathsf{Cob}_{d,(k,N)}^\mathcal{F}
\]
described by 
\[
((a^0,W,g)\prec (a^1,W,g))\mapsto (a^0,a^1,W_{[a^0,a^1]},g_{[a^0,a^1]})
\]
where $W_{[a^0,a^1]} := pr_W^{-1}(\mathbb{R}^{k-1}\times[a^0_k,a^1])$ and $g_{[a^0_k,a^1]} := g_{\mid_{W_{[a^0,a^1]}}}$.  The case $a^0=a^1$ gives the functor on objects.  We declare a point-set map $X\rightarrow mor\text{ }\mathsf{Cob}^\mathcal{F}_{d,(k,N)}$ to be ``continuous" if the pull back
\[
X\times_{mor\text{ }\mathsf{Cob}^\mathcal{F}_{d,(k,N)}} mor\text{ }\mathsf{D}^{\mathcal{F},\perp}_{d,(k,N)}\rightarrow mor\text{ }\mathsf{D}^{\mathcal{F},\perp}_{d,(k,N)}
\]
is ``continuous".  We do similarly for the objects of $\mathsf{Cob}^\mathcal{F}_{d,(k,N)}$.  It is straight forward to verify that $\mathsf{Cob}^\mathcal{F}_{d,(k,N)}$ becomes a quasi-topological category in this way.  Note that there is an equivalence of quasi-topological categories $\mathsf{Cob}^\mathcal{F}_{d,(1,N)} \cong  \mathsf{Cob}^\mathcal{F}_{d,N}$ determined by choosing a homeomorphism $int(I^{N-1})\cong \mathbb{R}^{N-1}$.

\begin{remark}
One could alternatively work with an elaboration of $\mathsf{Cob}^\mathcal{F}_{d,(k,N)}$ regarded as a topological $k$-category whose $1$-category of $(k-1)$-morphisms $mor^{k-1}(\emptyset,\emptyset)$ is $\mathsf{Cob}^\mathcal{F}_{d,N}$.  The higher category structure can be sketched as follows.  

For $a^0,a^1\in\mathbb{R}^k$ write $[a^0,a^1]:=\prod^k_1 [a^0_i,a^1_i]\subset\mathbb{R}^k$ as the product of closed intervals where it is agreed that $[s,t] = \emptyset$ if $s>t$.  An $l$-morphism is a $4$-tuple $(a^0,a^1, W,g)$ such that $a^0_i = a^1_i $ when $i\leq k-l$, and $W$ is \textit{cylindrical} near $[a^0,a^1]$; in particular $[a^0,a^1]$, as a manifold with corners, is transverse to the projection $pr_W:W\rightarrow \mathbb{R}^k$.  The source maps are given by 
\[
(a^0_1,...,a^0_k ; a^1_{k-l+1},...,a^1_k) \mapsto (a^0_1,...,a^0_k;a^1_{k-l+2},...,a^1_k)
\]
The target maps are similar.  Composition is given by concatenating appropriate intervals.

In particular, we can regard $\mathsf{Cob}^\mathcal{F}_{d,(d,N)}$ as a topological $d$-category which in the limit $N\rightarrow \infty$ is an (obviously equivalent) alternative to the $(\infty,d)$-category from \cite{hopkins-lurie:tfts}.
\end{remark}

\begin{lem}[$B\alpha$]
As quasi-topological categories, the functor $\alpha$ induces a weak homotopy equivalence 
\[
B\mathsf{D}^{\mathcal{F},\perp}_{d,(k,N)}\xrightarrow{B\alpha} B\mathsf{Cob}_{d,(k,N)}^\mathcal{F}.
\]

\end{lem}

\begin{proof}
We will show that $\alpha$ induces a level-wise weak homotopy equivalence of simplicial nerves
\[
N_l\mathsf{D}^{\mathcal{F},\perp}_{d,(k,N)}\xrightarrow{\simeq} N_l \mathsf{Cob}_{d,(k,N)}^\mathcal{F}.
\]

\begin{figure}
\begin{center}
	\includegraphics[width=3.5in,viewport=0 70 725 580,clip=true]{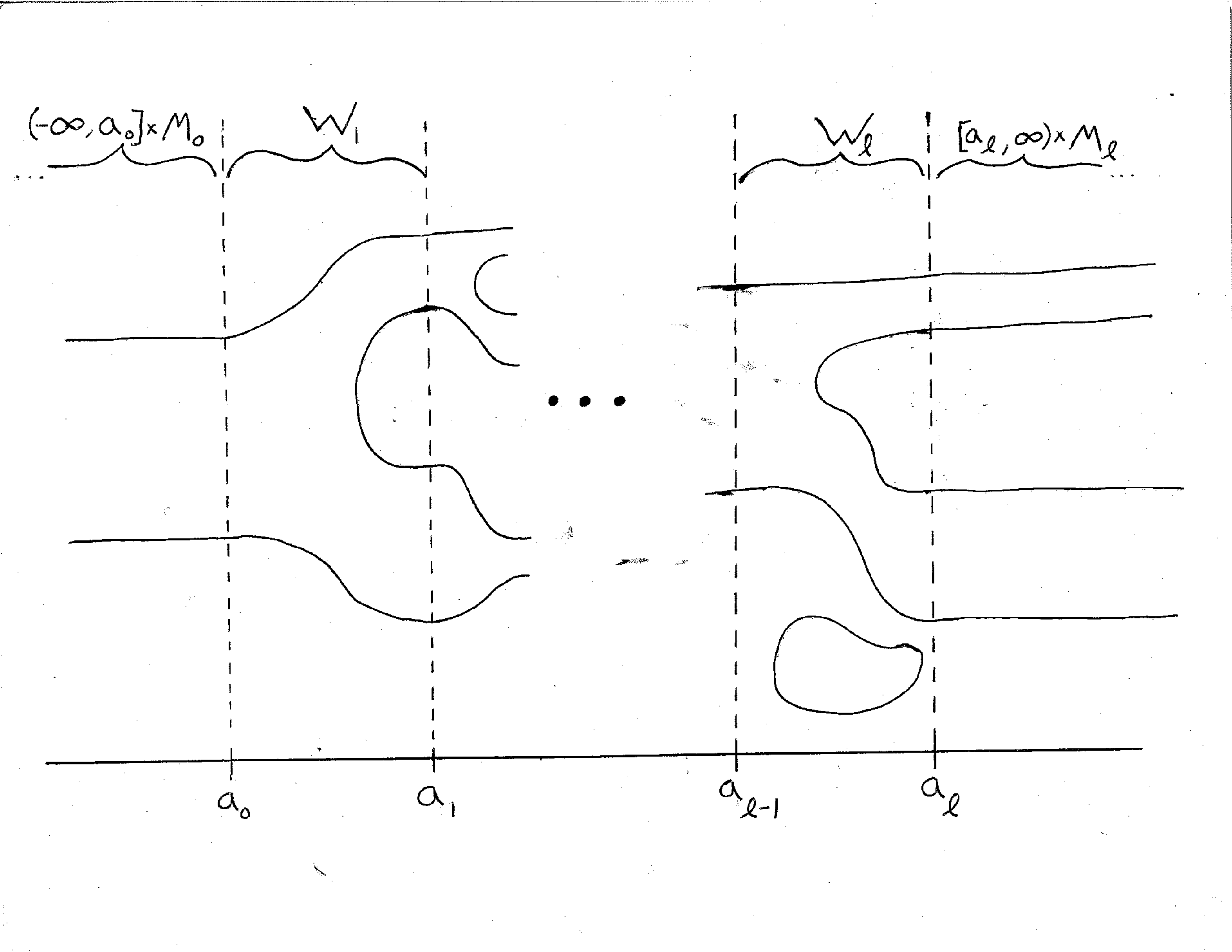}
	\caption{Gluing $W$}
	\label{fig:Gluing-W}
\end{center}
\end{figure}

There is a simplicial inclusion 
\[
N_\bullet \mathsf{Cob}_{d,(k,N)}\hookrightarrow N_\bullet\mathsf{D}^{\perp}_{d,(k,N)}
\]
given by sending 
$
((a^0,a^1,W_1),...,(a^{l-1},a^l,W_l))
$ to 
\[
((a^0\leq...\leq a^l),(-\infty,a_0]\times M_0\cup W_1\cup...\cup W_l \cup [a_l,\infty)\times M_l)
\]
where $t(W_i)=(a^i,M_i)=s(W_{i+1})$ (see Figure~\ref{fig:Gluing-W}).  We show that a similar map exists in compact families while we include $\mathcal{F}$-structure.

Choose once and for all a \textit{smooth} family of homeomorphisms $\psi_\epsilon:(-\infty,\epsilon)\rightarrow \mathbb{R}$, $\epsilon\in (0, \infty]$, with $\psi_\infty =id_\mathbb{R}$ and such that for each $\epsilon>0$, $\psi_\epsilon$ is the identity in a neighborhood of $(-\infty,0]$.

Let $K$ be a compact manifold with a smooth map $K\xrightarrow{f} N_l \mathsf{Cob}^\mathcal{F}_{d,(k,N)}$ in the sense of~\ref{subsec: smooth families}.  Because nowhere in the argument is $l$ relevant, we will take $l=1$.  Recall that for $W\in mor\text{ }\mathsf{Cob}_{d,(k,N)}$ a collared $d$-manifold,  $\mathcal{F}(W)=colim_\epsilon \mathcal{F}(W_\epsilon)$ where $W_\epsilon$ is an $\epsilon$-collaring of $W$.  It follows from the compactness of $K$ that there exists $\epsilon_f>0$ such that $f$ is represented by 
\begin{equation}\label{eq: f epsilon}
f_{\epsilon_f}:K\rightarrow \{(W,g)\mid W\in mor\text{ }\mathsf{Cob}_d\text{ and } g\in\mathcal{F}(W_{\epsilon_f})\}.
\end{equation}
Using the map~(\ref{eq: f epsilon}) along with the homeomorphisms $\psi_\epsilon$, it is possible to achieve a map $\tilde{f}_{\epsilon_f}:K\rightarrow mor\text{ }\mathsf{D}^{\perp}_{d,(k,N)}$ such that $\alpha \circ \tilde{f}_{\epsilon_f} = f$.  Indeed, the $\psi_\epsilon$ allow for us to extend an $\mathcal{F}$-structure on $W_\epsilon$ to an $\mathcal{F}$-structure on $W_\infty\in D_{d,(k,N)}$; the smoothness of $\phi_\epsilon$ in $\epsilon$ ensures we can do so in smooth families.  We have established an inverse map 
\[
[K,mor\text{ }\mathsf{Cob}^\mathcal{F}_{d,(k,N)}]\rightarrow [K,mor\text{ }\mathsf{D}^{\mathcal{F},\perp}_{d,(k,N)}]
\]
and it is clear that the composite of this map with $\alpha$ is the identity on the set $[K,mor\text{ }\mathsf{Cob}^\mathcal{F}_{d,(k,N)}]$.

Now, for each $(a^0\leq a^1)$, choose $\epsilon>0$ small and let 
\[
\lambda_{t,(a^0\leq a^1)}:\mathbb{R}\rightarrow \mathbb{R}
\]
be a smooth non-decreasing function given in most places ($\epsilon$-away from four points) by
\[
\lambda_{t,(a^0\leq a^1)}(x)=
\begin{cases} 
x - 1/(1-t), 
& 
x\in(-\infty,a^0_k-\epsilon - 1/(1-t)]   
\\
a^0_k,
&
x\in [a^0_k+\epsilon -1/(1-t),a^0_k-\epsilon]
\\
x,
&
x\in [a^0_k+\epsilon  , a^1_k - \epsilon]
\\
a^1_k,
&
x\in [a^1_k + \epsilon , a^1_k-\epsilon + 1/(1-t)]
\\
x + 1/(1-t) ,
&
x\in [a^1_k +\epsilon + 1/(1-t) , \infty).
\end{cases}
\]
(see Figure~\ref{fig:lambda-t}).  Choose such a $\lambda_t$ to be smooth in $t\in [0,1]$.

\begin{figure}
\begin{center}
	\includegraphics[width=3.5in,viewport=0 20 540 330,clip=true]{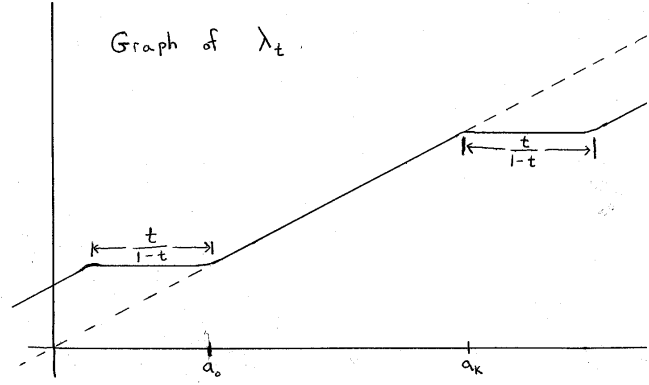}
	\caption{$\lambda_t$}
	\label{fig:lambda-t}
\end{center}
\end{figure}

For each $t\in [0,1]$, consider the self-map $\phi_t $ of $\mathbb{R}^{k-1}\times \mathbb{R}\times I^k$ given by 
\[
\phi_{t,(a^0\leq a^1)}(u,x,v)=(u,\lambda_t(x),v).
\]
Define the self-map
$\Phi_t$ of $mor\text{ }\mathsf{D}^{\mathcal{F},\perp}_{d,(k,N)}$ by 
\[
\Phi_t((a^0\leq a^1),W,g)=((a^0\leq a^l),\phi_{t,(a^0\leq a^l)}^{-1}(W), \phi_{t,(a^0\leq a^l)}^*g).
\]
It is clear that composition with $\Phi_t$ provides the necessary homotopy to conclude that the composition  
\[
[K,mor\text{ }\mathsf{D}^{\mathcal{F},\perp}_{d,(k,N)}]\rightarrow [K,mor\text{ }\mathsf{Cob}_{d,(k,N)}^\mathcal{F}]\rightarrow [K,mor\text{ }\mathsf{D}^{\mathcal{F},\perp}_{d,(k,N)}].
\]
is the identity set-map.
\end{proof}

We record the following proposition.  For $W^d$ a boundaryless $d$-manifold, write $B^\mathcal{F}_{(k,N)}(W) = (Emb^{prop}(W,\mathbb{R}^k\times  I^{N-k})\times \mathcal{F}(W))/{Diff}(W)$ for the space of properly embedded $(W^\prime,g)$ in $\mathbb{R}^k\times  I^{N-k}$ with $W^\prime \cong W$.   If $W$ has collared boundary, replace $\mathbb{R}^k$ by $\mathbb{R}^{k-1}\times [a,a^\prime]$ and consider collared embeddings.  For each $a\in\mathbb{R}^k$ and closed $(d-k)$-manifold $M$ there is a standard inclusion $B^\mathcal{F}_{(0,k)}(M)\xrightarrow{\iota} ob\text{ }\mathsf{Cob}^\mathcal{F}_{d,(k,N)}$.  
For $(a^i,M_i,g_i)\in ob\text{ }\mathsf{Cob}^\mathcal{F}_{d,(k,N)}$, $i=0,1$ with $a^0 < a^1$, there results a pull back square
\[
\xymatrix{
\iota^* mor\text{ }\mathsf{Cob}^\mathcal{F}_{d,(k,N)}  \ar[r]  \ar[d]
&
mor\text{ }\mathsf{Cob}^\mathcal{F}_{d,(k,N)}  \ar[d]^{s,t}
\\
B^\mathcal{F}_{(0,k)}(M_0)\times B^\mathcal{F}_{(0,k)}(M_1)  \ar[r]^{\iota\times \iota}
&
ob\text{ }\mathsf{Cob}^\mathcal{F}_{d,(k,N)}\times ob\text{ }\mathsf{Cob}^\mathcal{F}_{d,(k,N)}
}
\]
where the right vertical map is the source-target map. 
Recall the notion of smoothness from \textsection\ref{subsec: smooth families}.

\begin{prop}\label{prop: fibration}
Suppose $\mathcal{F}$ is a fiberwise equivariant sheaf.  Let $(a^i,M_i,g_i)\in ob\text{ }\mathsf{Cob}^\mathcal{F}_{d,(k,N)}$, $i=0,1$, with $a^0\leq a^1$.  The source-target map
\[
\iota^* mor\text{ }\mathsf{Cob}^\mathcal{F}_{d,(k,N)} \rightarrow B^\mathcal{F}_{(0,k)}(M_0)\times B^\mathcal{F}_{(0,k)}(M_1)
\]
has the homotopy lifting property with respect to smooth compact families.  
\end{prop}

\begin{proof}
Suppose for now that $\mathcal{F}\equiv *$ is trivial.  
We wish to complete the following diagram
\[
\xymatrix{
\{0\}\times X  \ar[d]   \ar[r]^{\tilde{f}}
&
\iota^* mor\text{ }\mathsf{Cob}_{d,(k,N)}  \ar[d]^{s,t}
\\
I\times X  \ar[r]^{f}  \ar@{.>}[ur]
&
B_{(0,k)}(M_0)\times B_{(0,k)}(M_1)
}
\]
where the horizontal maps are smooth and $X$ is compact.  Write $f^s_r$, $r\in I$, for the composition of $f$ with the projection on to the first factor.  Do similarly for $f^t_r$ with the second factor.  Corresponding to $f^t_r$ is the $(q+1+d-k)$-dimensional submanifold $E^t\subset I\times X^q\times(\{a^1\}\times I^{N-k})$.  Consider the space (manifold)
\[
\tilde{E}^t:=\mathbb{R}^{k-1}\times \{(r^\prime,r,e)\mid r^\prime\in[0,r]\text{ and } e\in E_{r^\prime}\}.
\]
There is a canonical inclusion 
$\tilde{E}^t\hookrightarrow \mathbb{R}^{k-1}\times [a^1_k,a^1_k + 1]\times E$ 
given by $(x,r^\prime,r,e)\mapsto (x,r,e)$.  Moreover, the map 
$\tilde{E}^t\rightarrow I\times X$, given by $(r^\prime,r,e)\mapsto (r,pr_X(e))$, makes 
\[
\tilde{E}^t\subset I\times X\times(\mathbb{R}^{k-1}\times [a^1_k,a^1_k +1]\times I^{N-k})
\]
into the total space of a smooth bundle of collared $d$-manifolds in $\mathbb{R}^{k-1}\times [a^1_k,a^1_k+1]\times I^{N-k}$.  Such a bundle corresponds to a map
\[
\tilde{f}^t_r:I\times X\rightarrow B_N(\mathbb{R}^{k-1}\times [a^1_k,a^1_k+1]\times M_0)\hookrightarrow mor\text{ }\mathsf{Cob}_{d,(k,N)}.
\]
Notice that $s\circ \tilde{f}^t_r \equiv f^t_0$ and $t\circ\tilde{f}^t_r = f^t_r$.

One could similarly define the space (manifold) 
\[
\tilde{E}^s:=\{(r^\prime,r,e)\mid 1-r^\prime\in[0,r]\text{ and }e\in E_{r^\prime}\} \subset I\times X\times (\mathbb{R}^{k-1}\times [a^0_k -1,a^0_k]\times I^{N-k})
\]
corresponding to a map $\tilde{f}^s_r$ satisfying $s\circ \tilde{f}^s_r = f^t_r$ and $t\circ \tilde{f}^s_r \equiv f^s_0$.  There is then a map 
\[
I\times X\rightarrow mor\cob\times_{ob} mor\cob\times_{ob} mor\cob
\]
given by $(r,x)\mapsto (\tilde{f}^s_r(x), \tilde{f}(x), \tilde{f}^t_r(x))$.  It is easily verified that upon composing with the composition map, $mor\times_{ob}mor\times_{ob} mor\rightarrow mor$, then rescaling intervals, we get our desired dotted arrow.

Now suppose $\mathcal{F}$ is an arbitrary fiberwise equivariant sheaf.  In this case, $\mathcal{F}$ is determined by a fibration $B_\theta\xrightarrow{\theta} BO(d)$ whose restriction to $Gr_{d,N}$ is again written $B_\theta\rightarrow Gr_{d,N}$; see \textsection\ref{subsubsec: alternative description}.  So a smooth bundle of $\mathcal{F}$-manifolds $X\xrightarrow{f} B^\mathcal{F}_{(0,k)}(M_0)$ over compact $X$ corresponds to a subspace $E_\theta \subset X\times(\mathbb{R}^{k-1}\times\{a\}\times I^{N-k})\times B_\theta$ such that the projection away from $B_\theta$ is a smooth bundle over $X$.  Following the lines above, we can associate to such a bundle $E_\theta\rightarrow I\times X$ the space
\[
\tilde{E}^t_\theta:=\{(r^\prime,r,e)\mid r^\prime\in[0,r]\text{ and } e\in (E_\theta)_{r^\prime}\}.
\]
It should be apparent how the argument from here is the same as that above.
\end{proof}

\begin{remark}
The above proposition is also true when $\mathcal{F}$ is \textit{flexible} in the sense of~\cite{gromov:h-principle}, page 76.  The above lemma is not true for general $\mathcal{F}$.  Such a non-example is achieved by the $2$-dimensional equivariant sheaf of holomorphic curves in a fixed complex manifold $Y$, $\mathcal{F}(W) = \{(J,h)\mid h\in Hol(W,J),Y)\}$.  Indeed, it is rarely possible to extend a given perturbation of the boundary of a holomorphic curve to the entire curve.  
\end{remark}

\subsubsection*{Comparing $D^\mathcal{F}_{d,(k+1,N)}$ and $\Omega B\mathsf{Cob}^\mathcal{F}_{d,(k,N)}$}\label{subsubsec: comparing D and omega B}

For $n\in\mathbb{Z}$, write $\overline{n}\in\mathbb{R}^k$ for the $k$-tuple $(0,...,0,n)$.  There are the represented functors
\[
mor(-,(\overline{n},\emptyset)): {\mathsf{Cob}_{d,(k,N)}^\mathcal{F}}^{op}\rightarrow \mathsf{QTop}
\]
as the Q-space of morphisms in $\mathsf{Cob}_{d,(k,N)}^\mathcal{F}$ to the empty manifold and the empty $\mathcal{F}$-structure.  Regard $\emptyset$ as a morphism $(\overline{n},\emptyset)\xrightarrow{\emptyset} (\overline{n+1},\emptyset)$.  There is a directed system of Q-spaces, natural in $(a,M,g)\in ob\text{ }\mathsf{Cob}_{d,(k,N)}^\mathcal{F}$, given by post-composition
\[
\dotsi\xrightarrow{\circ \emptyset} mor((a,M,g), (\overline{n},\emptyset))\xrightarrow{\circ \emptyset} mor((a,M,g),(\overline{n+1},\emptyset))\xrightarrow{\circ \emptyset}\dotsi
\]

In this way, define the functor
\[
F:{\mathsf{Cob}_{d,(k,N)}^\mathcal{F}}^{op}\rightarrow \mathsf{QTop}
\]
given by $-\mapsto hocolim(\dotsi\xrightarrow{\circ\emptyset} mor(-,(\overline{n},\emptyset))\xrightarrow{\circ\emptyset} \dotsi)$.  
Consider the wreath category
\[
\mathsf{Cob}_{d,(k,N)}^\mathcal{F}\wr F
\]
whose objects are pairs $(o,x)$ with $x\in F(o)$ for $o\in ob\text{ }\mathsf{Cob}_{d,(k,N)}^\mathcal{F}$ and a morphism $(o_0,y)\xrightarrow{\alpha} (o_1,x)$ is given by $o_0\xrightarrow{\alpha} o_1 \in mor\text{ }\mathsf{Cob}_{d,(k,N)}^\mathcal{F}$ providing $F(\alpha)(x)=y$.  There is a bijection 
\[
ob\text{ }\mathsf{Cob}^\mathcal{F}_{d,(k,N)}\wr F \leftrightarrow \{(W,g)\in mor\text{ }\mathsf{Cob}^\mathcal{F}_{d,(k,N)}\mid t(W,g)=(\overline{n},\emptyset)\text{ for some $n\in\mathbb{N}$}\}
\]
and an inclusion
\[
mor\text{ }\mathsf{Cob}^\mathcal{F}_{d,(k,N)}\wr F \subset mor\text{ }\mathsf{Cob}^\mathcal{F}_{d,(k,N)}
\]
from which $\mathsf{Cob}^\mathcal{F}_{d,(k,N)}\wr F $ inherits the structure of a quasi-topological category. There is an obvious forgetful functor
\[
\mathsf{Cob}_{d,(k,N)}^\mathcal{F}\wr F \rightarrow \mathsf{Cob}_{d,(k,N)}^\mathcal{F}.
\]

\begin{lem}[$\delta$]\label{lem: deloop}
The map 
\[
B\mathsf{Cob}_{d,(k,N)}^\mathcal{F}\wr F \rightarrow B\mathsf{Cob}_{d,(k,N)}^\mathcal{F} 
\]
induced by the forgetful functor is a quasi-fibration over the path-component of $\emptyset$ with fiber homeomorphic to the Q-space $D^\mathcal{F}_{d,(k-1,N)}$.
\end{lem}

The proof of Lemma~\ref{lem: deloop} occupies \textsection\ref{subsec: proof of lemma}.

Because the functor $F$ is built from represented functors, the classifying space $B\mathsf{Cob}_{d,(k,N)}^\mathcal{F}
\wr F$ is contractible via the self-functor described by 
\[
(a,M,g,W)\xrightarrow{W} (\overline{n},\emptyset)
\]
for $n$ large enough.  We immediately obtain the

\begin{cor}\label{cor:delooping}
The inclusion $\delta$ of the fiber into the homotopy-fiber of the above quasi-fibration is a weak homotopy equivalence.  That is, 
\[
\delta: D^\mathcal{F}_{d,(k-1,N)} \xrightarrow{\simeq} \Omega B\mathsf{Cob}_{d,(k,N)}^\mathcal{F}.
\]
\end{cor}

We have established the zig-zag of weak homotopy equivalences
\[
B\mathsf{Cob}^\mathcal{F}_{d,(k,N)}\xleftarrow{\alpha} B\mathsf{D}^{\mathcal{F},\perp}_{d,(k,N)}\xrightarrow{\beta} B\mathsf{D}^{\mathcal{F},\pitchfork}_{d,(k,N)}\xrightarrow{\gamma}D^\mathcal{F}_{d,(k,N)}\xrightarrow{\delta} \Omega B\mathsf{Cob}^\mathcal{F}_{d,(k+1,N)}
\]
resulting in a weak homotopy equivalence 
\[
B\mathsf{Cob}_{d,N}^\mathcal{F}\simeq \Omega^{N-1} D^\mathcal{F}_{d,(N,N)}= \Omega^{N-1}\Psi_d^\mathcal{F}(\mathbb{R}^N).
\]
It remains to identify the weak homotopy type of $\Psi^\mathcal{F}_d(\mathbb{R}^N)$.

\subsection{The weak homotopy type of $\Psi^\mathcal{F}_d(\mathbb{R}^N)$}\label{subsec: homotopy type of psi}

\subsubsection*{Statement of the Weak Homotopy Type}

Let $Th$ denote the Thom space functor.
Let $\gamma_{d,N}$ denote the tautological bundle over $Gr_{d,N}$ of $d$-planes in $\mathbb{R}^N$ with total space $U_{d,N}$, and let $\gamma^\perp_{d,N}$ denote its orthogonal compliment with total space $U_{d,N}^\perp$.  Observe that $U_{d,N}^\perp$ can be identified with the space of affine $d$-planes in $\mathbb{R}^N$ as seen by the assignment $(V,v)\mapsto v+V\subset\mathbb{R}^N$.  As such, there is an evident map 
\[
U^\perp_{d,N}\rightarrow \Psi_d(\mathbb{R}^N).
\]
There is a continuous extension of this map to the one-point compactification 
\[
Th(\gamma^\perp_{d,N})\rightarrow \Psi_d(\mathbb{R}^N)
\]
by insisting $\infty\mapsto \emptyset$.

Consider the fibration 
\[
\overline{p}_\mathcal{F}:Fr_{d,N}\times_{O(d)} \mathcal{F}(\mathbb{R}^d)\rightarrow Gr_{d,N}
\]
with fiber $\mathcal{F}(\mathbb{R}^d)$
similar to that discussed in \textsection\ref{subsec: trad sheaves of sections}.  For explicitness, 
$\overline{p}^*_\mathcal{F}U^\perp_{d,N}=\{(V,v,g)\mid V\in Gr_{d,N}, \text{ } \mathbb{R}^N\ni v\perp V, \text{ and } g\in \mathcal{F}(V)\}$. 
There is a similar map as that directly above,
\[
Th(\overline{p}_\mathcal{F}^*\gamma^\perp_{d,N})\rightarrow \Psi_d^\mathcal{F}(\mathbb{R}^N),
\]
given by $(V,v,g)\mapsto (V+v,g)$, $\infty\mapsto \emptyset$.

\begin{lem}\label{lem: thom space}
This map $Th(\overline{p}_\mathcal{F}^*\gamma^\perp_{d,N})\rightarrow \Psi_d^\mathcal{F}(\mathbb{R}^N)$ is a weak homotopy equivalence.
\end{lem}

To conclude the proof of the main theorem, bring to attention the above fibration of interest $\overline{p}_\mathcal{F}: B^\mathcal{F}\rightarrow BO(d)$ having fiber $\mathcal{F}(\mathbb{R}^d)$.  The main theorem concerns the fibration $p_\mathcal{F}:Gr^\mathcal{F}_{d,N}\rightarrow Gr_{d,N}$ having fiber $Stalk_0(\mathcal{F}_{\mid_{\mathbb{R}^d}})$.   These two fibrations are weakly homotopy equivalent from Proposition~\ref{prop: stalk}.  As a corollary of Lemma~\ref{lem: thom space} we have

\begin{theorem*}[Main Theorem]
There is a weak homotopy equivalence
\[
B\mathsf{Cob}^\mathcal{F}_{d,N} \simeq \Omega^{N-1} Th(p^*_\mathcal{F}\gamma^\perp_{d,N}).
\]
\end{theorem*}

\subsubsection*{Strategy for Proving Lemma~\ref{lem: thom space}}

We consider the following open subsets of $\Psi_d^\mathcal{F}(\mathbb{R}^N)$:
\[
V_0=\{(W,g)\in\Psi_d^\mathcal{F}(\mathbb{R}^N)\mid 0\notin W\}
\]
and 
\[
V_{min}=\{(W,g)\in\Psi_d^\mathcal{F}(\mathbb{R}^N)\mid \lVert - \rVert \text{ has a unique nondegenerate minimum on } W\}.
\]
Observe that the map at hand $Th(\overline{p}_\mathcal{F}^*\gamma^\perp_{d,N})\rightarrow \Psi_d^\mathcal{F}(\mathbb{R}^N)$ restricts appropriately to yield a morphism of (horizontal) diagrams
\[
\xymatrix{
(\overline{p}^*_\mathcal{F}U^\perp_{d,N}  \ar[d]
&
\overline{p}^*_\mathcal{F}(U^\perp_{d,N}\setminus\{0\})  \ar[l]  \ar[r]   \ar[d]
&
\{\infty\})  \ar[d]
\\
(V_{min}
&
V_{min}\cap V_0   \ar[l]  \ar[r]
&
V_0)
}
\]
We establish that each vertical map is a weak homotopy equivalence.  It follows that the map of homotopy colimits of these (horizontal) diagrams is a weak homotopy equivalence.

To finish, we  recognize 
\[
Th(\overline{p}^*_\mathcal{F}\gamma^\perp_{d,N})=hocolim(\overline{p}^*_\mathcal{F}U^\perp_{d,N}\hookleftarrow \overline{p}^*_\mathcal{F}U^\perp_{d,N}\setminus\{0\}\rightarrow *)
\]
and, in this situation,
\[
\Psi_d^\mathcal{F}(\mathbb{R}^N)=V_0\cup V_{min}\xleftarrow{\simeq} hocolim(V_{min}\leftarrow
V_{min}\cap V_0\rightarrow  
V_0).
\]

\subsubsection*{The Weak Homotopy Types of the $V_-$'s}

\begin{lem}[Homotopoy Type of $V_0$]

There is a homotopy equivalence 
\[
V_0\simeq *.
\]

\end{lem}

\begin{proof}

For $t>0$, let $\lambda_t:\mathbb{R}^N\rightarrow \mathbb{R}^N$ be scaling by $1/t$.  It it apparent that $\lambda^*_t:\Psi^\mathcal{F}_d(\mathbb{R}^N)\rightarrow \Psi^\mathcal{F}_d(\mathbb{R}^N)$ restricts to a self-map of $V_0$ and as such, it is possible to continuously extend $\lambda_t$ at $t=0$ by $\lambda_0(W) \equiv \emptyset$.   Because $\lambda_1$ is the identity map, we have demonstrated a homotopy equivalence $V_0\simeq * = \{\emptyset\} $.  
\end{proof}

\begin{lem}[The Weak Homotopy Type of $V_{min}$]

The inclusion
\[
\overline{p}^*_\mathcal{F}U^\perp_{d,N}\hookrightarrow V_{min},
\]
given by $(V,v,g)\mapsto (V+v,g)$ as above,
is a weak homotopy equivalence.  Indeed, $v\in(V+v)$ is the unique closest point in $V+v$ to $0\in\mathbb{R}^N$.

\end{lem}

\begin{proof}

Choose once and for all a continuous family of embeddings $\phi_r:\mathbb{R}^N\rightarrow\mathbb{R}^N$, $r\in (0,\infty]$, with $\phi_\infty = id_{\mathbb{R}^N}$ and such that for each $r$, $\phi_r$ is the identity in a neighborhood of the origin and has image $D^N_r\subset\mathbb{R}^N$, the open $r$-disk about the origin.  For $r^\prime\in[r,\infty]$, the $\phi_{r^\prime}$ themselves realize a homotopy from $\Psi^\mathcal{F}_d(\mathbb{R}^N)\xrightarrow{\phi^*_r}\Psi^\mathcal{F}_d(\mathbb{R}^N)$ to the identity map.

Write $W_{\mid_r}=\phi^{-1}_r(W-p)$.  Let $pr_{T_p W}: W\rightarrow T_p W$ be the projection in $\mathbb{R}^N$ onto the \textit{affine} subspace $T_p W\subset\mathbb{R}^N$.  The inverse function theorem implies the existence of $r>0$ such that 
\[
pr_{T_0 W_{\mid_r}}:W_{\mid_r}\rightarrow T_0 W_{\mid_r}
\]
is a diffeomorphism.  Write $\pi_r$ for the inverse diffeomorphism.   Let $0<r_W\leq \infty$ be the supremum of such $r > 0$.  The map $V_{min}\rightarrow (0,\infty]$, given by the assignment $(W,g)\mapsto r_W$, is upper-semicontinuous.  For $(W,g)\in V_{min}$ and $r>0$ as above, the map 
\[
Pr_{t,W_{\mid_r}}:I\times W_{\mid_r}\rightarrow \mathbb{R}^N,
\]
given by 
\[
w\mapsto  (1-t) pr_{T_0 W_{\mid_r}}(w) + t w,
\]
is an isotopy of embeddings of $W_{\mid_r}$ from $pr_{T_0 W_{\mid_r}}$ to $W_{\mid_r}\hookrightarrow \mathbb{R}^N$.

Write ${\phi_r}_{\mid_W}$ for the embedding $W_{\mid_r}\rightarrow W$ given by $w\mapsto \phi_r(w) +p$.  We obtain the sequence of maps of Q-spaces
\begin{equation}\label{eq: composition}
\mathcal{F}(W)
\xrightarrow{{{\phi_r}^*_{\mid_W}}} 
\mathcal{F}(W_r)
\xrightarrow{\pi^*_r} 
\mathcal{F}(T_0 W_{\mid_r})
\xleftarrow{{{\phi_r}^*_{\mid_{T_p W}}}} 
\mathcal{F}(T_pW).
\end{equation}
We have already indicated that each of the maps above is a homotopy equivalences.  Indeed, the second map $\pi^*_\epsilon$ is a homotopy equivalence induced from the homotopy $Pr_{t,W_{\mid_r}}$ above; ${\phi_{r^\prime}}_{\mid_-}$, $r^\prime \in [-,\infty]$, provides an explicit homotopy inverse for the first and third maps above.  With this choice of homotopy inverse, the resulting composition, written as
\[
\Phi^*_r:\mathcal{F}(W)\rightarrow \mathcal{F}(T_p W),
\]
in ~\eqref{eq: composition}
is a homotopy equivalence.
The homotopy class of $\Phi_r$ is independent of the (valid) choice of $r>0$  as illustrated by the family $\phi_r$.

Now let $K$ be a finite, and therefore compact, CW complex with a map $K\xrightarrow{f} V_{min}$.  For $x\in K$, write $f(x)=(W_x,g_x)$.  Because the assignment $K\rightarrow (0,\infty]$ given by $x\mapsto \epsilon_{W_x}$ is upper-semicontinuous, for each  such map $f$, there is an $r_f>0$ such that for every $x\in K$ the projection $pr_{T_p W_x}$ is a diffeomorphism.  In this way, from $f\in Map(K,V_{min})$ we obtain a map $\tilde{f}\in Map(K,\overline{p}^*_\mathcal{F}U^\perp_{d,N})$ by 
\[
x\mapsto (T_p W_x , 
\Phi^*_{r_f}(g_x))\in \overline{p}^*_\mathcal{F}U^\perp_{d,N}.
\]
The homotopy class of $\tilde{f}$ is independent of the choice of $r_f>0$.

The assignment $f\mapsto \tilde{f}$ is describes a well-defined set-map $[K, \overline{p}^*_\mathcal{F}U^\perp_{d,N}]\rightarrow [K,V_{min}]$.  It is clear that the resulting composition 
\[
[K, \overline{p}^*_\mathcal{F}U^\perp_{d,N}]
\rightarrow 
[K,V_{min}]
\rightarrow 
[K,  \overline{p}^*_\mathcal{F}U^\perp_{d,N}]
\]
is the identity set-map.

We will now demonstrate a homotopy from $f$ to $\tilde{f}$ from which it will follow that the composition 
\[
[K, V_{min}]\rightarrow [K,\overline{p}^*_\mathcal{F}U^\perp_{d,N}]\rightarrow [K,V_{min}]
\]
is also the identity set-map, thus proving the weak homotopy equivalence $\overline{p}^*_\mathcal{F}U^\perp_{d,N}\simeq V_{min}$.  But we have already taken the care to demonstrate such a homotopy.  Indeed, in the definition of $\tilde{f}$, choosing a homotopy inverse to the homotopy equivalence $\Phi_{r_f}$ provides a homotopy $\tilde{f}\simeq f$. 
\end{proof}

\begin{lem}[Weak Homotopy Type of $V_{min}\cap V_0$]
The inclusion
\[
\overline{p}^*_\mathcal{F}(U^\perp_{d,N}\setminus \{0\})\hookrightarrow V_{min}\cap V_0,
\]
given by $(V,v,g)\mapsto (V+v,g)$ as above,
is a weak homotopy equivalence.
\end{lem}

\begin{proof}

This follows from the proof of the above lemma and the following observation.  For $K$ a compact space, the weak homotopy inverse $[K,V_{min}]\rightarrow [K, \overline{p}^*_\mathcal{F}U^\perp_{d,N}]$ described above as the construction of $\tilde{f}$, restricts to a (weak) homotopy inverse 
\[
[K,V_{min}\cap V_0]\rightarrow [K, \overline{p}^*_\mathcal{F}(U^\perp_{d,N}\setminus \{0\})].
\]
Indeed, nowhere in the construction of $\tilde{f}$ from $f$ is there a dependence on $p_x\in W_x$, the unique closest point in $W_x$ to $0\in\mathbb{R}^N$.  
\end{proof}

\subsection{Proof of Lemma~\ref{lem: deloop}}\label{subsec: proof of lemma}

\subsubsection*{Identifying the Fiber}

The fiber of $B\mathsf{Cob}^\mathcal{F}_{d,(k,N)}\wr F \rightarrow B\mathsf{Cob}^\mathcal{F}_{d,(k,N)}$ over $(0,\emptyset)\in ob\text{ }\mathsf{Cob}^\mathcal{F}_{d,(k,N)}$ is the Q-space of morphisms
\[
F(0,\emptyset)=\{(\overline{0},\overline{n},W,g)\in mor\text{ }\mathsf{Cob}^\mathcal{F}_{d,(k,N)}\mid pr^{-1}_k(\{0,n\}))=\emptyset\}
\]
where $pr_k:\mathbb{R}^k\times I^{N-k} \rightarrow\mathbb{R}$ is projection onto the $k$th coordinate.  A diffeomorphism $(0,\infty)\cong (0,1)$ induces a homeomorphism from $F(0,\emptyset)$ to the Q-space
\[
\{(W,g)\subset \mathbb{R}^k\times I^{N-k}\mid pr_k(W) \subset (0,\infty)\}.
\]
This latter Q-space is $D^\mathcal{F}_{d,(k+1,N)}$ by definition.

\subsubsection*{A Preliminary Observation}

In the proof of the lemma, we will make  repeated implicit use of the following observation.


Consider the diagram 
\[
\xi:=(ob\text{ }\mathsf{Cob}^\mathcal{F}_{d,(k,N)}\wr F  
\xrightarrow{forget}
ob\text{ }\mathsf{Cob}^\mathcal{F}_{d,(k,N)}).
\]
The fiber of $\xi$ over $(a,M,g)$ is the Q-space $F(a,M,g)$.  Consider a diagram of the same shape 
\[
End(\xi) \rightarrow ob\text{ }\mathsf{Cob^\mathcal{F}}_{d,(k,N)}\times ob\text{ }\mathsf{Cob}^\mathcal{F}_{d,(k,N)}
\]
whose fiber over $((a^0,M_0,g_0),(a^1,M_1,g_1))$ is the Q-space of maps $Map(F(a^1,M_1,g_1),F(a^0,M_0,g_0))$.  Because of the construction of $F$ and because our cobordism category is quasi-topological, there is a continuous map
\[
mor_{\mathsf{Cob}^\mathcal{F}_{d,(k,N)}}\rightarrow End(\xi)
\]
over $ob\text{ }\mathsf{Cob}^\mathcal{F}_{d,(k,N)}\times ob\text{ }\mathsf{Cob}^\mathcal{F}_{d,(k,N)}$.  It follows from this that a path of morphisms $(W_t,g_t)$, $t\in[0,1]$ induces a continuous family of maps $F(t(W_t,g_t))\xrightarrow{F(W_t,g_t)} F(s(W_t,g_t))$.


\subsubsection*{Strategy}

We use a group completion argument along the lines of~\cite{mcduff-segal:group-completion} and~\cite{dwyer:centralizer} (3.12) and later~\cite{tillmann:mapping-class-group}.   This is a version of Quillen's Theorem B in the setting of topological categories (simplicial categories).  As so, to show the map in the lemma is a quasi-fibration, we need only verify that  $F(o_1)\xrightarrow{F(m)} F(o_0)$ is a homotopy equivalence for each $o_0\xrightarrow{m} o_1 \in mor\text{ }\mathsf{Cob}^\mathcal{F}_{d,(k,N)}$.  In what follows, we show that this criterion holds for morphisms $(a^0,M_0,g_0)\xrightarrow{(W,g)}(a^1,M_1,g_1)$ with either $M_0$ or $M_1$, possibly both, the empty manifold.  Breaking into these cases is only to make the general cases more conceptually clear.  In fact, the general case is simply an `amalgamation' of the previous cases with the details being the same.  The general case will not be presented.

We verify the required criterion as follows.  Take $(a^0,M_0,g_0)\xrightarrow{(W,g)}(a^1,M_1,g_1)$ to be a morphism in $\mathsf{Cob}^\mathcal{F}_{d,(k,N)}$.  Write $\overline{x}:=(0,...,0,x)\in\mathbb{R}^k$.
We can assume $a^0=\overline{0}$ and $a^1=\overline{1}$.  We demonstrate a path in $mor{\cob}$ from $(W,g)$ to some $(W_+,g_+)$.  From the preliminary observation, such a path prescribes a homotopy from $F(\overline{1},M_1,g_1)\xrightarrow{F(W,g)}F(\overline{0},M_0,g_0)$ to $F(W_+,g_+)$.  We construct a morphism $(\overline{-1},M_{-1},g_{-1})\xrightarrow{(W_{-1},g_{-1})}(\overline{0},M_0,g_0)$, with $(M_{-1},g_{-1})$ built easily from $(M_1,g_1)$, and a path in $mor\cob$ from $(W_+,g_+)\circ (W_-,g_-)$ to a `trivial' morphism $(\overline{W},\overline{g})$.  It will be clear from the sense in which $(\overline{W},\overline{g})$ is trivial that $F(\overline{W},\overline{g})$ is homotopic to $id_{F(\overline{1},M_1,g_1)}$.  It follows that $F(W_-,g_-)$ is a left homotopy inverse to $F(W_+,g_+)$.  The situation presented will be  symmetric under exchanging the indicies $\pm $ and thus $F(W_-,g_-)$ is also a right homotopy inverse and therefore $F(W,g)$ is a homotopy equivalence.

\subsubsection*{Useful tools and simplifications}

Write the coordinates of $\mathbb{R}^k\times I^{N-k}$ as $(x,t)$ with corresponding coordinate vectors as $e_{x_i}$ and $e_{t_i}$.  Let
\[
pr_i:\mathbb{R}^k\times I^{N-1}\rightarrow \mathbb{R},
\]
be projection onto the $x_i$-coordinate and
\[
pr_{\mathbb{R}^{k-1}}:\mathbb{R}^k\times I^{N-k}\rightarrow\mathbb{R}^{k-1}
\]
be projection onto the first $k-1$ coordinates.  Write 
\[
R:\mathbb{R}^k\times I^{N-k} \rightarrow  \mathbb{R}^k\times I^{N-k}
\]
for rotation by angle $\pi$ in the $2$-plane $span\{e_{x_{k-1}},e_{x_k}\}$.

\begin{figure}
\begin{center}
	\includegraphics[width=4in,viewport=70 120 720 540,clip=true]{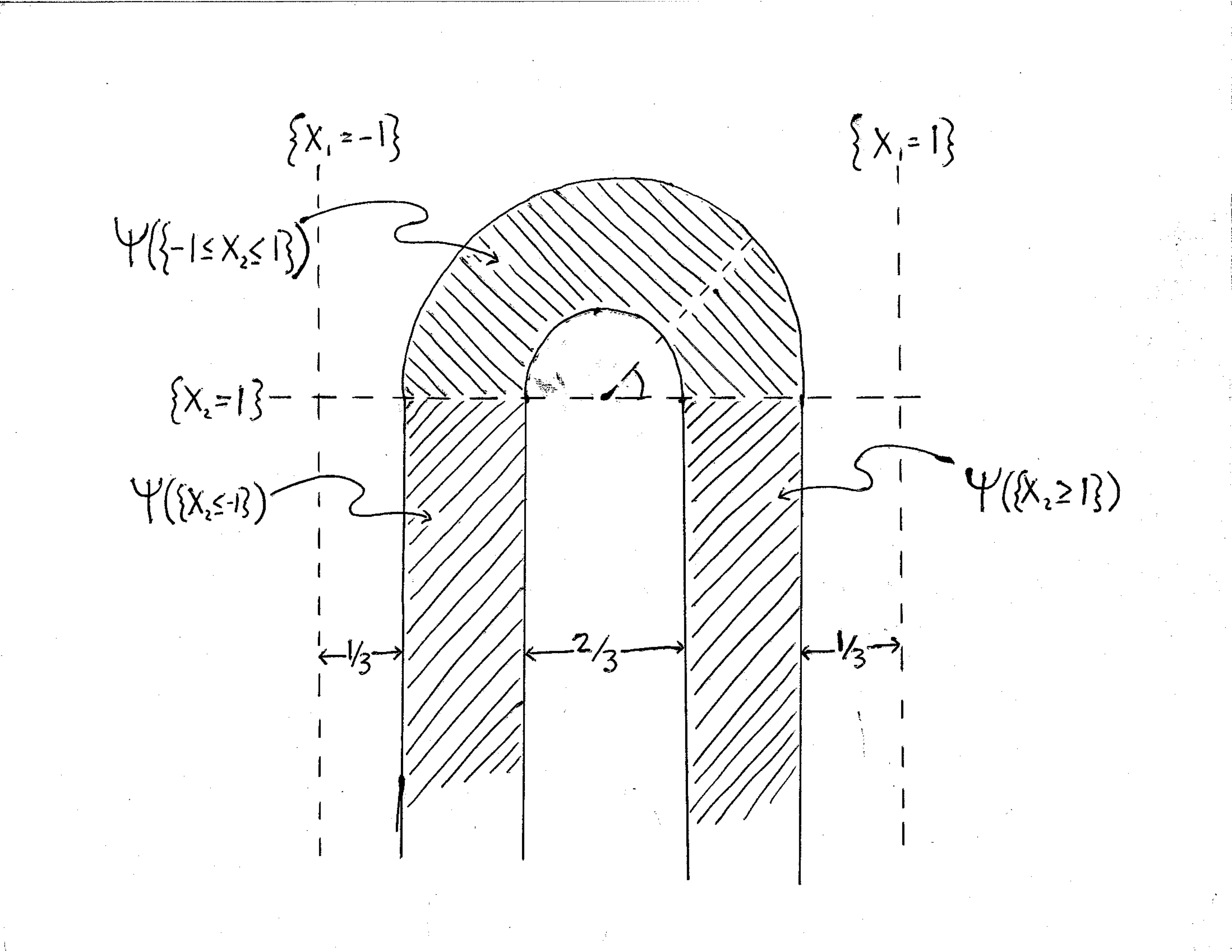}
	\caption{$\psi$}
	\label{fig:psi}
\end{center}
\end{figure}

Choose an embedding $\psi$ from $\mathbb{R}^{k-1}\times [-1,1]\times I^{N-k}$ into itself which is the identity on all but the $(x_{k-1},x_k)$-coordinates and is described on $\{(x_{k-1},x_k)\}$ as follows.   Take $\delta>0$ to be sufficiently small.  For $x_{k-1}\leq -1$,
\[
(x_{k-1},x_k) \mapsto (x_{k-1}+2,x_k/3+1/3),
\]
for $x_{k-1} \geq 1$, 
\[
(x_{k-1},x_k) \mapsto  R_{\mid_{\{(x_{k-1},x_k)\}}} (x_{k-1}-2,x_k/3+1/3),
\] 
for $-1+\delta<x_k<1-\delta$, 
\[
(x_{k-1},x_k)\mapsto ((x_k/3+1/3)sin((\pi/2)(x_{k-1}-1)) , (x_k/3+1/3)cos((\pi/2)(x_{k-1}-1))),
\]
and for $1-\delta\leq \pm x_{k-1} \leq 1$,
\[
\mp \partial_{x_{k-1}} \psi> 1/2.
\]
See figure Figure~\ref{fig:psi}.  Note that the image of $\psi$ has $k$th coordinate bounded above by $2$.

Choose once and for all a continuous family 
\[
\phi_r:\mathbb{R}^{k-1}\rightarrow\mathbb{R}^{k-1}  \text{\hspace{.7cm}}  (r\in(0,\infty])
\]
of embeddings with $\phi_\infty =id_{\mathbb{R}^{k-1}}$ and such that for each $r>0$, $\phi_r$ has image the open $r$-disk about $0\in\mathbb{R}^{k-1}$ and is the  identity in a neighborhood of $0\in\mathbb{R}^{k-1}$.  Let $\Phi_r: = \phi_r\times id_\mathbb{R}\times id_{I^{N-k}}$ be the resulting self-map of $\mathbb{R}^k\times I^{N-k}$.  It will be  used repeatedly and without saying that the $\phi_{r^\prime}$, $r^\prime\in[r,\infty]$, themselves realize a homotopy from $id_{\mathbb{R}^k\times I^{N-k}}$ to $\Phi_r$.

For a cobordism $(a^0,M_0,g_0)\xrightarrow{(W,g)}(a^1,M_1,g_1)$ in $\mathsf{Cob}^\mathcal{F}_{d,(k,N)}$ we can assume without loss in generality that $a^0=\overline{0}$ and $a^1=\overline{1}$.  Let $a\in\mathbb{R}^{k-1}$ be transverse to the projection 
\[
pr:=pr_W:W\rightarrow \mathbb{R}^{k-1}.
\]
Without loss of generality assume $a=0$.  The manifold $P^{d-k+1} := pr^{-1}(a)\subset \{a\}\times [0,1] \times I^{N-k}$ is compact and collared.  The inverse function theorem implies the existence of $r>0$ such that for $(W^\prime,g^\prime):=(\Phi^{-1}_r(W),\Phi^*_r g) $, the projection 
\[
pr:W^\prime\rightarrow P\times\mathbb{R}^{k-1}
\]
is a diffeomorphism.  In fact, there is an isotopy of embeddings 
\[
t\mapsto (1-t)\iota + (t)pr
\]
from the inclusion $\iota : W^\prime \hookrightarrow \mathbb{R}^{k-1}\times [0,1]\times I^{N-k}$ to $pr$.  This homotopy concatenated with the homotopy induced from $\Phi_{r^\prime}$, $r^\prime\in[r,\infty]$, describes a path in $mor\text{ }\mathsf{Cob}^\mathcal{F}_{d,(k,N)}$ from $(W,g)$ to $(\mathbb{R}^{k-1}\times P,p_*g^\prime)$.   We can thus assume $W=\mathbb{R}^{k-1}\times P$ is `straight' in the $\mathbb{R}^{k-1}$-direction. Notice that this procedure does not change the underlying $(d-k)$-manifolds of the source and target.

\subsubsection*{Case $(\overline{0},\emptyset)\xrightarrow{(W,g)}(\overline{1},\emptyset)$}

\begin{figure}
\begin{center}
	\includegraphics[width=3in,viewport=70 0 680 600,clip=true]{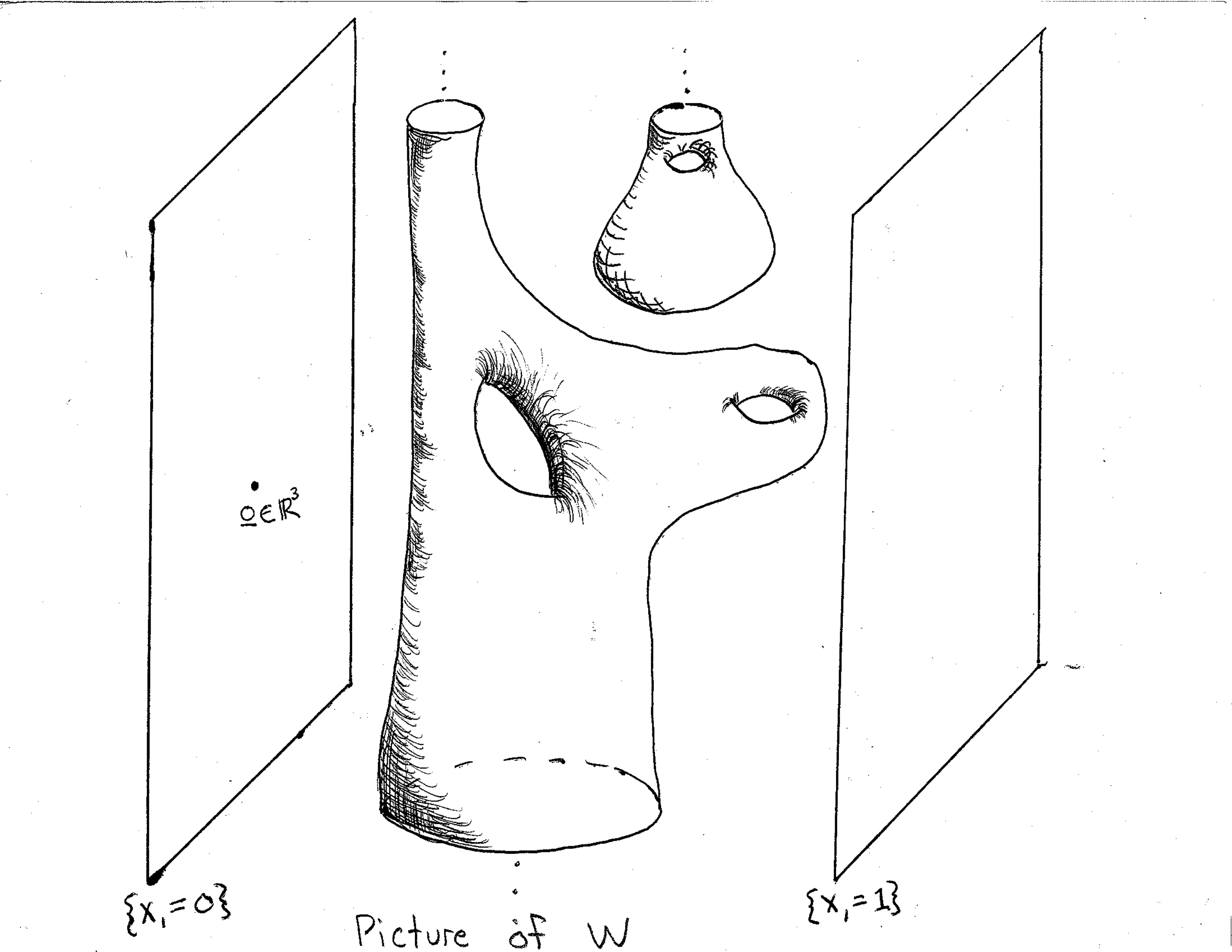}
	\caption{$W$}
	\label{fig:W}
\end{center}
\end{figure}

\begin{figure}
\begin{center}
	\includegraphics[width=4in,viewport=20 15 737 510,clip=true]{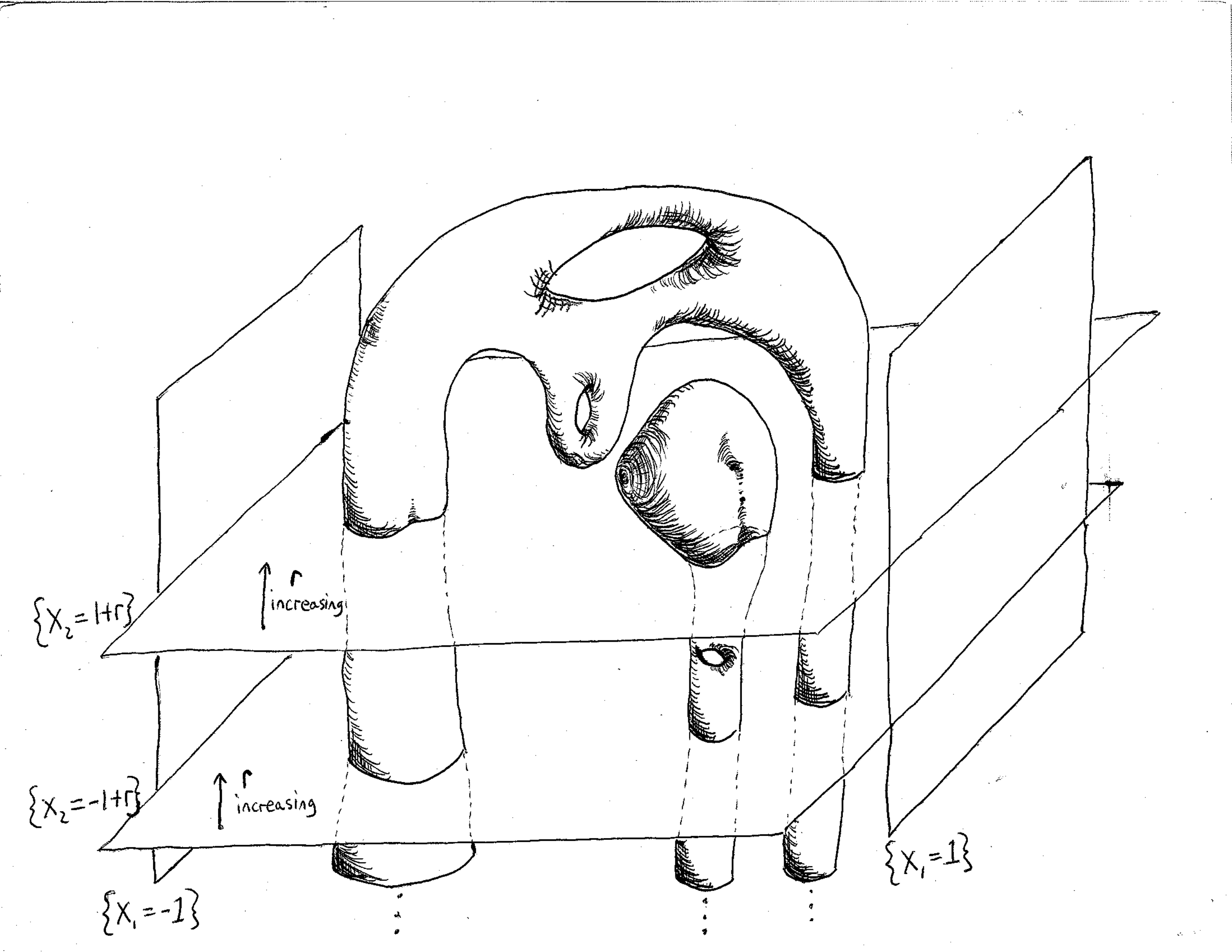}
	\caption{$\psi(W)$}
	\label{fig:psi-W}
\end{center}
\end{figure}

\begin{figure}
\begin{center}
	\includegraphics[width=3.2in,viewport=65 20 710 565,clip=true]{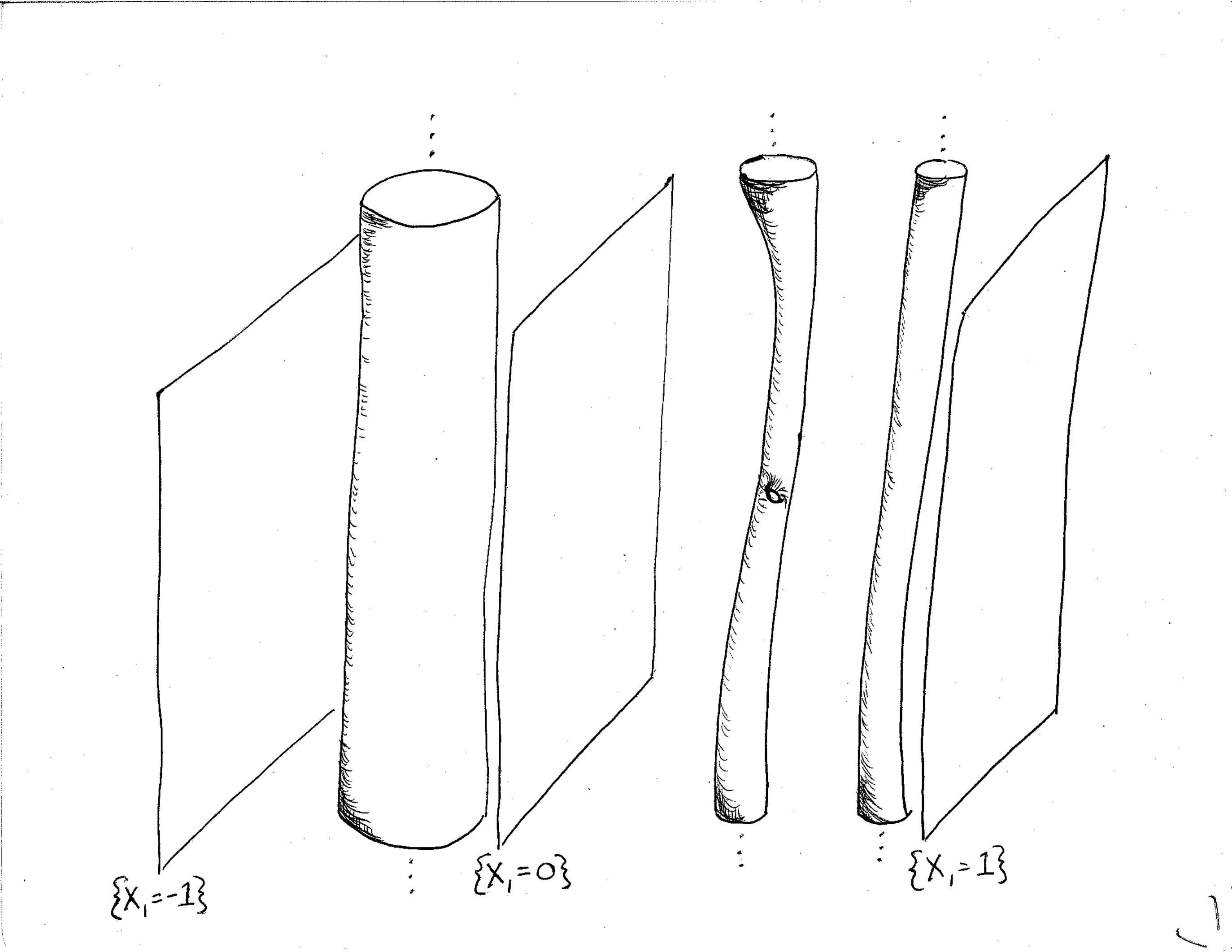}
	\caption{$\overline{W}_0 = W_+ \circ W_-$}
	\label{fig:W-split}
\end{center}
\end{figure}

Consider the $d$-dimensional $\mathcal{F}$-manifolds
\[
(\overline{W}_r,\overline{g}_r):= (\Phi^{-1}_1(\psi(W)-re_{x_{k-1}}),\Phi^*_1(\_ -r e_{x_{k-1}})^* \psi_* g),
\]
$r\in[0,3]$.
Define $(W_+,g_+)\in mor_{\mathsf{Cob}^\mathcal{F}_{d,(k,N)}}((\overline{0},\emptyset),(\overline{1},\emptyset))$ as
\[
(W_+,g_+):= (\overline{W}_0\cap pr^{-1}_k[0,\infty), {\overline{g}_0}_{\mid_{pr^{-1}_k[0,\infty)}})
\]
and $(W_-,g_-)\in mor_{\mathsf{Cob}^\mathcal{F}_{d,(k,N)}}((\overline{-1},\emptyset),(\overline{0},\emptyset))$ as
\[
(W_-,g_-):= (\overline{W}_0\cap pr^{-1}_k[0,\infty), {\overline{g}_0}_{\mid_{pr^{-1}_k[0,\infty)}});
\]
see Figures~\ref{fig:W},~\ref{fig:psi-W}, and~\ref{fig:W-split}.

There is a path in $mor\text{ }\mathsf{Cob}^\mathcal{F}_{d,(k,N)}$ from $(W,g)$ to $(W_+,g_+)$ as realized by scaling by a factor of $1/3$ about $x_1=1/2$ in the $x_k$-coordinate, translating by $2$ in the $x_{k-1}$-coordinate, then using the family $\Phi_{r^\prime}$, $r^\prime\in[1,\infty]$.

Because $\psi(W)$ has $x_{k-1}$-coordinate bounded above, the family $(\overline{W}_r,\overline{g}_r)$, $r\in[0,3]$, realizes a path in $mor$ from the composition $(\overline{W}_0,\overline{g}_0) = (W_+,g_+)\circ (W_-,g_-)$ to $\emptyset$.

Using the identification  $F(\nu,\emptyset)\cong D^\mathcal{F}_{d,(k-1,N)}$, $\nu= -1,0,1$, it is clear that  $F(\emptyset)\simeq id_{D^\mathcal{F}_{d,(k-1,N)}}$ and we have achieved a homotopy from 
\[
F(\emptyset)\xrightarrow{ F(W_-,g_-)\circ F(W_+,g_+)} F(\emptyset)
\]
to $id_{F(\emptyset)}$.   As promised, this demonstrates that  $F(W_-,g_-)$ is a left homotopy inverse to $F(W,g)$.  To show that $F(W_-,g_-) $ is also a right homotopy inverse  amounts to the above procedure performed on $(W_-,g_-)$ in place of $(W,g)$ while making use of the obvious fact that our rotation satisfies $R\circ R = id_{\mathbb{R}^k\times I^{N-k}}$.

\subsubsection*{Case $(\overline{0},M_0,g_0)\xrightarrow{(W,g)}(\overline{1},\emptyset)$}

\begin{figure}
\begin{center}
	\includegraphics[width=4in,viewport=8 30 710 600,clip=true]{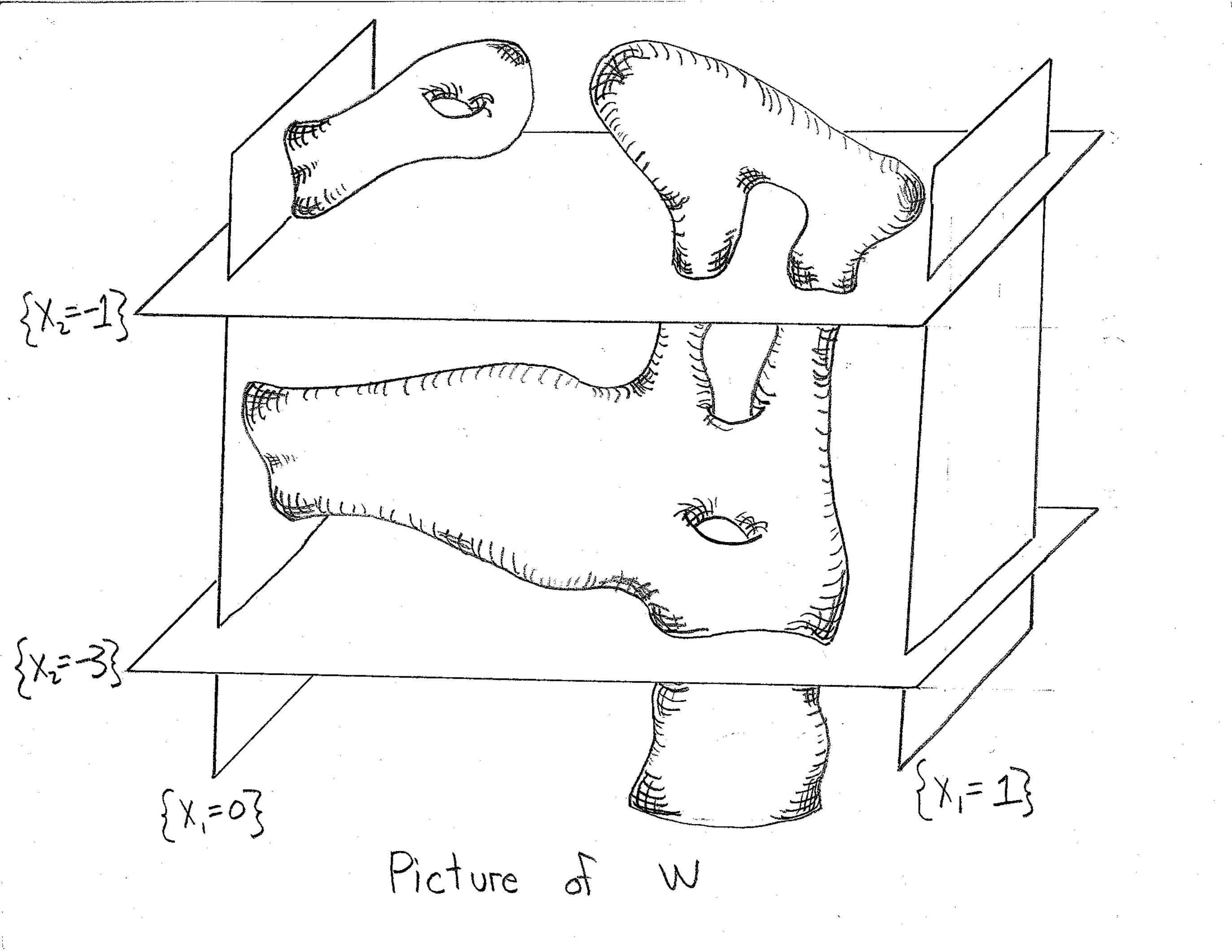}
	\caption{$W$}
	\label{fig:W-bdry}
\end{center}
\end{figure}

\begin{figure}
\begin{center}
	\includegraphics[width=4in,viewport=0 20 780 570,clip=true]{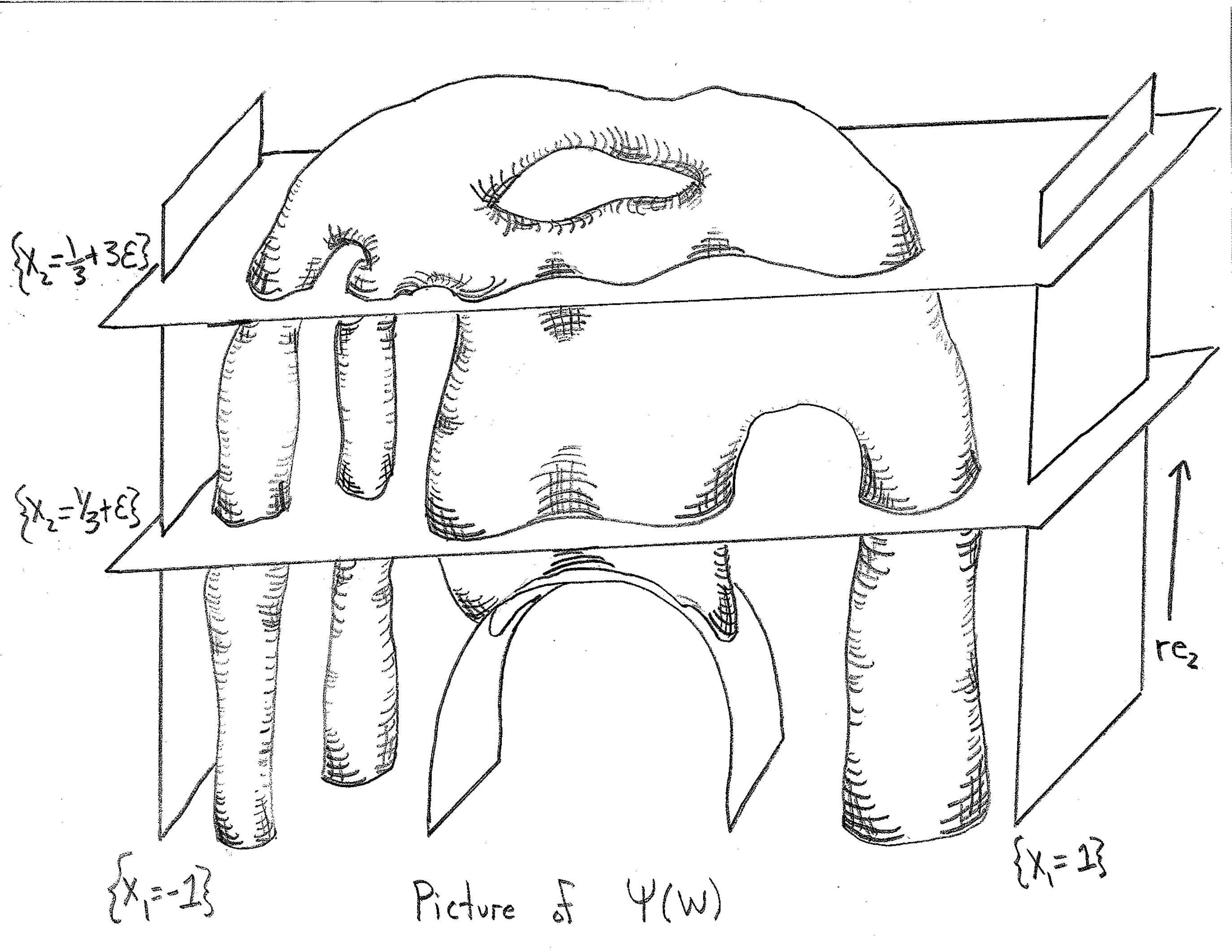}
	\caption{$\psi(W)$ and choice of $\epsilon>0$}
	\label{fig:psi-W-bdry}
\end{center}
\end{figure}

\begin{figure}
\begin{center}
	\includegraphics[width=3.2in,viewport=30 10 780 560,clip=true]{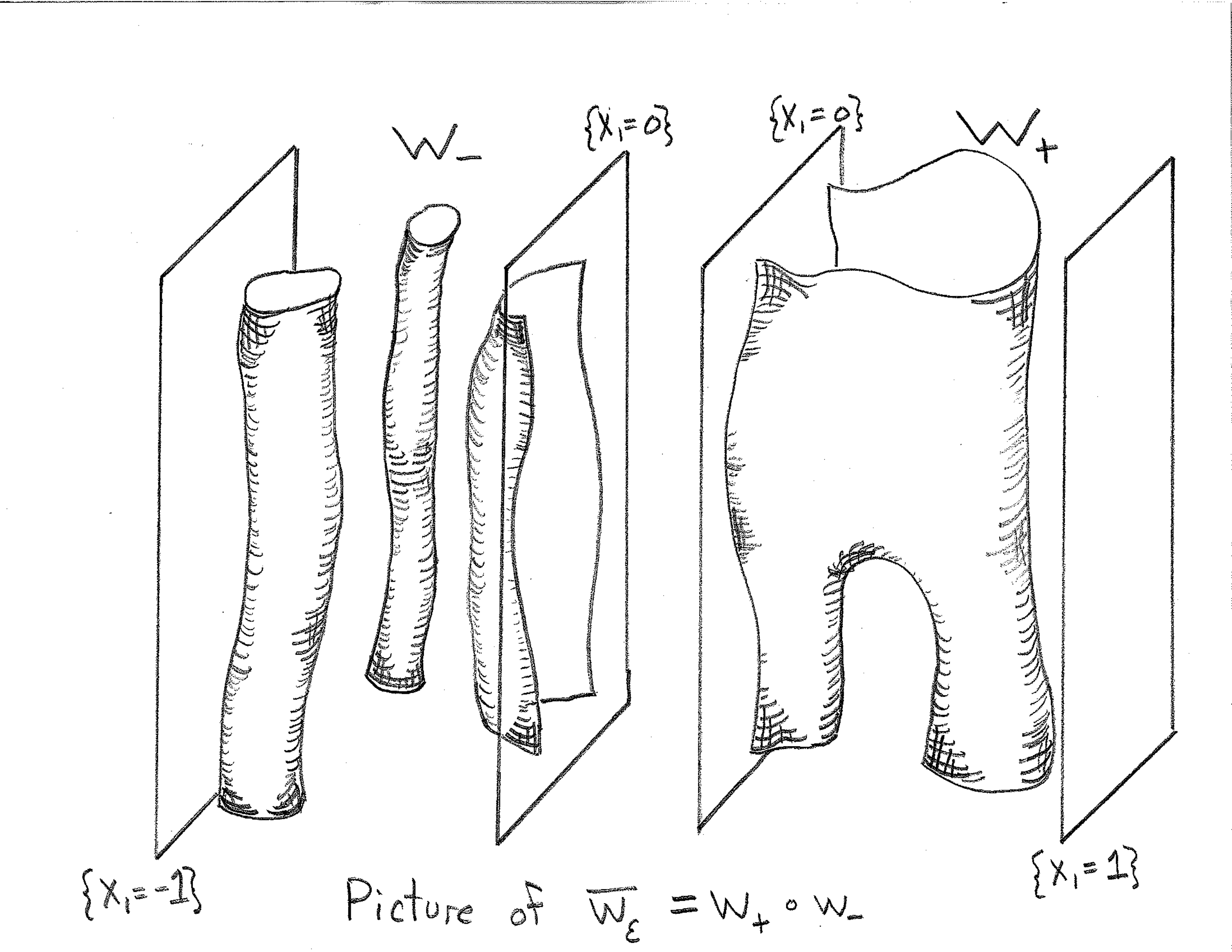}
	\caption{$\overline{W}_0 = W_+ \circ W_-$}
	\label{fig:W-bdry-split}
\end{center}
\end{figure}

Define 
\[
(W_1,g_1):= (\Phi^{-1}_1(W + 2e_{x_{k-1}}),\Phi^*_1(\_+2e_{x_{k-1}})^* g)\in mor_{\mathsf{Cob}^\mathcal{F}_{d,(k,N)}}((0,M^{\prime\prime},g^{\prime\prime}),(1,\emptyset)).
\]
Similar to the previous case, there is a path in $mor$ from $(W,g)$ to $(W_1,g_1)$ as illustrated by translating by $2$ in the $x_2$-coordinate then using the family $\Phi_{r^\prime}$, $r^\prime \in[1,\infty]$.

For $\epsilon>0$ consider the $d$-dimensional $\mathcal{F}$-manifolds
\[
(\overline{W}_r,\overline{g}_r):=
(\Phi^{-1}_\epsilon (\Psi(W)-(r+\epsilon+1/3)e_{x_{k-1}}), \Phi^*_\epsilon (\_ - (r+\epsilon+1/3)e_{x_{k-1}})^* \psi_* g),
\]
$r\in[\epsilon,3]$;
see Figures~\ref{fig:W-bdry} and~\ref{fig:psi-W-bdry}.   Because $W$ is collared, for $\epsilon>0$ small enough, $(\overline{W}_\epsilon,\overline{g}_\epsilon) \in mor_{\mathsf{Cob}^\mathcal{F}_{d,(k,N)}}((\overline{-1},\emptyset),(\overline{1},\emptyset))$ and 
$
pr_k:\overline{W}_\epsilon\rightarrow \mathbb{R}
$
has $0\in\mathbb{R}$ as a regular value.  Define $(W_+,g_+)\in mor_{\mathsf{Cob}^\mathcal{F}_{d,(k,N)}}((\overline{0},M^{\prime},g^{\prime}),(\overline{1},\emptyset))$ as
\[
(W_+,g_+):=(\overline{W}_\epsilon\cap pr^{-1}_k[0, \infty) , {\overline{g}_\epsilon}_{\mid_{pr^{-1}_k[0,\infty)}})
\]
and $(W_-,g_-)\in mor_{\mathsf{Cob}^\mathcal{F}_{d,(k,N)}}((\overline{-1},\emptyset),(\overline{0},M^{\prime},g^{\prime}))$ as
\[
(W_-,g_-):=(\overline{W}_\epsilon\cap pr^{-1}_k(-\infty, 0] , {\overline{g}_\epsilon}_{\mid_{pr^{-1}_k(-\infty,0]}});
\]
see Figure~\ref{fig:W-bdry-split}.

\begin{figure}
\begin{center}
	\includegraphics[width=3in,viewport=80 140 500 570,clip=true]{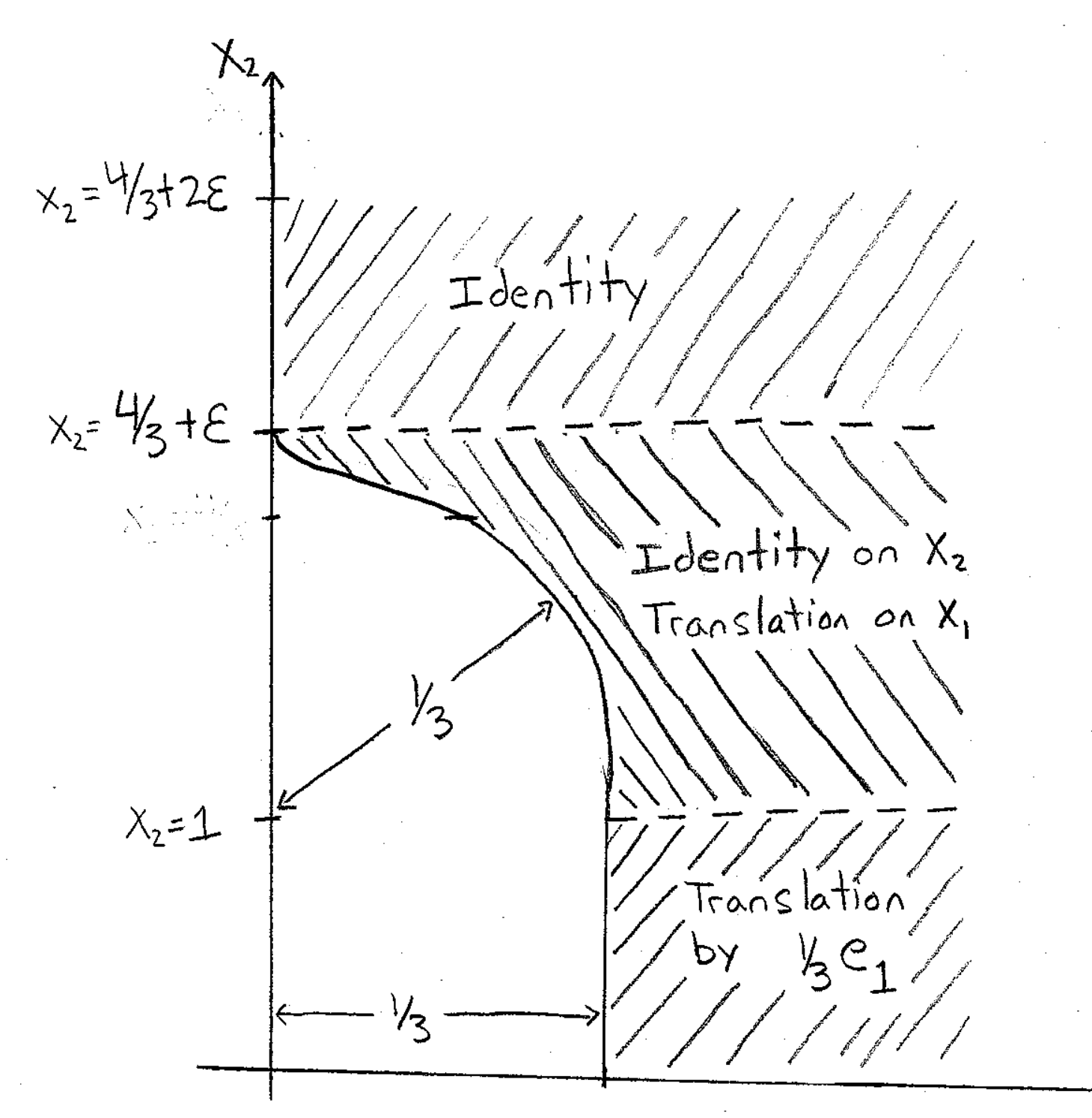}
	\caption{$\Theta$}
	\label{fig:theta}
\end{center}
\end{figure}

Choose a diffeomorphism $\Theta$ from $\mathbb{R}^{k-1}\times [0,1]\times I^{N-k}$ onto the region shaded in Figure~\ref{fig:theta} which is the identity on all but the $(x_{k-1},x_k)$-coordinates and is described on $\{(x_{k-1},x_k)\}$ as follows.  For $x_{k-1}\leq -1$,
\[
(x_{k-1},x_k)\mapsto (x_k,x_{k-1}/3 +1/3),
\]
and for $ x_k\geq 1$,
\[
(x_{k-1},x_k)\mapsto (x^\prime_{k-1},x^\prime_k)\text{\hspace{.3cm} with \hspace{.3cm}} 4/3+ \epsilon \leq x^\prime_{k-1} \leq 4/3 +3\epsilon;
\]
see Figure~\ref{fig:theta} for the rest of the description of $\Theta$.

Upon possibly re-choosing a smaller $\epsilon>0$, there is an ambient isotopy of diffeomorphisms of $\mathbb{R}^{k-1}\times [0,1]\times I^{N-k}$ restricting to an isotopy as collared $\mathcal{F}$-manifolds from the composite embedding 
\[
(W_+,g_+)\xrightarrow{(\_ + (2\epsilon+4/3)e_{x_{k-1}})\circ\Phi_\epsilon}(\psi(W),\psi_*g)\xrightarrow{\Theta^{-1}} (\Theta^{-1}\psi(W),\Theta^*\psi_*g)
\]
to a diffeomorphism of collared $\mathcal{F}$-manifolds
\[
(W_+,g_+)\xrightarrow{\cong} (\Theta^{-1}\psi(W),\Theta^*\psi_*g).
\]
In particular, there is a path in $mor$  between the two $\mathcal{F}$-manifolds above.
Because $\Theta^{-1}\circ \psi$ is translation by $2e_{x_{k-1}}$ on $\{x_{k-1}\leq -1\}$, 
\[
((\Theta\circ \Phi_1)^{-1}\psi(W),(\Theta\circ\Phi)^*\psi_*g)= (W_1,g_1)
\]
above.  It follows that there is a path in $mor$ from $(W,g)$ to $(W_+,g_+)$ and thus the maps $F(W,g)$ and  $F(W_+,g_+)$ are homotopic.

Evidently, because $\psi(W)$ has $x_{k-1}$-coordinate bounded above, the family $(\overline{W}_r,\overline{g}_r)$, $r\in[\epsilon,3]$, is a path of morphisms from the composition  $(\overline{W}_\epsilon,\overline{g}_\epsilon)=(W_+,g_+)\circ (W_-,g_-)$ to $\emptyset$.
Again, using the identification $F(\nu,\emptyset)\cong D^\mathcal{F}_{d,(k-1,N)}$, 
we have demonstrated that $F(W_-,g_-)$ is a left homotopy inverse to $F(W,g)$.  But 
\[
(\overline{-1} ,\emptyset)\xrightarrow{(W_-,g_-)} (\overline{0},M^\prime,g^\prime)
\]
is a morphism \textit{from} the empty $(d-1)$-$\mathcal{F}$-manifold.  So to show that $F(W_-,g_-)$ is a right homotopy inverse to $F(W_+,g_+)$ we pass to the remaining case.

\subsubsection*{Case $(0,\emptyset)\xrightarrow{(W,g)}(1,M_1,g_1)$}

\begin{figure}
\begin{center}
	\includegraphics[width=4in,viewport=0 72 792 540,clip=true]{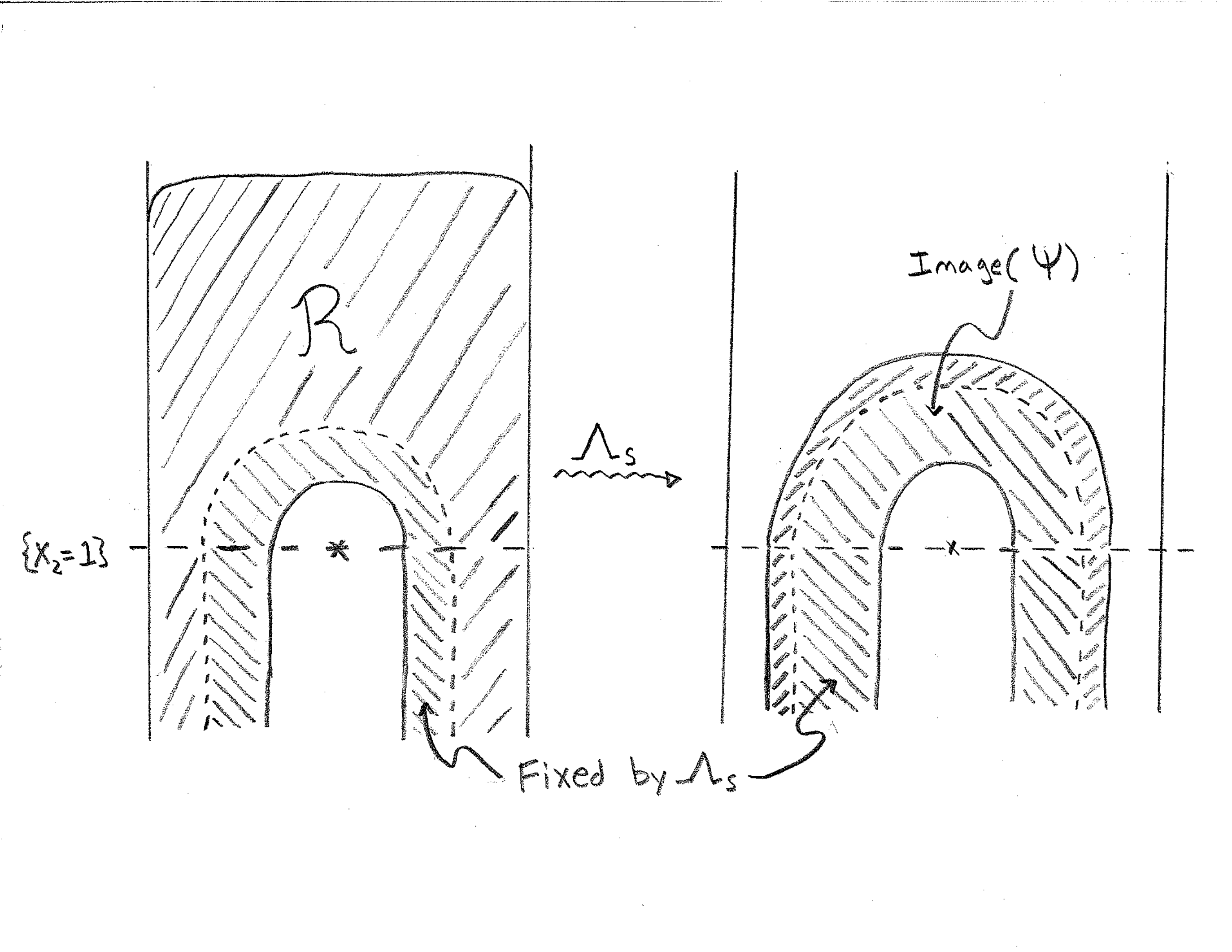}
	\caption{$\mathcal{R}$}
	\label{fig:R}
\end{center}
\end{figure}

Choose $\delta>0$ so that 
\[
W\cap pr^{-1}_k[1-\delta,1] = \mathbb{R}^{k-1}\times [1-\delta]\times M_1
\]
is a collar of $M_1=\partial W_1 \subset W_1$.  Let $\mathcal{R}\subset \mathbb{R}^{k-1}\times  [-1,1]\times I^{N-k}$ be the region in Figure~\ref{fig:R}.    There is an obvious smooth straight-line deformation 
\[
\Lambda_s:\mathcal{R}\rightarrow\mathcal{R},\hspace{.5cm}  s\in[0,1],
\]
of $\Lambda_0(\mathcal{R})=\mathcal{R}$ onto $\Lambda_1(\mathcal{R})=\psi(\mathbb{R}^{k-1}\times [1-\delta,1]\times I^{N-k})$
which is the identity near $\psi(\{x_k = 1-\delta\})$.

Consider the $d$-dimensional $\mathcal{F}$-manifold 
\[
(\overline{W}_r,\overline{g}_r):=(( \Lambda_1\circ\Phi_1)^{-1}(\psi(W)-re_{x_{k-1}}),(\Lambda_1\circ\Phi_1)^* (\_ - re_{x_{k-1}})^*\psi_* g),
\]
$r\in[0,3]$.
Define $(W_+,g_+) \in mor_{\mathsf{Cob}^\mathcal{F}_{d,(k,N)}}((0,\emptyset),(1,M^{\prime},g^{\prime}))$ as
\[
(W_+,g_+):=(\overline{W}_0\cap pr^{-1}_k[0, \infty),{\overline{g}_0}_{\mid_{pr^{-1}_k[0,\infty)}})
\]
and $(W_-,g_-)\in  mor_{\mathsf{Cob}^\mathcal{F}_{d,(k,N)}}((\overline{-1},R(M^{\prime}),R^*g^{\prime}),(\overline{0},\emptyset))$ as
\[
(W_-,g_-):=(\overline{W}_0\cap pr^{-1}_k(-\infty,0],{\overline{g}_0}_{\mid_{pr^{-1}_k(-\infty,0]}});
\]
see Figure~\ref{fig:W-bdry-split} for a similar picture.

There is a path in $mor$ from $(W,g)$ to $(W_+,g_+)$ as realized by scaling by a factor of $1/3$ about $x_k=1/2$ in the $x_k$-coordinate, translating by $2$ in the $x_{k-1}$-coordinate, using the family $\Phi_{r^\prime}$, $r^\prime\in[1,\infty]$, and finally using the isotopy $\Lambda_s$.  Thus the maps $F(W,g)$ and $F(W_+,g_+)$ are homotopic.

The family $(\overline{W}_r,\overline{g}_r)$, $r\in[0,3]$, realizes a path in  $mor\text{ }\mathsf{Cob}^\mathcal{F}_{d,(k,N)}$ from the composition $(\overline{W}_0,\overline{g}_0)=(W_-,g_-)\circ(W_+,g_+)$  to a product morphism $(\mathbb{R}^{k-1}\times [-1,1]\times M^{\prime\prime},g^{\prime\prime})$; see Figures~\ref{fig:psi-W-bdry} and~\ref{fig:W-bdry-split} for a similar picture.   There is a homeomorphism $F(\overline{-1},M_0,g_0)\cong F(\overline{1},M_0,g_0)$ seen by translation by $2e_{x_k}$.  Under this identification, it is clear that composition with the product $(\mathbb{R}^{k-1}\times [-1,1]\times M_0,g_0)$ is a homotopy equivalence.  We have thus demonstrated that $F(W_-,g_-)$ is a left homotopy inverse to $F(W,g)$.  From the previous case, it follows that in fact, $F(W_-,g_-)$ is also a \textit{right} homotopy inverse and so $F(W,g)$ is a homotopy equivalence.

\subsubsection*{General case $(\overline{0},M_0,g_0)\xrightarrow{(W,g)}(\overline{1},M_1,g_1)$}

Here it is merely pointed out that the procedures from the previous two cases can be done simultaneously.  An amalgamation of these two procedures proves that $F(W,g)$ is a weak homotopy equivalence in this general case.

\section{Interpretations of $B\mathsf{Cob}^\mathcal{F}_d$}\label{subsec: interpretations}

We interpret the set (group) $[X,B\mathsf{Cob}_d^\mathcal{F}]$ as concordance classes of certain families of $\mathcal{F}$-manifolds in an appropriate sense.  As such, the main theorem is a statement about families of $\mathcal{F}$-manifolds and each of the examples below can be, and perhaps should be, regarded in this way.  In this sense, from the main theorem, the topology of $\Omega^\infty MT\mathcal{F}(d)$ naturally reflects how $\mathcal{F}$-manifolds are organized among each other.

We will draw another interpretation of $B\mathsf{Cob}^\mathcal{F}_d$ as a sort of (group) completion of its morphism spaces $\mathcal{F}(W)//{Diff}(W)$ otherwise known as the moduli spaces of $\mathcal{F}$-structures on $W$.  From the main theorem, $\Omega^\infty MT\mathcal{F}(d)$ could, perhaps should, then be interpreted as a parametrizing space for `stable' characteristic classes of bundles with $d$-dimensional $\mathcal{F}$-manifold fibers.

\subsection{Stable families of manifolds}

Given  a locally finite open cover $\mathcal{U}=\{U_\alpha\}$ of a (paracompact) space $X$, construct a category $X_{\mathcal{U}}$ as follows.  Declare
\[
ob X_\mathcal{U} = \coprod_\alpha U_\alpha\text{ \hspace{1cm} and \hspace{1cm} } mor X_\mathcal{U} = \coprod_{\alpha,\beta} U_\alpha \cap U_\beta.
\]
The source and target maps of $X_\mathcal{U}$ are the determined by the inclusions 
\[
U_\alpha\xleftarrow{s}U_\alpha\cap U_\beta\xrightarrow{t} U_\beta,
\]
the identity map by $U_\alpha \xrightarrow{i} U_\alpha\cap U_\alpha$,  and composition by the inclusion 
\[
(U_\alpha\cap U_\beta) \times_{U_\beta} (U_\beta\cap U_\gamma) = U_\alpha\cap U_\beta \cap U_\gamma \xrightarrow{\circ} U_\alpha\cap U_\gamma.
\]

Let $\mathsf{C}$ be a topological category.
Define a  \textit{pair $(\mathcal{U},F)$ on $X$} to be the data of a locally finite open cover $\mathcal{U}$ of $X$  and a functor $F:X_\mathcal{U}\rightarrow \mathsf{C}$.  A \textit{concordance} from the pair $(\mathcal{U}_0,F_0)$ to the pair $(\mathcal{U}_1,F_1)$ is a pair $(\mathcal{V}, H)$ on $[0,1]\times X$ such that $\mathcal{V}$ restricts to $\mathcal{U}_i$ on $\{i\}\times X$ ($i=0,1$) and $H_{\mid_{\{i\}\times X}} = F_i$.

There is a bijection
\[
[X,B\mathsf{C}]\cong \{\text{pairs } (\mathcal{U},F)\text{ on $X$}\}/concordance,
\]
see~\cite{segal:classifying-spaces} for a first account and~\cite{weiss:classify} for a more recent viewpoint.  To supply some intuition toward a proof of this bijection, a pair $(\mathcal{U},F)$ on $X$ determines a map $X\rightarrow B\mathsf{C}$, well-defined up to homotopy, as follows.  The locally finiteness of the cover $\mathcal{U}$ ensures that the set $A_x:=\{\alpha\mid x\in U_\alpha\}$ is finite.  Regard the point $x$ in the $finite$-fold intersection $x\in \cap_{\alpha\in A_x} U_\alpha$ as a simplex in the simplicial nerve $N_\bullet X_\mathcal{U}$.  Choose a partition of unity $\{\lambda_\alpha\}$ subordinate to $\mathcal{U}$.  Declare $x\mapsto (F(x),\Sigma_{A_x} \lambda_\alpha(x))\in B\mathsf{C}$.  The homotopy class of this map is independent of choices.  Conversely, there is a canonical open cover of $B\mathsf{C}$ obtained from its simplicial definition; pull back and a Yoneda-type argument finish the bijection.

Applying this construction to the cobordism category $\mathsf{Cob}^\mathcal{F}_d$ establishes the interpretation of $[X,B\mathsf{Cob}^\mathcal{F}_d]$ as concordance classes of \textit{stable} families of $d$-dimensional $\mathcal{F}$-manifolds over $X$.  
More precisely, upon making some choices and using the bijection above, associate to a map $X\rightarrow B\mathsf{Cob}^\mathcal{F}_d$  a pair $(\mathcal{U},F)$ on $X$.  From this pair build a space $E\subset X\times\mathbb{R}\times\mathbb{R}^\infty$ by declaring over $x\in X$, $E_x = F_{A_x}(x)\in mor\text{ }\mathsf{Cob}^\mathcal{F}_d$.  (Here, regard $F_{A_x}(x)\in N_\bullet\mathsf{Cob}^\mathcal{F}_d$ as a morphism (a cobordism with $\mathcal{F}$-structure) via composition in $\mathsf{Cob}^\mathcal{F}_d$).  The collection $\{E_x\}$ fits together into a topological space $E\subset\mathbb{R}\times\mathbb{R}^{\infty}$.  By construction, the map $E\xrightarrow{pr_X} X$ has $d$-dimensional smooth $\mathcal{F}$-manifolds as fibers.  Regard $E\subset X\times\mathbb{R}\times\mathbb{R}^\infty$ as a \textit{stable} family of $\mathcal{F}$-manifolds in this way.

\begin{terminology}
The phrase ``family'', or ``stable family'', will be reserved for the notion described above as distinguished from notion of a bundle when the context is clear.
\end{terminology}

When $X$ is a manifold, the ideas of~\cite{galatius-madsen-tillmann-weiss:cobordism} and~\cite{madsen-weiss:mumford} provide a better description of $[X,B\mathsf{Cob}_d^\mathcal{F}]$ as concordance classes of codimension $-d$ submersions over $X$ having vertical $\mathcal{F}$-structure.  It is this point of view that drives the methods in the proof of the main theorem, namely, the quasi-topology of the sheaf $\Psi^\mathcal{F}_d$.
In the end, due to the inclusion
\[
\mathcal{F}(M^{d-1})//{Diff}(M)\hookrightarrow ob\text{ }\mathsf{Cob}_d^\mathcal{F}\hookrightarrow B\mathsf{Cob}^\mathcal{F}_d,
\]
given a \textit{fiber bundle} $E\xrightarrow{pr_X} X$ having fibers $(d-1)$-dimensional $\mathcal{F}$-manifolds, there is a homotopy class of a map $X\rightarrow B\mathsf{Cob}_d^\mathcal{F}$.  This map classifies the stable family $E\times\mathbb{R}\rightarrow X$ of $d$-dimensional $\mathcal{F}$-manifolds as in the discussion above.  Thus $B\mathsf{Cob}^\mathcal{F}_d$ classifies families more general than smooth fiber bundles.  The section to follow is helpful in this direction.

\subsection{The cobordism category of fibered manifolds}

Fix an embedded compact $r$-manifold $X^r\subset\mathbb{R}^\infty$.  Consider a category whose objects are roughly smooth fiber bundles $M^{r+d-1}\rightarrow X^r$ with compact $(d-1)$-dimensional fiber and whose morphisms are concordances of such.  More precisely, define $\mathsf{Cob}_d(X)$ as follows.
Declare
\[
ob\text{ }\mathsf{Cob}_d(X) =\{(a, M^{r+d-1})\}
\]
where $ a\in Map(X,\mathbb{R})$ and $M\subset X\times\mathbb{R}^{N-1}$ with $M\xrightarrow{pr_X} X \text{ a smooth fiber bundle}$;
and 
\[
mor\text{ }\mathsf{Cob}_d(X) = \{identities\}\amalg \{(a,W^{r+d})\}
\]
where $ a\in Map(X,\mathbb{R}^2_+)$ and $W\subset [0,1]\times X\times\mathbb{R}^{N-1}$ with $ W\xrightarrow{pr_X} [0,1]\times X$ a smooth fiber bundle which is \textit{cylindrical} over a neighborhood of $\{\nu\}\times X$, $\nu=0,1$.  The source, target, and identity maps are apparent.  Composition is given by concatenation of concordances along with reparametrizing $[0,1]$ using the maps $a\in Map(X,\mathbb{R}^2_+)$.

This category $\mathsf{Cob}_d(X)$ is a topological category in the same way that $\mathsf{Cob}_d$ is a topological category.  It is possible to describe the weak homotopy type of $\mathsf{Cob}_d(X)$ by reconsidering its object and morphism spaces.  Recall the notion of a smooth map $X\rightarrow B_N(M)$ from \textsection\ref{subsec: smooth families}.  Let $Map^{sm}(X,B_N(M))$ denote the space of smooth maps; it is a subset of the space of continuous maps $Map(X,B_N(M))$ and is topologized accordingly.  Observe that
\[
ob\text{ }\mathsf{Cob}_d(X) \cong Map^{sm}(X,ob\text{ }\mathsf{Cob}_d)
\]
and
\[
mor\text{ }\mathsf{Cob}_d(X)\cong Map^{sm}(X,mor\text{ }\mathsf{Cob}_d).
\]
Moreover, the structure maps of $\mathsf{Cob}_d(X)$ are induced from the structure maps of $\mathsf{Cob}_d$.  This is close to the general situation of obtaining from a (topological) category $\mathsf{C}$ the (topological) category $Map(X,\mathsf{C})$.

It has already been noted (see \textsection\ref{subsec: smooth families}) that the inclusion 
\[
Map^{sm}(X,B_N(M))\hookrightarrow Map(X,B_N(X))
\]
is a weak homotopy equivalence.  It follows that the map of simplicial spaces
\[
N_\bullet \mathsf{Cob}_d(X)   \hookrightarrow   Map(X,N_\bullet\mathsf{Cob}_d)
\]
is level-wise a weak homotopy equivalence and therefore 
\[
B\mathsf{Cob}_d(X) \xrightarrow{\simeq} BMap(X, N_\bullet \mathsf{Cob}_d) 
\]
is a weak homotopy equivalence.  The same holds true if we include $\mathcal{F}$-structure into the discussion.

There is an obvious map $BMap(X,N_\bullet \mathsf{Cob}_d)\rightarrow Map(X,B\mathsf{Cob}_d)$.  This map is far from an equivalence; the homotopy fiber measures the obstruction for a submersion over the compact manifold $X$ to be concordant to a fiber bundle over $X$.  Indeed, a point in $BMap(X,N_\bullet \mathsf{Cob}_d)$ is represented, up to homotopy, by a point in $N_\bullet \mathsf{Cob}_d(X)$.  In particular, by a bundle $W\rightarrow X$ over $X$.

\subsection{Group completion}\label{subsubsec: group completion}

The details for making this subsection precise can be found in Rezk's paper on Segal spaces (\cite{rezk:segal-spaces}) or Lurie's survey of topological field theories (\cite{hopkins-lurie:tfts}, \textsection $2.1$).  The purpose of this subsection is to justify regarding the classifying space construction $B\mathsf{Cob}^\mathcal{F}_d$ as a group completion construction.  Informally, a group completion of the moduli spaces of $\mathcal{F}$-structures on $d$-manifolds.

Roughly, a pre-$\infty$-category is a simplicial space $\mathcal{X} = \{X_k\}$ for which the diagram
\[
\xymatrix{
X_{k+l} \ar[d] \ar[r]
&
X_l  \ar[d]
\\
X_k \ar[r]
&
X_0
}
\]
is homotopy Cartesian; that is, the universal map $X_{k+l} \rightarrow X_k\times_{X_0} X_l$ is a weak homotopy equivalence.  The left vertical map is that induced by the simplicial map $\{0,...,l\} \hookrightarrow \{0,...,l+k\}$ as inclusion of the first $l$ letters while the upper horizontal map is that induced by the inclusion of the last $k$ letters.  A functor between pre-$\infty$-categories is a morphism of simplicial spaces.  A functor between pre-$\infty$-categories is an \textit{equivalence} if it is a weak homotopy equivalence on each space of $k$-simplicies. The nerve functor gives an embedding of the category of topological categories into the category of pre-$\infty$-categories (actually landing in $\infty$-categories; see below).  We thus regard the $0$-simplicies of a pre-$\infty$-category as the objects, the $1$-simplicies as the morphisms, and the $k$-simplicies ($k\geq 2$) as higher morphisms.

We say that a $1$-simplex $f$ of  a pre-$\infty$-category is \textit{invertible} if there is a $2$-simplex $\sigma\in X_2$ with $\partial_2 \sigma = f$ and with $\partial_1 \sigma$ in the path component of $s_0(X_0)$ in $X_1$.  Notice that if $\mathcal{X} = N_\bullet\mathsf{C}$ for some topological category $\mathsf{C}$, then $f$ is invertible in this sense if $f$ is invertible in the ordinary sense.  Write $Z\subset X_1$ for the collection of $1$-simplicies.  By definition the degeneracy map $s_0:X_0 \rightarrow X_1$ factors through $Z$.  An $\infty$-category is a pre-$\infty$-category such that $s_0:X_0\rightarrow Z$ is a weak homotopy equivalence.

We say that an $\infty$-category is an $\infty$-groupoid if every morphism is invertible.  Just as $\Omega X$ is the prototypical example of a group-like $A_\infty$-space, the prototypical example of an $\infty$-groupoid is, for $X$ a (reasonable) space, the simplicial space $\mathcal{P}X$ described by 
\[
\mathcal{P}X_k = Map(\Delta^k,X).  
\]
This is called the \textit{fundamental groupoid of $X$}, or the \textit{path category on $X$}.  For $\mathcal{X}$ an $\infty$-category with classifying space $B\mathcal{X}$, there is a canonical functor 
\[
\mathcal{X}\rightarrow \mathcal{P}B\mathcal{X}.  
\]
This functor is initial with respect to functors $\mathcal{X}\rightarrow \mathcal{G}$ for $\mathcal{G}$ an $\infty$-groupoid.  In this way we regard $B\mathcal{X}$, actually the arrow $\mathcal{X}\rightarrow \mathcal{P}B\mathcal{X}$, as the groupoid completion of $\mathcal{X}$.  This universal property determines the weak homotopy type of $B\mathcal{X}$.  This construction is familiar in the case when $\mathcal{X} = N_\bullet M$ where $M$ is a topological monoid.  In this case, the groupoid completion of $\mathcal{X}$ recovers $M\rightarrow \Omega BM$ as expected.

Rephrase this notion of groupoid completion in a particular instance and more informally as follows.  Consider a topological category $\mathsf{C}$.  Let $Y$ be a space.  Suppose given a collection of maps, 
\[
mor_\mathsf{C}(o_0, o_1)\rightarrow Map((I,\{0\},\{1\}),(Y,\{y_0\},\{y_1\})) \simeq \Omega Y,
\]
for each pair $o_0,o_1\in\mathsf{C}$ such that composition in $\mathsf{C}$ agrees with with concatenation of paths in $Y$ in an appropriate sense.  Then there exists a map
\[
\Omega B\mathsf{C}\rightarrow \Omega Y
\]
which is unique up to homotopy.
In this way think of $\Omega B\mathsf{C}$ as parametrizing `stable' characteristic classes of the morphism spaces $mor_\mathsf{C}(o_0,o_1)$.  The word `stable' is appropriate because the right hand side of the map $mor_\mathsf{C}(o_0,o_1)\rightarrow  \mathcal{P}_{o_0,o_1} B\mathsf{C} \simeq \Omega B\mathsf{C}$ is independent of $(o_0,o_1)$.  Thus, the topology of $ \Omega B\mathsf{C}$ is a sort of amalgamation of the topology of the morphism spaces of $\mathsf{C}$.

Now take $\mathsf{C}=\mathsf{Cob}^\mathcal{F}_d$.  A component of the morphism space 
\[
mor_{\mathsf{Cob}^\mathcal{F}_d}((M_0,g_0),(M_1,g_1))
\]
is homotopy equivalent to the moduli space 
\[
\mathcal{F}(W;\partial)//{Diff}(W;\partial)
\]
of $\mathcal{F}$-structures on $W =(M_0\xrightarrow{W} M_1)$ fixed on the boundary.  From the discussion above,
$\Omega B\mathsf{Cob}^\mathcal{F}_d\simeq \Omega^{\infty}MT\mathcal{F}(d)$
parametrizes `stable' characteristic classes of moduli spaces of $\mathcal{F}$-structures.  It is in this sense that $ \Omega B\mathsf{Cob}^\mathcal{F}_d$ is the `group completion' of the moduli spaces $\mathcal{F}(W)//{Diff}(W)$ of $\mathcal{F}$-structures.

\begin{remark}
As an important example, take $d=2$ and $\mathcal{F}$ to be the sheaf of complex structures; see \textsection\ref{subsec: sheaf of complex stctrs} below.  In this case these ideas are consistent with the Mumford conjecture and its proof (see~\cite{madsen-weiss:mumford}).  Namely, the cohomology ring of $\Omega B\mathsf{Cob}^\mathcal{F}_d$ \textit{is} the cohomology of $BDiff^{or}(W^2)$ in the stable range from Harer stability (\cite{harer:stability}). 

On the other hand, recently Ebert (\cite{ebert:vanishing}) proved, for $W^3$ closed and oriented, $BDiff^{or}(W)\rightarrow \Omega B\mathsf{Cob}^{\mathcal{O}r}_3$ is the zero map on rational cohomology.  So it is only formal that $\Omega B\mathsf{Cob}^\mathcal{F}_d$ parametrizes stable characteristic classes.  
\end{remark}

The main theorem then informs us that there are no \textit{geometric} stable characteristic classes of moduli spaces which are not obtained from tangent data.  Indeed, the main theorem states that the morphism $\mathcal{F}\rightarrow\tau\mathcal{F}$ induces an equivalence on the groupoid completions of the resulting cobordism categories.

\section{Fiberwise equivariant sheaves}\label{subsec: fiberwise equivariant sheaves}

Consider a fibration $B_\theta \xrightarrow{\theta} BO(d)$.  
Take $\mathcal{F}_\theta: {\mathsf{Emb}_d}^{op}  \rightarrow   \mathsf{QTop}$ to be the associated \textit{fiberwise} equivariant sheaf of $\theta$-structures.  Recall from \textsection\ref{subsubsec: alternative description} that $\mathcal{F}$ is described by
\[
\mathcal{F}_\theta(W) = \{l\in Map(W,B_\theta)\mid \theta\circ l = \tau_W\}.
\]
This situation gives rise to a large class of examples of equivariant sheaves; only a few of which will be mentioned below with the rest left to the reader's imagination.  Such examples are understood in~\cite{galatius-madsen-tillmann-weiss:cobordism} and from our main theorem are the prototypical example to which all others compare.  Indeed, for $\mathcal{F}$ any equivariant sheaf, the associated $\tau\mathcal{F}$ is a fiberwise equivariant sheaf.

When $\theta:Y\times BO(d)\rightarrow BO(d)$ is a trivial fibration, an $\mathcal{F}_\theta$-structure on $W^d$ is simply a map $W\rightarrow Y$.  In this case we denote $\mathcal{F}_\theta$ as $map_Y$.  In light of the main theorem it is not difficult to verify
\[
B\mathsf{Cob}^{map_Y}_d \simeq \Omega^{\infty-1}MTO(d)\wedge Y_+.
\]
More interesting examples come from a representation $G\xrightarrow{\rho} O(d)$.  In this case take $B_\theta = EG\times_\rho Fr_d = BG\xrightarrow{\theta = B\rho} BO(d)$.  As a special case, for $G=U(d/2)$, $\mathcal{F}_\theta(W):=\mathcal{J}(W)$ is the space of almost-complex structures on $W$ and 
\[
B\mathsf{Cob}^\mathcal{J}_d \simeq \Omega^{\infty-1}MTU(d/2).
\]
Other examples include $G=SO(d)$ yielding the orientation sheaf, denoted $\mathcal{O}r$, and $G=Spin(d)$ to yield the sheaf of spin structures.

\begin{terminology}
A \textit{linear} $\mathcal{F}_\theta$-structure on a vector space $V^d\in BO(d)$ is a point in $\theta^{-1}(V)$, the fiber of the fibration $B_\theta\rightarrow BO(d)$ over $V$.  This is to be separate, though related, to the already established notion of a  $\mathcal{F}_\theta$-structure on $V\subset\mathbb{R}^\infty$ regarded as a submanifold (forgetting the vector space structure).  The difference is that an $\mathcal{F}_\theta$-structure on $V$ is a continuous choice of \textit{linear} $\mathcal{F}_\theta$-structures on the vector spaces $T_v V\cong V$ for each $v\in V$.  
\end{terminology}

\section{Immersions and embeddings}

\subsection{Immersions}

Fix a smooth embedded $n$-manifold $Y^n\subset\mathbb{R}^\infty$ and define the equivariant sheaf $imm_Y$  by 
\[
imm_Y(M)=Imm(M,Y),
\]
the space of immersions of $M$ into $Y$, topologized with the $C^\infty$ Whitney topology.  It is easy to verify that $imm_Y$ is indeed an equivariant sheaf.

Because $Y\subset\mathbb{R}^\infty$, there is a canonical metric on $Y$ and we can thus assume $\tau_Y$ has structure group in $O(d)$.  Write $P_{\tau_Y}\rightarrow Y$ for the  principal $O(n)$-bundle associated to $\tau_Y$.  There is an action of $O(n)$ on the Grassmann $Gr_{d,n}$ of $d$-planes in $\mathbb{R}^n$ by acting on the ambient $\mathbb{R}^n$.  There results a bundle
\[
P_{\tau_Y}\times_{O(n)} Gr_{d,n}\rightarrow Y
\]
which we denote by $Gr_d(\tau_Y)$ having total space denoted by $Gr_d(TY)$.  There is a canonical map $Gr_d(TY)\rightarrow BO(d)$ given by 
\[
(V\subset T_yY)\mapsto (V\subset T_y Y \subset T_0\mathbb{R}^\infty \cong \mathbb{R}^\infty).
\]
Use superscript notation to denote Thom spectrum functor.

\begin{theorem}
There is a weak homotopy equivalence
\[
B\mathsf{Cob}^{imm_Y}_d\simeq \Omega^{\infty-1}(Gr_d(TY))^{-\gamma_d}.
\]
\end{theorem}

\begin{proof}
From the main theorem, we need only describe the fibration $p_\mathcal{F}: B{imm_Y} \rightarrow BO(d)$.  By definition, 
\[
B{imm_Y}=\{(V^d,g)\mid g\in Stalk_0({imm_Y}_{\mid_V})\}.
\]

Speaking generally for a moment, given two vector bundles $\alpha$ and $ \beta$ over possibly different base spaces, define $Inj(\alpha,\beta)$ to be the space of bundle morphisms $\alpha\rightarrow \beta$ which are fiberwise injective; this space is topologized using a compact-open topology.
Regard $V\in BO(d)$ as a vector bundle over a point $*$.  As such, a point in $Inj(V,\tau_Y)$ is the data of a point $y\in Y$ and an injection of vector spaces $V\hookrightarrow T_yY$.

There is an inclusion 
\[
Stalk_0({imm_Y}_{\mid_V})   \hookrightarrow \{(y,f)\mid f:V\hookrightarrow T_y Y\}   = Inj(V,\tau_Y).
\]
The inverse function theorem informs us that every immersion is  locally standard.  It follows without difficulty that the inclusion above is a homotopy equivalence.  There is an action of $O(d)$ on $Inj(\mathbb{R}^d,\tau_Y)$ by pre-composition.
The infinite Stiefel space $Fr_d$ of orthogonal $d$-frames in $\mathbb{R}^\infty$ with its $O(d)$-action as a model for $EO(d)$.  So
\[
B{imm_Y} \simeq Inj(\mathbb{R}^d,\tau_Y)//O(d):= Fr_d\times_{O(d)} Inj(\mathbb{R}^d,\tau_Y) 
\]
over $BO(d)$.  This is what one would expect from the theory of immersions via h-principles.

There is the apparent 
quotient map $Inj(\mathbb{R}^d,\tau_Y)\rightarrow Gr_d(\tau_Y)$.  Because this action of $O(d)$ on $Inj(\mathbb{R}^d,\tau_Y)$ is free and $Fr_d$ is contractible, the projection 
\[
B{imm_Y}\xrightarrow{\simeq} Gr_d(TY)
\]
is a homotopy equivalence
over $BO(d)$.  The result follows from the main theorem.  
\end{proof}

Let $Simm_Y$ be the equivariant sheaf of \textit{stable} immersions defined as the colimit
\[
Simm_Y:= colim_k \text{ }imm_{\mathbb{R}^k\times Y};
\]
the colimit is over the directed system induced from the inclusions $\{0\}\times\mathbb{R}^k\times Y\hookrightarrow\mathbb{R}^{k+1}\times Y$.

The statement of the above theorem is more explicit when $\tau_Y$ is trivial.  
Indeed, in this case $Gr_d(\tau_Y) = (Y\times Gr_{d,n}\rightarrow Y)$ and if follows that
\[
B\mathsf{Cob}^{imm_Y}_d\simeq \Omega^{\infty-1} (Th(\gamma^\perp_{d,n})\wedge Y_+).
\]  
In this parallelizable situation then, 
\[
B\mathsf{Cob}_d^{Simm_Y}\simeq \Omega^{\infty-1}(MTO(d)\wedge Y_+)\simeq B\mathsf{Cob}_d^{map_Y}.
\]
Recognize the space in the middle as one whose homotopy groups are the $MTO(d)$-homology of $Y$.

\subsection{Isometric immersions}

Let $Y=(Y^n,g_Y)$ be Riemannian $n$-manifold with $n\geq d$.  Let $met$ be the equivariant sheaf of Riemannian metrics, $met(W^d) = \{\text{Riemannian metrics on }W\}$, topologized as a subspace of the space of $2$-tensors on $W$.
Consider the equivariant sheaf $
iso_Y:{\mathsf{Emb}_d}^{op}\rightarrow \mathsf{Top}
$
described by 
\[
iso_Y(W)=\{(g,f)\mid f\in imm_Y(W)\text{ such that } f^*g_Y = g \in met(W)\}.  
\]
The obvious forgetful morphism $iso_Y\rightarrow imm_Y$ is a homeomorphism since any such $g\in met(W)$ is canonically determined by the immersion $f$.  It follows from the above section that 
\[
B\mathsf{Cob}^{iso_Y}_d\simeq \Omega^{\infty-1}(Gr_d(TY))^{-\gamma_d}.
\]

\subsection{The embedded cobordism category}

The author learned the content of the following section from Randal-Williams in his Oxford transfer thesis (see~\cite{randal-williams:embedding}).

For $Y^n$ a smooth manifold as before, let 
\[
emb_Y:{\mathsf{Emb}_d}^{op}\rightarrow \mathsf{Top}
\]
be the functor (presheaf) described by $emb_Y(M)= Emb(M,Y)$, the space of embeddings of $M$ into $Y$ topologized with the $C^\infty$ Whitney topology.  This functor $emb_Y$ is \textit{not} a sheaf.  The sheafification of $emb_Y$ is the sheaf $imm_Y$ above.

Consider the cobordism category $\mathsf{Cob}^{emb}_d(Y)$ whose objects are  closed $(d-1)$-submanifolds of $Y$ and whose morphisms are compact collared $d$-submanifolds of $[a_0,a_1]\times Y$.  More precisely, as a topological category
\[
ob\text{ }\mathsf{Cob}^{emb}_d(Y) = \mathbb{R}^\delta\times(\coprod_{[M^{d-1}]} Emb(M,Y)/{Diff}(M))
\]
and 
\[
mor\text{ }\mathsf{Cob}^{emb}_d(Y) =\{identities\}\amalg (\mathbb{R}^2_+)^\delta\times (\coprod_{[W^{d}]} Emb(W,[0,1]\times Y)/{Diff}(W)).
\]
This is \textit{not} an example of a cobordism category with $\mathcal{F}$-structure for an equivariant sheaf $\mathcal{F}$.  Nevertheless, it is possible to identify its weak homotopy type.

Recall the $\mathsf{Top}$-valued sheaf $\Psi_d$ from~\S\ref{subsec: the sheaf}.  Recall the construction of the sheaf of sections
\[
\Gamma(\Psi_d(\tau_Y)_\mid)
\]
from the beginning of~\S\ref{subsec: trad sheaves of sections}.  The h-principle comparison morphism results in a map of spaces (of global sections)
\[
\Psi_d(Y)\xrightarrow{exp^*} \Gamma(\Psi_d(\tau_Y)).
\]
Labeling this map as $exp^*$ is meant to remind the reader of the map's construction.

More generally, for $U\subset Y$ open with the closure $\overline{U}\subset Y$, define 
\[
\Psi_d(Y,U):=\{W^d\in\Psi_d(Y)\mid W\subset U\}.
\]
Similarly, define $\Gamma_{U}(\Psi_d^{fib}(Y))$ to be those sections which send $Y\setminus U$ to the base point $\emptyset\in \Psi_d$.  The map $exp^*$ restricts to a map
\[
\Psi_d(Y,U)\xrightarrow{exp^*} \Gamma_{U}(\Psi_d(\tau_Y)).
\]
Gromov's h-principle says that this map $exp^*$ is a weak homotopy equivalence provided the sheaf $\Psi_d$ is \textit{micro-flexible}; see~\cite{gromov:h-principle}.  Randal-Williams~\cite{randal-williams:embedding} showed that $\Psi_d$ is indeed micro-flexible and this theorem of Gromov's thus applies.

Upon replacing $Y$ by $\mathbb{R}\times Y$ and for $U=\mathbb{R}\times cpt(Y)$ with $cpt(Y)$ the interior of a compact core of $Y$, there is the following

\begin{theorem}
There is a weak homotopy equivalence
\[
B\mathsf{Cob}^{emb}_d(Y) \simeq \Psi_d(\mathbb{R}\times Y,\mathbb{R}\times cpt(Y)).
\]
\end{theorem}

The proof of the above theorem is the same as the more specialized case when $Y=\mathbb{R}^{N-1}$ and $cpt(Y)=(0,1)^{N-1}$ which was done in \textsection\ref{sec: proof of main theorem} of this paper.  We conclude with the theorem of Randal-Williams; see~\cite{randal-williams:embedding} for details.

\begin{theorem}[Randall-Williams]\label{thm: embedding cob cat}
There is a weak homotopy equivalence
\[
B\mathsf{Cob}^{emb}_d(Y)\simeq \Gamma_{cpt(Y)}(\Psi_d(\epsilon^1\oplus \tau_Y)).
\]
\end{theorem}

To make this more explicit still, the $O(n)$-equivariant map $Th({\gamma^d_n}^\perp)\xrightarrow{\simeq}\Psi_d(\mathbb{R}^n)$ is a weak homotopy equivalence from Lemma~\ref{lem: thom space}.  For $P_{\tau_Y}$ the associated principal $O(n)$-bundle for $\tau_Y$, we arrive at the homotopy equivalence of bundles over $BO(d)$
\[
Th(\gamma^\perp_{d,{\tau_Y}}):=(P_{\tau_Y}\times_{O(n)} Th(\gamma^\perp_{d,n})\rightarrow Y) \xrightarrow{\simeq} \Psi_d(\tau_Y).
\]
There is the resulting weak homotopy equivalence
\[
B\mathsf{Cob}^{emb}_d(Y) \simeq \Gamma_{cpt(Y)}(Th(\gamma^\perp_{d,{\epsilon^1\oplus\tau_Y}}))
\]
which is most explicit when $\tau_Y$ is stably trivial.  In this case, the right hand side of the above weak homotopy equivalence is homotopy equivalent to $Map_{cpt(Y)} (Y,Th(\gamma^\perp_{d,n+1}))$.

\subsection{The stably embedded cobordism category}

Upon replacing $Y$ with $\mathbb{R}^k\times Y$, Theorem~\ref{thm: embedding cob cat} takes the form 
\[
B\mathsf{Cob}^{emb_{\mathbb{R}^k\times Y}}_d\simeq \Omega^k  \Gamma_{cpt(Y)}(Th(\gamma^\perp_{d,\epsilon^1\oplus\tau_Y})).
\]
For each $t\in\mathbb{R}$, the inclusion $\{t\}\times \mathbb{R}^k\times Y\hookrightarrow \mathbb{R}\times \mathbb{R}^k\times Y$ induces a functor of cobordism categories 
\[
\mathsf{Cob}_d^{emb_{\mathbb{R}^k\times Y}}\rightarrow \mathsf{Cob}_d^{emb_{\mathbb{R}^{k+1}\times Y}}.
\]
In this way we obtain maps
\[
\Sigma B\mathsf{Cob}^{emb_{\mathbb{R}^k\times Y}}_d\rightarrow B\mathsf{Cob}^{emb_{\mathbb{R}^{k+1}\times Y}}_d
\]
which become the structure maps for a spectrum whose $k^{th}$ space is $B\mathsf{Cob}^{emb_{\mathbb{R}^k\times Y}}_d$.  Theorem~\ref{thm: embedding cob cat} identifies the weak homotopy type of this spectrum as the infinite loop space of sections of the fiberwise spectrum $MTO^{\epsilon^1\oplus\tau_Y}(d)$ which we define as follows.

Let $\alpha=(P\rightarrow Y)$ be a principal $O(n)$.  There is the bundle  
\[
Th(\gamma^\perp_{d,\epsilon^k\oplus\alpha}):=(P\times_{O(n)} Th(\gamma^\perp_{d,k+n})\rightarrow Y)
\]
with the canonical section $y\mapsto \infty$.   Define the $(k+n)^{th}$ term of the fiberwise spectrum $MTO^{\alpha}(d)$ as $Th(\gamma^\perp_{d,\epsilon^k\oplus \alpha})$.  The structure maps of this fiberwise spectrum, $(S^1\times Y )\wedge_Y  MTO^{\alpha}(d)_l \rightarrow MTO^{\alpha}(d)_{l+1}$, are simply induced by the maps
\[
\Sigma Th(\gamma^\perp_{d,\epsilon^k\oplus \alpha}) \xrightarrow{\simeq} Th(\gamma^\perp_{d,\epsilon^{k+1}\oplus \alpha})
\]
on fibers.

We summarize the previous two paragraphs as follows.  Let $\mathsf{Cob}^{Semb}_d(Y)$ denote the  categorical colimit of the categories $\mathsf{Cob}^{emb_{\mathbb{R}^k\times Y}}_d$ under the functors induced from the inclusions $\{0\}\times\mathbb{R}^k\times Y\hookrightarrow \mathbb{R}^{k+1}\times Y$.

\begin{theorem}
There is a weak homotopy equivalence
\[
B\mathsf{Cob}^{Semb}_d(Y)\simeq \Omega^{\infty-n-1} \Gamma(MTO^{\tau_Y}(d)).
\]

\end{theorem}

Explicitly, when $Y$ is parallelizable this weak homotopy equivalence takes the form
\[
B\mathsf{Cob}^{Semb}_d(Y)\simeq Map_{cpt(Y)}(Y,\Omega^{\infty-n-1} MTO(d))\simeq Map_{cpt(Y)}(Y,\Omega^{-n}B\mathsf{Cob}_d).
\]
In this parallelizable situation we recognize the middle space as one whose homotopy groups are the $MTO(d)$-cohomology groups of the one-point compactification of $Y$.

\begin{remark}
It is important to note the discrepancy in the notation used above.
Recall the equivariant sheaf $Simm_Y$ of \textit{stable} immersions.  Define similarly the functor $Semb_Y$, already alluded to, described by 
\[
Semb_Y= colim_k \text{ }emb_{\mathbb{R}^k\times Y};
\]
the colimit being over the apparent directed system.
This functor $Semb_Y$ is actually an equivariant sheaf and as such it is equivalent to its sheafification which is $Simm_Y$.  The category labeled $\mathsf{Cob}^{Semb}_d(Y)$ is \textit{not} the geometric cobordism category with structure $\mathcal{F}= Semb_Y$.  
\end{remark}

\subsection{Comparing the immersion and embedded cobordism categories}

Choose an embedding $Y\hookrightarrow \mathbb{R}^{N-1}$ and an $\epsilon>0$ such that an $\epsilon$-neighborhood of $Y\subset \{0\}\times \mathbb{R}^{N-1}\hookrightarrow\mathbb{R}^N$ is a tubular neighborhood.  Regarding an embedding as an example of an immersion amounts to a functor $\mathsf{Cob}^{emb}_d(Y)\rightarrow \mathsf{Cob}_{d,N}^{imm_{\mathbb{R}\times Y}}$ realizing, up to weak homotopy equivalence, as a map 
\[
\Gamma_{cpt(Y)}(Th(\gamma^\perp_{d,\epsilon^1\oplus\tau_Y}))\rightarrow \Omega^{N-1} Th(p^*_{imm_{\mathbb{R}\times Y}}\gamma^\perp_{d,N}).
\]
The homotopy fiber of this map then measures the obstruction for a family of $d$-manifolds immersed  into $Y$ to be concordant to a family of $d$-manifolds embedded into $Y$.

Again, for explicitness, examine this map when $\tau_Y$ is stably trivial.  Write $Y^*$ for the one-point compactification of $Y$.  Using superscript notation to denote the Thom space,  the map at hand  becomes 
\[
Map_{cpt(Y)} (Y,(Gr_{d,n+1})^{\gamma^\perp_{d,n+1}})\rightarrow \Omega^{N-1} ((Gr_{d,n+1})^{\gamma^\perp_{d,N}}\wedge Y^*)
\]
This map sends $y\mapsto (V_y,v_y)$ to the assignment $t\mapsto \emptyset$ if $t\in\{0\}\times\mathbb{R}^{N-1}$ is at least $\epsilon $ away from $Y$ and  $t\mapsto (y_t,V_{y_t},v_{y_t} + y_t - t)$ if $t$ is within $\epsilon$ of $Y$ having  $y_t$ as the closest point on $Y$ to $t$.  It is interesting to further specialize to $Y=\mathbb{R}^{N-1}$ where, up to weak homotopy, the map at hand becomes 
\[
Map_{cpt}(\mathbb{R}^{N-1},Th(\gamma^\perp_{d,N}))\xrightarrow{=} \Omega^{N-1} Th(\gamma^\perp_{d,N});
\]
this identity should be regarded as tautological.

\subsection{The immersed cobordism category}

Once again, fix an embedded manifold $Y^n\subset\mathbb{R}^N$.  Consider the cobordism category $\mathsf{Cob}^{imm}_d(Y)$ whose objects are closed $(d-1)$-manifolds immersed in $Y$ and whose morphisms are compact collared $d$-manifolds immersed in $[a_0,a_1]\times Y$.  More precisely, as a topological category 
\[
ob\text{ }\mathsf{Cob}^{imm}_d(Y) = \mathbb{R}^\delta\times (\coprod_{[M^{d-1}]} Imm(M,Y)/{Diff}(M))
\]
and 
\[
mor\text{ }\mathsf{Cob}^{imm}_d(Y) =\{identities\}\amalg (\mathbb{R}^2_+)^\delta \times (\coprod_{[W^d]} Imm(W,[0,1]\times Y)/{Diff}(W)).
\]
This is \textit{not} an example of a cobordism category with $\mathcal{F}$-structure for an equivariant sheaf $\mathcal{F}$.  Nevertheless, we identify its weak homotopy type.  We will only handle the case $Y = \mathbb{R}^n$, the general case is dealt with as described above in regards to $\mathsf{Cob}^{emb}_d(Y)$.  It is also possible to include extra structure $\mathcal{F}$ though this also will not be discussed.

Write $Sp:\mathsf{Top}_*\rightarrow \mathsf{Top}_*$ for the symmetric product functor 
\[
(X,*)\mapsto ((\coprod_k X^k/\Sigma_n)/\sim , *)
\]
where $\Sigma_n$ acts on $X^k$ (with $X^0 := *$) by permuting the factors and $[(*,x_1,...,x_k)] \sim [(x_1,...,x_k)]$ generates the equivalence relation $\sim$.

\begin{theorem}\label{thm:immersed}
There is a weak homotopy equivalence
\[
B\mathsf{Cob}^{imm}_d(\mathbb{R}^n) \simeq \Omega^{n}Sp(Th(\gamma^\perp_{d,n+1})).
\]
\end{theorem}

\begin{remark}
Invoking the theorem of Dold and Thom, the above theorem makes the homotopy groups of $B\mathsf{Cob}^{imm}_d(\mathbb{R}^n)$ explicit as the homology of $Gr_{d,n+1}$.  In particular, for $d=n$ in the oriented case, one observes that the cobordism group of codimension $1$ immersed $d$-manifolds is $\mathbb{Z}$.  
\end{remark}

\begin{proof}[proof (sketch)]
The argument will only be lightly sketched.  The details are similar but more complicated to those in the proof of the main theorem.  We begin with an observation.

Let $W^d$ be an immersed manifold in $\mathbb{R}^{n+1}$.  Then there is an $\epsilon>0$ such that the intersection of $W$ with an $\epsilon$-ball $B_\epsilon(p)\subset\mathbb{R}^{n+1}$ is standard in the following sense.  Choose a diffeomorphism $\phi_{\epsilon,p} : \mathbb{R}^{n+1} \rightarrow B_\epsilon(p)$.   Then there is a diffeomorphism of $\mathbb{R}^{n+1}$ restricting to $\phi_{\epsilon,p}^{-1}(W) \approx V_1 \cup ... \cup V_k$ for some affine $d$-planes $V_i \subset \mathbb{R}^{n+1}$.  Give the collection of such unordered tuples of affine $d$-planes the topology of $Sp(Th(\gamma^\perp_{d,n+1}))$ with the empty $d$-plane assigned to the base-point $\infty\in Th(\gamma^\perp_{d,n+1})$.

Now suppose $W\subset \mathbb{R}^{n+1-k}\times I^k\subset \mathbb{R}^{n+1}$ is bounded in the last $k$ coordinates.  Up to homotopy, the assignment $p\mapsto \phi^{-1}_{\epsilon,p}(W)$ describes a map $\{0\}\times\mathbb{R}^k \rightarrow Sp(Th(\gamma^\perp_{d,n+1}))$ which extends to a based map on the one-point compactification $S^k$ of $\{0\}\times\mathbb{R}^k$.  This assignment $W\mapsto \Omega^k  Sp(Th(\gamma^\perp_{d,n+1}))$ can be done in appropriately topologized compact families of such bounded immersed $W$'s.  This assignment is natural enough to yield a zig-zag
\[
B\mathsf{Cob}^{imm}_d(\mathbb{R}^n) \leftarrow ... \rightarrow \Omega^n Sp(Th(\gamma^\perp_{d,n+1}))
\]
as in the proof of the main theorem.

The foremost subtlety in adapting the proof of the main theorem to show that these maps are weak homotopy equivalences is (1) in transversality issues with immersed manifolds, and (2) in showing $\Psi^{imm}_d(\mathbb{R}^{n+1}):=\{W\subset\mathbb{R}^{n+1}\mid W\text{ is immersed}\}$, appropriately topologized, is weakly equivalent to $Sp(Th(\gamma^\perp_{d,n+1}))$.  The first point can be handled with sufficient care.  The second point amounts to an elaboration of \textsection\ref{subsec: homotopy type of psi}.

\end{proof}

\subsection{Singular maps a la jet spaces}

Take $J^k(\mathbb{R}^d,\mathbb{R}^n)$ to be the space of equivalence classes of smooth maps from $(\mathbb{R}^d,\{0\})$ to $(\mathbb{R}^n,\{0\})$ where two maps are declared equivalent if their Taylor expansions are identical through degree~$k$.  There is an obvious action of ${Diff}(\mathbb{R}^d)^{op}\times {Diff}(\mathbb{R}^n)$ on $J^k(\mathbb{R}^d,\mathbb{R}^n)$.  It is possible then to form the bundle $J^k(W^d,Y^n) \rightarrow W^d\times Y^n$ whose fiber over $(w,y)$ is $J^k(\mathbb{R}^d,\mathbb{R}^n)$.  Local Taylor expansions of a map $f:W \rightarrow Y$ together give a section $j^k(f)$ of this bundle.

Let $\mathcal{S}\subset J^k(\mathbb{R}^d,\mathbb{R}^n)$ be a ${Diff}(\mathbb{R}^d)^{op}\times {Diff}(\mathbb{R}^n)$-invariant subspace of the space of $k$-jets.  There is then a bundle $\mathcal{S}^Y_W \rightarrow W\times Y$ whose fiber over a point $(w,y)$ is $\mathcal{S}$.  A smooth map $f: W^d\rightarrow Y^n$ is said to be \textit{$\mathcal{S}$-singular} if $j^k(f)$ lands in $\mathcal{S}^Y_W$.  A section of $\mathcal{S}^Y_W \rightarrow W\times Y$ is said to be \textit{holonomic} if it is of the form $j^k(f)$ for some $f:W\rightarrow Y$.

Notice that the condition of being $\mathcal{S}$-singular is local.  So the space $\mathcal{S}_Y(W)$ of $\mathcal{S}$-singular maps $W\rightarrow Y$ prescribes a $d$-dimensional equivariant sheaf,
\[
\mathcal{S}_Y:{\mathsf{Emb}_d}^{op} \rightarrow \mathsf{Top}.
\]  
The space $\mathcal{S}_Y(W)$ is the space of holonomic sections of $J^k(W,Y) \rightarrow W$ (see~\cite{eliashberg:h-principle}).

Being embedded in euclidean space, $Y\subset \mathbb{R}^\infty$ inherits a canonical metric,thus ensuring that the structure group of the tangent bundle $\tau_Y$ lies in $O(n)$.  Write $P_{\tau_Y}\rightarrow Y$ for the principal $O(n)$-bundle associated to  $\tau_Y$.  The $Diff$-action on $\mathcal{S}$ restricts to an action of $O(d)^{op}\times O(n)$ on $\mathcal{S}$.

The space $\tau\mathcal{S}_Y(W)$ is the space $F\mathcal{S}_Y(W)$ of \textit{formal $\mathcal{S}$-singular maps}, or equivalently, the space of all sections of $J^k(W,Y) \rightarrow W$.  The main theorem shows that the classical h-principle comparison map $\mathcal{S}_Y(W)\rightarrow F\mathcal{S}_Y(W)$ for each $W$ amounts to a weak homotopy equivalence 
\[
B\mathsf{Cob}^{\mathcal{S}_Y}_d \rightarrow B\mathsf{Cob}^{F\mathcal{S}_Y}_d.
\]  
The following theorem makes this weak homotopy type explicit.

\begin{theorem}
There is a weak homotopy equivalence
\[
B\mathsf{Cob}^{\mathcal{S}_Y}_d \simeq \Omega^{\infty-1} ((P_{\tau_Y}\times_{O(n)} \mathcal{S})//O(d))^{-\gamma_d}.
\]
\end{theorem}

\begin{proof}
Interested in describing the homotopy type of $B\mathsf{Cob}^{\mathcal{S}_Y}_d$, consider the space $Stalk_0({\mathcal{S}_Y}_{\mid_{\mathbb{R}^d}})$.  A point in this stalk is represented by an $\mathcal{S}$-singular map $\mathbb{R}^d\xrightarrow{f} Y$.  After possibly re-choosing the representative, via the preferred choice of exponential map on $Y$ such a map is the data of a point $y\in Y$ and a $\mathcal{S}$-singular map $\mathbb{R}^d\rightarrow T_y Y$.

This describes a weak homotopy equivalence
\[
Stalk_0({\mathcal{S}_Y}_{\mid_{\mathbb{R}^d}})   \xrightarrow{\simeq} P_{\tau_Y}\times_{O(n)} \mathcal{S}.
\]
The map above is equivariant with respect to the action of $O(d)$ on $\mathcal{S}$.  We arrive at a weak homotopy equivalence
\[
B_{\mathcal{S}_Y} \simeq (P_{\tau_Y}\times_{O(n)} \mathcal{S})//O(d)
\]
over $BO(d)$.  The result follows from the main theorem.
\end{proof}

Taking $\mathcal{S}\subset J^1(\mathbb{R}^d,\mathbb{R}^n)$ to be the space of injective maps ($d\leq n$) recovers $imm_Y$, the sheaf of immersions into $Y$, from above.

Taking $\mathcal{S}\subset J^1(\mathbb{R}^d,\mathbb{R}^n)$ to be the space of surjective maps ($d\geq n$) recovers a sheaf of submersions to $Y$ labeled as $sub_Y$.  Interestingly, provided $Y$ is closed, the space of \textit{closed} morphisms of $\mathsf{Cob}^{sub_Y}_d$ is a space of smooth bundles over $Y$ with $d$-dimensional fiber.       
The space $((P_{\tau_Y}\times_{O(n)} \mathcal{S})//O(d))^{-\gamma_d}$ appearing in the theorem above is homotopy equivalent to the space of data $(v,V,y,f)$ where $V^d$ is a $d$-plane in $\mathbb{R}^\infty$, $v\perp V$,  $y\in Y$, and $f: V\rightarrow T_y$ is a surjection of vector spaces.  For explicitness, when $Y$ is parallelizable,
\[
B\mathsf{Cob}^{sub_Y}_d\simeq \Omega^{\infty-1} ((Fr_{n,d}//O(d))^{-\gamma_d}\wedge Y_+).  
\]

\subsection{A category of flow lines}

Fix a Riemannian manifold $Y$ with a Morse function $Y\xrightarrow{f} \mathbb{R}$.  As an application of the above subsection, consider a category whose space of objects are roughly points in $Y$ and whose morphisms are Morse flows.  We make this more precise.

Take $d=1$ and $\mathcal{F}$ to be given by $\mathcal{F}(W) = flow_Y(W) :=\{(v,\gamma, \lambda)\}$ where $v: W\rightarrow Fr_1(TW)$ is a framing of $W$, regarded as a non-vanishing vector field, and $(\gamma,\lambda): W\rightarrow Y\times (-\infty,0]$ is a piecewise smooth map whose non-smooth image points are critical points of $f$ and whose locus of smooth points satisfies the differential equation
\[
\lambda  D\gamma(v) = \nabla f.
\]
The category of interest is $\mathsf{Cob}^{flow_Y}_1$.

\begin{theorem}
There is a weak homotopy equivalence
\[
B\mathsf{Cob}^{flow_Y}_1 \simeq Q Y_+.
\]
\end{theorem}

\begin{remark}
This result is consistent with results of R.L. Cohen,  J.D.S. Jones, and G.B. Segal as follows.  Consider a category $Flow_Y$ whose space of objects is the (discrete) set of critical points of $f:Y\rightarrow \mathbb{R}$.  The space of morphisms $mor_{Flow_Y}(p_0,p_1)$ is the (compactified) moduli space of flow lines between $p_0$ and $p_1$.  Roughly, a morphism is an unparametrized flow line between critical points, or a collection of concatenating unparametrized flow lines between intermediate critical points.  The canonical map $BFlow_Y\rightarrow Y$ is a homotopy equivalence.  See~\cite{cohen-jones-segal:morse} for details and a proof of this fact.  A choice of an embedding $Y\hookrightarrow \mathbb{R}^\infty$ yields a functor $Flow_Y \rightarrow \mathsf{Cob}^{flow_Y}_1$.  Upon taking classifying spaces, this realizes the standard inclusion $Y\rightarrow QY_+$.  This is similar to the situation $\mathsf{P}Y \rightarrow \mathsf{Cob}^{\mathcal{O}r\times map_Y}_1$ realizing to $Y \rightarrow QY_+$ where $\mathsf{P}Y$ is the path category on $Y$.
\end{remark}

\begin{remark}
Instead of a Riemannian manifold $Y$ with a Morse function one could consider a manifold $Y$ with a generic vector field $V$. 
\end{remark}

\begin{proof}
There is a forgetful morphism of equivariant sheaves $flow_Y \rightarrow fr\times map_{Y}$ to the sheaf of framings and smooth maps to $Y$.  This induces a map on stalks $Stalk_0({flow_Y}_{\mid_\mathbb{R}}) \rightarrow Stalk_0({fr\times map_{Y}}_{\mid_\mathbb{R}}) = \mathbb{Z}/2\times Y$.  There is a homotopy inverse described as follows.

Let $C\subset Y$ be the discrete set of critical points of $Y$.  Send $y\in Y\setminus C$ to the germ $(-\epsilon,\epsilon)\xrightarrow{\gamma} Y$ with $\gamma(0)=y$, satisfying the differential equation at hand for $\lambda(s) = -\sqrt{\lVert \nabla f(\gamma(s))\rVert}$.  Note that such a germ is uniquely defined by this data.  Extend this map to all of $Y$ by sending $y\in C$ to the constant map $\gamma(s) \equiv y$.  This extension is continuous and establishes a homotopy inverse.
The result then follows from the main theorem.
\end{proof}

\section{Holomorphic maps}\label{sec:holomorphic-maps}

\subsection{The sheaf of complex structures}\label{subsec: sheaf of complex stctrs}

Take $\mathcal{F}: {\mathsf{Emb}_d}^{op}  \rightarrow   \mathsf{Top}$ to be the sheaf of complex structures.  That is, declare
\[
\mathcal{F}(W) = \mathbb{C}(W):=\{\mathcal{A} = \text{ maximal holomorphic atlas on }W\}.
\]
To describe the topology of this set observe the bijection between  $\mathbb{C}(W)$ and integrable almost-complex structures on $W$.  Topologize the set $\mathbb{C}(W)$ as a subspace of  the space of almost-complex structures denoted $\mathcal{J}(W)$.  This inclusion $\mathbb{C}(W^d)\hookrightarrow\mathcal{J}(W^d)$ is sufficiently natural to describe a morphism of equivariant sheaves.

\begin{theorem}
The inclusion morphism of equivariant sheaves $\mathbb{C}\hookrightarrow \mathcal{J}$ induces a weak homotopy equivalence
\[
B\mathsf{Cob}_d^\mathbb{C}\xrightarrow{\simeq} B\mathsf{Cob}_d^\mathcal{J}\simeq \Omega^{\infty-1} MTU(d/2).
\]
\end{theorem}

\begin{proof}
We are interested in comparing the stalks
\[
Stalk_0(\iota^*_{\mathbb{R}^d}\mathbb{C})\rightarrow Stalk_0(\iota^*_{\mathbb{R}^d}\mathcal{J}).
\] 
The equivariant sheaf $\mathcal{J}$ is a \textit{fiberwise} equivariant sheaf coming from the fibration $BU(d/2)\rightarrow BO(d)$ by its very definition (see example (1)).  Moreover, $\mathcal{J}$ is exactly the tangential equivariant sheaf $\tau\mathbb{C}$ associated to the equivariant sheaf $\mathbb{C}$.  Indeed, on the one hand $Stalk_0(\iota^*_V\mathcal{J})$ is the space of linear complex structures on the vector space $T_0 V \cong V$; on the other hand,  an element $i\in Stalk_0(\iota^*_V\tau\mathbb{C})$ is an equivalence class of integrable almost-complex structures on neighborhoods of $0\in V$ which, because any complex structure is locally standard, can be identified with linear complex structures on $T_0 V$.  The result follows from the main theorem.
\end{proof}

\subsection{The sheaf of holomorphic maps}\label{subsec: sheaf of holo maps}

Fix a complex manifold $Y=(Y,i)$.  Consider the equivariant sheaf $\mathcal{J}\times map_Y:{\mathsf{Emb}_d}^{op}\rightarrow \mathsf{Top}$ determined by
\[ 
\mathcal{J}\times map_Y(M):=\{(J,f)\mid J\in\mathcal{J}(M)\text{ and } f\in Map(M,Y)\}
\]
with the obvious topology.  This sheaf $\mathcal{J}\times map_Y$ is a \textit{fiberwise} equivariant sheaf coming from the fibration $Y\times BU(d/2)\rightarrow BO(d)$.

Consider a similar equivariant sheaf 
$\mathbb{C}_Y:{\mathsf{Emb}_d}^{op}\rightarrow \mathsf{Top}$
described by 
\[
\mathbb{C}_Y(W):=\{(\mathcal{A},h)\mid \mathcal{A}\text{ is a maximal holomorphic atlas on $W$ and } h\in Hol((W,\mathcal{A}),Y)\}.  
\]
This set is topologized as a subspace of $\mathcal{J}\times map_Y(W)$.
Because the above data is local, it is straight forward to verify that the sheaf-gluing property holds for $\mathbb{C}_Y$.

\begin{theorem}
The inclusion morphism of equivariant sheaves $\mathbb{C}_Y\hookrightarrow \mathcal{J}\times map_Y$ induces a weak homotopy equivalence
\[
B\mathsf{Cob}_d^{\mathbb{C}_Y}\simeq B\mathsf{Cob}_d^{\mathcal{J}\times map_Y}\simeq \Omega^{\infty-1} MTU(d/2)\wedge Y_+.  
\]
\end{theorem}

\begin{proof}
From the main theorem, it is sufficient to show
\begin{equation}\label{eq: complex stalk}
Stalk_0({\mathbb{C}_Y}_{\mathbb{R}^d}) \xrightarrow{\simeq }Stalk_0({\mathcal{J}\times map_Y}_{\mathbb{R}^d}).
\end{equation}
A point in the space on the right-hand side is the data of a linear almost-complex structure $J$ on $\mathbb{R}^d$ and $y\in Y$.  Regard $\mathbb{R}^d$ as a vector bundle over a point.  
A point in the space on the left hand side is the data of morphism of almost-complex vector bundles $(\mathbb{R}^d,J)\rightarrow (\tau_Y,i)$ which is the data of a point $y\in Y$ and a map of almost-complex vector spaces $\phi:(\mathbb{R}^d,J) \rightarrow (T_y Y,i)$.  The assignment $(J,y)\mapsto ((\mathbb{R}^d,J)\xrightarrow{0} (T_y Y, i))$, the zero-map of vector bundles, is a homotopy inverse to the map in \ref{eq: complex stalk} as seen by the homotopy $\phi_t = (1-t)\phi$, $t\in[0,1]$.
\end{proof}

Let $(Y,\omega)$ be a symplectic manifold.  Let $j\in\mathcal{J}(Y)$ be an almost-complex structure which is \textit{compatible} with $\omega$; that is, the $2$-tensor $\omega(-,j -)$ is a metric on $Y$.   It is possible to perform the analysis above on  the symplectic almost-complex manifold $(Y,\omega,j)$  and consider the category of $j$-holomorphic maps to $Y$.  Precisely, consider the equivariant sheaf $\mathcal{J}_Y$ determined by 
\[
\mathcal{J}_Y(M)=\{(J,h)\mid J\in\mathcal{J}(M)\text{ and } h\in Hol((M,J),(Y,j))\},
\]
topologized appropriately.  The same proof above shows the weak homotopy equivalence 
\[
B\mathsf{Cob}_d^{\mathcal{J}_Y}\simeq B\mathsf{Cob}_d^{\mathcal{J}\times map_Y}\simeq \Omega^{\infty-1} MTU(d/2)\wedge Y_+.
\]

\subsection{The category of holomorphic curves}

Fix $Y=(Y,i)$ a complex manifold.  The discussion to follow holds equally well for $Y=(Y,\omega,j)$ a symplectic almost-complex manifold.
Let $\mathcal{O}r$ be the \textit{fiberwise} equivariant sheaf of orientations coming from the fibration $BSO(d)\rightarrow BO(d)$.  There is a morphism of embedding sheaves $\mathcal{J}\rightarrow \mathcal{O}r$ due to the observation that a linear complex structure on a vector space canonically determines an orientation on that vector space.  When $d=2$, this morphism of embedding sheaves is a homotopy equivalence.  The result of the above~\S\ref{subsec: sheaf of holo maps} implies the weak homotopy equivalence 
\[
B\mathsf{Cob}^{\mathbb{C}_Y}_2\xrightarrow{\simeq} B\mathsf{Cob}^{\mathcal{O}r\times map_Y} _2\simeq \Omega^{\infty-1} MTSO(2)\wedge Y_+.
\]

Also, when $d=2$ the morphism of embedding sheaves $\mathbb{C}\rightarrow \mathcal{J}$ is a  homeomorphism because in $2$-dimensions all almost-complex structures are integrable.  Fix a morphism $F\in mor\text{ }\mathsf{Cob}_2$  of genus $g$ with $p$ incoming and $q$ outgoing boundary components.  The orientation preserving (collared) diffeomorphism group $Diff^+(F)$ acts on $\mathbb{C}_Y(F)$ given explicitly by pulling back complex structures and pre-composing with maps to $Y$.  The space of morphisms in $\mathsf{Cob}_2^{\mathbb{C}_Y}$ with underlying source and target oriented $1$-manifolds $s(F)$ and $t(F)$  has 
\[
Emb(F,[0,1]\times \mathbb{R}^\infty)\times_{Diff^+(F)} \mathbb{C}_Y(F)\simeq \mathbb{C}_Y(F)//Diff^+(F)
\]
as a component.  The projection 
\[
\mathbb{C}_Y(F)//Diff^+(F)\rightarrow \mathbb{C}_Y(F)/Diff^+(F)
\]
is a homotopy equivalence provided $p+q>0$ and is a rational homotopy equivalence otherwise.

The quotient space $\mathbb{C}_Y(F)/Diff^+(F) $  is the moduli space $\mathcal{M}_{g,p,q}(Y)$ of holomorphic curves in the sense of Gromov-Witten theory.  We thus have natural maps, well-defined up to homotopy (only rationally when $p=q=0$), 
\[
\mathcal{M}_{g,p,q}(Y)\rightarrow \Omega^\infty MTSO(2)\wedge Y_+
\]
for all $g,p,q\geq 0$ such that gluing holomorphic curves 
\[
\mathcal{M}_{g_0,p,q}(Y)\times_{(LY)^q}\mathcal{M}_{g_1,q,r}(Y)\rightarrow \mathcal{M}_{g_0+g_1+q-1,p,r}(Y)
\]
agrees, up to homotopy, with concatenation of loops in $\Omega^\infty MTSO(2)\wedge Y_+$.  Clearly, the right hand side above is independent of $g$, $p$, and $q$ and thus only reflects the `stable' topology of the moduli spaces in the sense of \textsection\ref{subsubsec: group completion}.  In fact, rationally, the loop space $\Omega^\infty MTSO(2) \wedge Y_+$ is initial among all such collections of maps $\mathcal{M}_{g,p,q}\rightarrow \Omega X$ to a loop space (see~\textsection\ref{subsubsec: group completion}).  In this way, the loop-space $\Omega^\infty MTSO(2)\wedge Y_+$ parametrizes `stable' rational characteristic classes of the moduli spaces of holomorphic curves in $Y$.  This is summarized as the following

\begin{statement*}
The `stable' moduli space of holomorphic curves in a complex manifold $Y$ has the weak homotopy type of $\Omega^\infty MTU(d/2)\wedge Y_+$
\end{statement*}

It is of interest to know connectivity of the maps 
\[
\mathcal{M}_g(Y)\rightarrow \Omega^\infty MTSO(2)\wedge Y_+.
\]
Let $\mathcal{M}^\alpha_g(Y)\subset\mathcal{M}_g(Y)$ be the subspace consisting of holomorphic maps $h$ having homological class $h_*[F]=\alpha\in H_2(Y)$; this subspace is a union of connected components.  For $Y={\mathbb{C}P}^n$, it is shown in~\cite{ayala:param-segal} that  for $\alpha = d>>g>>1$, 
\[
H_k(\mathcal{M}^d_g(\cp);\mathbb{Q})\cong H_k(\Omega^\infty MTSO(2)\wedge {\mathbb{C}P}^n_+;\mathbb{Q}),
\]
the left hand side being well understood.  In general, one expects  $\Omega^\infty MTSO(2)\wedge Y_+$ to  approximate the moduli spaces $\mathcal{M}^\alpha_g(Y)$ for $\alpha$ and $g$ large in some sense.

\section{Symplectic cobordism categories}

\subsection{Symplectic structures}

Begin by considering the equivariant sheaf of symplectic structures, written $Symp:{\mathsf{Emb}_d}^{op}\rightarrow \mathsf{Top}$.  
Recall that a symplectic structure on an even dimensional manifold $W^d$ is a maximal rank closed $2$-form $\omega$.  Maximal rank here means that locally $\omega$ can be written as a $d\times d$ matrix having maximal rank $=d$.  So $Symp$ is given by 
\[
Symp(W)=\{\omega \text{ a symplectic form on }W\};
\]
topologized as a subspace of the space of sections of $\Lambda^2\tau^*_W$.
Note that a closed $2$-form $\omega$ on $W^d$ being maximal rank is equivalent to $\omega^{d/2}$ being a volume form on $W^d$.  In this way there is a morphism of embedding sheaves $Symp\rightarrow \mathcal{O}r$.

\begin{theorem}
There is a weak homotopy equivalence
\[
B\mathsf{Cob}_d^{Symp} \simeq \Omega^{\infty-1} MTU(d/2).
\]
\end{theorem}

\begin{proof}
There is a canonical symplectic structure $\omega_0=\Sigma_i d x_i \wedge d y_i$ on the vector space $\mathbb{R}^d\cong {\mathbb{R}^2}^{(d/2)}=\{(x_i,y_i)\}^{d/2}_{i=1}$.  Define the \textit{symplectic group} $Sp(d):=\{A\in GL_d(\mathbb{R})\mid A^* \omega_0 = \omega_0\}$.   It is a standard exercise to show that $BSp(d) \simeq BU(d/2)$.  
Write $Sp$ for the equivariant sheaf associated to the fibration $BSp(d)\rightarrow BO(d)$.  Explicitly, $Sp(W)$ consists of non-degenerate $2$-forms on $W$.  There is the obvious inclusion morphism of equivariant sheaves
\[
Symp\hookrightarrow Sp
\]
given by forgetting the closedness condition of a symplectic structure.  Moreover, this inclusion induces a homotopy equivalence 
\[
Stalk_0(Symp_{\mid_{\mathbb{R}^d}})\xrightarrow{\simeq}  Stalk_0 ( Sp_{\mid_{\mathbb{R}^d}})
\]
since, among other reasons,  any symplectic structure is locally standard (Darboux's theorem).
The theorem follows upon invoking the main theorem.
\end{proof}

\subsection{Lagrangian immersions}

Let $V^{2d}=(V^{2d},\omega)$ be a $2d$-dimensional symplectic vector space.  A $d$-dimensional subspace $U\subset V$ is said to be \textit{Lagrangian} if  $\omega_{\mid_U}=0$.  Let $Y^{2d}=(Y^{2d},\omega)\subset\mathbb{R}^\infty$ be a $2d$-dimensional embedded symplectic manifold and let $W^d$ be any smooth manifold.  A smooth map $W\xrightarrow{f} Y$ is said to be a \textit{Lagrangian immersion} if for each $w\in W$, $T_w W\xrightarrow{D_w f} T_{f(w)} Y$ is an embedding of vector spaces having image a Lagrangian subspace.  Related to the equivariant sheaf $imm_Y$ is the equivariant sheaf $lag_Y:{\mathsf{Emb}_d}^{op}\rightarrow \mathsf{Top}$ determined by 
\[
lag_Y(M)=\{M\xrightarrow{f}Y \mid f\text{ is a Lagrangian immersion}\}  
\]
topologized in the obvious way.

Denote by $Gr^{lag}_{d,2d}$ the Grassmann manifold of Lagrangian subspaces of $(\mathbb{R}^{2d},\omega_0)$.  There is an action of $Sp(2d)$ on $Gr^{lag}_{d,2d}$ given by acting on the ambient $\mathbb{R}^{2d}$.  Let $P_{\tau_Y}\rightarrow Y$ be the principal $Sp(2d)$-bundle associated to $(\tau_Y,\omega)$.  Write $Gr^{lag}_d(TY, \omega)$ for the space $P_{\tau_Y} \times_{Sp(2d)} Gr^{lag}_{d,2d}$.  Because $Y\subset\mathbb{R}^\infty$ is embedded, there is a map 
\[
Gr^{lag}_d(TY,\omega)\rightarrow BO(d)
\]
given by $(V^d\subset T_y Y)\mapsto( V^d\subset T_y \mathbb{R}^\infty\cong\mathbb{R}^\infty)$.

\begin{theorem}
There is a weak homotopy equivalence
\[
B\mathsf{Cob}^{lag_Y}_d\simeq \Omega^{\infty-1}(Gr^{lag}_d(TY,\omega))^{-\gamma_d}.
\]
\end{theorem}

\begin{proof}
We are interested in the weak homotopy type of the stalk $Stalk_0({lag_Y}_{\mid_{\mathbb{R}^d}})$.  
In general, for $\alpha$ a $d$-dimensional vector bundle and  $(\beta,\omega)$ a $2d$-dimensional symplectic vector bundle, let $Lag(\alpha,\beta)$ be the space of bundle morphisms $\alpha\rightarrow \beta$ which are fiberwise Lagrangian embeddings of vector spaces.  Regarding $\mathbb{R}^d$ as a vector bundle over a point $*$, it is not difficult to verify
\[
Stalk_0({lag_Y}_{\mid_{\mathbb{R}^d}})\xrightarrow{\simeq} Lag(\mathbb{R}^d,\tau_Y).  
\]
Explicitly, 
\[
Lag(\mathbb{R}^d, \tau_Y) = P_{\tau_Y}\times_{Sp(2d)} Fr^{lag}_{d,2d}
\]
where $Fr^{lag}_{d,2d}$ is the space of $d$-frames in $(\mathbb{R}^{2d},\omega_0)
$ which span Lagrangian subspaces.  The space $Fr^{lag}_{d,2d}$ is the tautological principal $O(d)$-bundle over the Grassmann space $Gr^{lag}_{d,2d}$ of Lagrangian subspaces of $\mathbb{R}^{2d}$.  The group $Sp(2d)$ acts transitively on $Gr^{lag}_{d,2d}$ with stabilizer of a point $O(d)\subset Sp(2d)$, the diagonal embedding.  We conclude that $Lag(\mathbb{R}^d, \tau_Y) = P_{\tau_Y}$.  
It follows that $Gr^{lag}_d(TY, \omega)$ is the space $P_{\tau_Y}/O(d)$.

In the end there is a weak homotopy equivalence
\[
B{lag_Y}=Fr_d\times_{O(d)} P_{\tau_Y}\xrightarrow{\simeq} Gr^{lag}_d(TY,\omega);
\]
over $BO(d)$ where
the middle term maps to $BO(d)$ by projection onto the first coordinate.  From the main theorem we arrive at the weak homotopy equivalence
\[
B\mathsf{Cob}^{lag_Y}_d\simeq \Omega^{\infty-1}(Gr^{lag}_d(T,\omega))^{-\gamma_d}.
\]
\end{proof}

\subsection{Contact forms}

Recall that a $1$-form $\alpha$ on an odd-dimensional manifold $M^d$ is said to be a \textit{contact form} if $d\alpha$ has maximal rank; that is, locally $d\alpha$ can be written as an $(odd)\times(odd)$ matrix having maximal rank.  A (co-oriented) \textit{contact structure} is a hyperplane distribution $\{\alpha = 0\}$ for some contact form $\alpha$.  A \textit{contact manifold} is a pair $(M,\alpha)$ where $\alpha$ is a contact form on $M$.  Write $Cont(M)$ for ths space of contact forms on $M$ topologized as a subspace of the space of smooth sections of $\Lambda\tau^*_M$.  Because a contact form is locally defined, it is routine to verify that 
\[
Cont:{\mathsf{Emb}_d}^{op}\rightarrow \mathsf{Top}
\]
is indeed an equivariant sheaf.

For $d$ odd, consider the \textit{fiberwise} equivariant sheaf $Hyp^{Symp}$ associated to the fibration $B_{hyp}\rightarrow BO(d)$ where
\[
B_{hyp} = \{(V,\alpha,\omega) \mid \alpha \in V^*\setminus\{ 0 \}  \text{ and } \omega \in Sp(ker(\alpha))\}
\]
with $\alpha$ co-oriented.  The association $\alpha  \rightsquigarrow \ker(\alpha)$ justifies the notation.  This set $B_{hyp}$ is topologized in an obvious way.  The projection $B_{hyp}\rightarrow BO(d)$ is a fiber bundle with fiber $Fr^{or}_{d-1,d}\times_{SO(d-1)} Sp^{lin}(\mathbb{R}^{d-1})$ where $Sp^{lin}(V)$ is the space of \textit{linear} symplectic structures on $V$.   This fiber is the total space of a bundle over $S^{d-1} \cong Gr^{or}_{d-1,d}$ with fiber $Sp^{lin}(\mathbb{R}^{d-1})\simeq \mathcal{J}_{or}(\mathbb{R}^{d-1})$.    
Again, because every contact form is locally standard, there is a weak homotopy equivalence
\[
Stalk_0(Cont_{\mid_{\mathbb{R}^d}})\xrightarrow{\simeq} Stalk_0(Hyp^{symp}_{\mid_{\mathbb{R}^d}}).
\]

To be explicit, specialize to the case $d=3$ when $\mathcal{J}_{or}(\mathbb{R}^2)\simeq *$ is contractible.  In this case, $Stalk_0(Hyp^{symp}_{\mid_{\mathbb{R}^d}})//O(3) \simeq S^2//O(3)$, the unit sphere bundle of $\gamma_3$.  So $S^2//O(3) \simeq BO(2)$.  The pull back of $\gamma_3$ to $BO(2)$ being $\epsilon\oplus\gamma_2$, from the main theorem there is the weak homotopy equivalence
\[
B\mathsf{Cob}^{Cont}_3\simeq \Omega^{\infty-2} MTO(2).
\]

\subsection{Restricting the objects of cobordism categories}\label{subsec: Restricting objects}

\subsubsection*{The construction}

Let $\mathsf{C}$ be a topological category with a continuous map $O\xrightarrow{f} ob\text{ }\mathsf{C}$.  There is the \textit{pull back} topological category $f^*\mathsf{C}$ whose space of objects is $O$ and whose space of morphisms is  the pull back 
\[
\xymatrix{
(f\times f)^* mor\text{ }\mathsf{C}  \ar[d]^{s,t}   \ar[r] 
&
mor\text{ }\mathsf{C}  \ar[d]^{s,t}
\\
O\times O  \ar[r]^{f\times f}
&
ob\text{ }\mathsf{C}\times ob\text{ }\mathsf{C}.
}
\]
The source, target, and identity maps are apparent and are continuous.
Composition is determined by composition in $\mathsf{C}$.

There is the following thesis which indicates how little the topology of $B\mathsf{C}$ reflects the topology of $ob\text{ }\mathsf{C}$.  Regard the topological category $\mathsf{C}$ as an $\infty$-category as in \textsection\ref{subsubsec: group completion} and recall the definition of a morphism being invertible in that context.  Two objects $o_0,o_1\in ob\text{ }\mathsf{C}$ are said to be \textit{equivalent} if there is an invertible morphism $o_0\rightarrow o_1$ in $\mathsf{C}$.  It is not difficult to verify that this is indeed an equivalence relation $\sim$ on $ob\text{ }\mathsf{C}$.

\begin{thesis}\label{thm: skeleton}

If $O\xrightarrow{f}ob\text{ }\mathsf{C}\rightarrow ob\text{ }\mathsf{C}/\sim$ is surjective then the universal morphism 
\[
Bf^*\mathsf{C}\xrightarrow{\simeq} B\mathsf{C}.
\]
is a homotopy equivalence.
\end{thesis}

The situation of the hypothesis of the theorem is not uncommon.  Indeed, if the source-target map $mor\rightarrow ob\times ob$ is a fibration, then establishing surjectivity of $\pi_0 O\rightarrow \pi_0 ob$ is sufficient.

\subsubsection*{Example: the path category}

Let $A\subset X$ be a pair of spaces.  Recall the \textit{path category} $\mathcal{P}_A X$from \textsection\ref{subsubsec: group completion} defined as follows.  Declare 
\[
ob\text{ }\mathcal{P}_A X := A
\hspace{.3cm}\text{ and }\hspace{.3cm}
mor\text{ }\mathcal{P}_A X :=  Map(([0,1],\{0,1\}),(X,A)).
\]
The source and target maps are given by the two evaluation maps $ev_\nu:Map([0,1],X)\rightarrow X$, $\nu= 0,1$ respectively.  Composition is given by concatenation of reparametrized paths.  Apparently, composition is not associative, but only up to a contractible choice of homotopies.  But this should not bother us since the path category $\mathcal{P}_A X$ is an $\infty$-category in the sense of \textsection\ref{subsubsec: group completion}.  When $A=X$ simply write $\mathcal{P}X$ and when $A=*$ write $\Omega X$. 
From Theorem~\ref{thm: skeleton} it follows that, for $X$ connected and $A\neq \emptyset$,
\[
B\mathcal{P}_A X \simeq B\mathcal{P} X \simeq B\Omega X \simeq X.
\]
This example captures much of the general behavior.

\subsection{Contact-symplectic cobordism category}

Recall that for $M\in\mathsf{Emb}_{d-1}$ with $d$ odd, $Symp(M):=colim_\epsilon Symp((-\epsilon,\epsilon)\times M)$.  In this way 
\[
Symp:\mathsf{Emb}^{op}_{d-1}\rightarrow \mathsf{QTop}
\]
becomes an equivariant sheaf for $(d-1)$ odd.  There is a morphism of $(d-1)$-dimensional equivariant sheaves 
\[
Cont\xrightarrow{S} Symp,
\]
called \textit{symplectization}, defined by the assignment 
\[
\alpha\mapsto d(e^t\alpha),\hspace{.5cm}t\in(-\epsilon,\epsilon).  
\]
It is worth remarking that this assignment is injective.

It is then possible to define the \textit{contact-symplectic} cobordism category 
\[
\mathsf{Cob}^{cSymp}_d:= S^*\mathsf{Cob}^{Symp}_d
\]
as in \textsection\ref{subsec: Restricting objects} above.  So the objects of $\mathsf{Cob}^{cSymp}_d$ are closed contact $(d-1)$-manifolds and the morphisms are (strong) symplectic cobordisms.

\begin{conj}\label{thm: contact-symplectic}
The universal map
\[
B\mathsf{Cob}^{cSymp}_d\rightarrow B\mathsf{Cob}^{Symp}_d \simeq \Omega^{\infty-1}MTU(d/2).
\]
is a weak homotopy equivalence.
\end{conj}

A sketch-proof of the conjecture is presented as Proposition~\ref{prop: contact} in the next section.

\subsection{The symplectic field theory cobordism category}

In this subsection, we identify the homotopy type of the source cobordism category of Symplectic Field Theory (SFT), which we denote by $\mathsf{Cob}_d^\mathcal{SFT}$.  We make first set things up.

\subsubsection*{Definitions}

The cobordism category $\mathsf{Cob}_d^\mathcal{SFT}$ can be roughly described as follows.  The objects are $4$-tuples $(M,\Omega,\lambda,J)$ where $M$ is an odd-dimensional oriented closed manifold, $\Omega$ is a \textit{Hamiltonian structure}, $\lambda$ is a \textit{stable framing} of $\Omega$, and $J$ is an almost-complex structure on $(M,\Omega,\lambda)$.  The morphisms are $3$-tuples $(W,\omega,J)$ where $W$ is an even-dimensional manifold regarded as a cobordism, $\omega$ is a symplectic structure on $W$, and $J$ is an almost-complex structure on $(W,\omega)$.   At the moment, it is not clear what the source and target maps are and how to compose morphisms.  We begin by  defining these terms.

A \textit{Hamiltonian} structure  on an odd dimensional oriented manifold $(M^{d-1},\sigma)$, $\sigma\in \mathcal{O}r(M)$, is a closed $2$-form $\Omega$ of maximal rank.  As with the definition of a symplectic structure, maximal rank here means that locally $\Omega$ can be written as an anti-symmetric $(d-1) \times (d-1)$ matrix having maximal rank $=d-2$.  Let $Ham(M)$ denote the space of Hamiltonian structures on $M$ topologized as a subspace of the space of sections $\Gamma(\Lambda^2\tau_M^*)$.

To each oriented Hamiltonian manifold $(M,\Omega)$ is associated its \textit{characteristic} line field $l=l(\Omega)=ker\Omega\hookrightarrow \tau_M$.  So $\Omega$ defines a linear symplectic structure (and thus orientation) on the vector bundle $\tau_M/l$.  A vector field $R=R(\Omega)$ is called \textit{characteristic} if it generates $l$ and the orientation class $[R,\Omega^{(d-1)/2}]$ agrees with that prescribed on $(M,\sigma)$.  A $1$-form $\lambda$ on $M$ is said to \textit{stabilize} the Hamiltonian manifold $(M,\Omega)$ if there is a characteristic vector field $R$ such that the Lie derivative $L_R \lambda =\lambda$ (such an $R$ is unique up to scaling by a constant.  Denote by $R_\lambda$ the unique $R$ with $\lambda(R)=1$.  Such a $1$-form yields a hyperplane distribution $\xi(\lambda)=\{\lambda = 0\}$; necessarily, $\Omega_{\mid_\xi}$ will be a (linear) symplectic structure on $\xi$.    As before, a \textit{compatible} almost-complex structure $J$ on a symplectic vector bundle $(\alpha,\omega)$ will be one such that the $2$-tensor $\omega(-,J-)$ is a metric on $\alpha$; the space of such almost-complex structures being denoted $\mathcal{J}(\alpha,\omega)$ or simply $\mathcal{J}(\omega)$ if the bundle $\alpha$ is understood.   A \textit{framing} of the Hamiltonian manifold $(M,\sigma,\Omega)$ is a pair $(\lambda,J)$ such that $\lambda$ stabilizes $(M,\sigma,\Omega)$ and $J\in \mathcal{J}(\xi,\Omega_{\mid_\xi})$.

Let $Ham^{fr}:\mathsf{Emb}^{op}_{d-1}\rightarrow \mathsf{Top}$, $d$ even, be the equivariant sheaf of framed Hamiltonian structures given as 
\[
Ham^{fr}(M)=\{(\sigma,\Omega,\lambda,J)\mid ((M,\sigma),\Omega,\lambda,J)\text{ is a framed  Hamiltonian manifold}\}.
\]
This set is topologized in an obvious way.  The sheaf and equivariance conditions are routine to verify.

\begin{example}
A contact structure $\alpha\in Cont(M)$ canonically induces a Hamiltonian structure $d\alpha\in Ham(M)$.  This provides a morphism of $(d-1)$-dimensional equivariant sheaves
\[
Cont\rightarrow Ham.
\]
Better still, a pair $(\alpha,J)\in Cont^\mathcal{J}(M)$, consisting of a contact structure and an almost complex structure on $\xi:=\{\alpha=0\}$ compatible with the \textit{linear} symplectic structure $d\alpha$ on $\xi$, induces a \textit{framed} Hamiltonian structure $(d\alpha,\alpha,J)\in Ham^{fr}(M)$.  This describes a morphism of equivariant sheaves
\[
Cont^\mathcal{J}\rightarrow Ham^{fr}.
\]
\end{example}

Let $Symp^\mathcal{J}:{\mathsf{Emb}_d}^{op}\rightarrow \mathsf{Top}$ be the $d$-dimensional equivariant sheaf given as 
\[
Symp^\mathcal{J}(W)=\{(\omega,J)\mid \omega\in Symp(W) \text{ and } J\in\mathcal{J}(W,\omega)\}.
\]
Again, because this data is locally defined and clearly ${Diff}(M)$ invariant, $Symp^\mathcal{J}$ is indeed an equivariant sheaf.  There is a morphism of $(d-1)$-dimensional equivariant sheaves
\[
S:Ham^{fr}\rightarrow Symp^\mathcal{J}
\]
given by 
\[
((M,\sigma),\Omega,\lambda,J)\mapsto (M, \Omega + d(e^t \lambda), J\oplus ((R_\lambda,\partial_t)\mapsto (\partial_t, -R_\lambda)))
\]
where $t\in(-\epsilon,\epsilon)$.  
This morphism is known as \textit{symplectization} and is a general procedure for obtaining a symplectic structure on $(-\epsilon,\epsilon)\times M$ from a framed Hamiltonian structure on $M$.

Using the construction in \textsection\ref{subsec: Restricting objects}, for $d$ even, define
\[
\mathsf{Cob}^\mathcal{SFT}_d := S^*\mathsf{Cob}^{Symp^\mathcal{J}}_d.
\] 
So the objects of the SFT cobordism category are framed Hamiltonian manifolds while the morphisms are symplectic almost-complex manifolds.

\subsubsection*{The weak homotopy type of $B\mathsf{Cob}^\mathcal{SFT}_d$}

\begin{lem}\label{lem: symplectization}
As a morphism of $(d-1)$-dimensional equivariant sheaves, symplectization
\[
Ham^{fr}\xrightarrow{S} Symp
\]
is level-wise a homotopy equivalence
\end{lem}

\begin{proof}
We describe a homotopy inverse.  Given $(\omega,J)\in Symp^\mathcal{J}(M)$, the $2$-form $\Omega:=\omega_{\mid_M}$ is a Hamiltonian structure on $M$.  This describes a morphism of $(d-1)$-dimensional equivariant sheaves
\[
Symp^\mathcal{J}\xrightarrow{r} Ham^\mathcal{J}.
\]
The forgetful morphism of $(d-1)$-dimensional equivariant sheaves
\[
Ham^{fr}\rightarrow Ham^\mathcal{J},
\]
given by forgetting the data of the characteristic $1$-form,
is level-wise a homotopy equivalence.  Indeed, the morphism is level-wise a fibration with fiber the affine space of characteristic $1$-forms.  This establishes the commutative diagram  
\[
\xymatrix{
&
Ham^{fr}   \ar[d]  \ar[dl]^s
\\
Symp^\mathcal{J}  \ar[r]^r
&
Ham^\mathcal{J}.
}
\]
The result follows.
\end{proof}

\begin{conj}\label{thm: SFT}
The universal map
\[
B\mathsf{Cob}_d^\mathcal{SFT} \xrightarrow{\simeq} B\mathsf{Cob}_d^{Symp^\mathcal{J}}
\]
is a weak homotopy equivalence.
\end{conj}

\begin{proof}[Sketch]
We appeal to~\S\ref{subsec: Restricting objects}
from which it suffices to show  $Ham^{fr}(M)\xrightarrow{S} Symp^\mathcal{J}(M) \rightarrow Symp^\mathcal{J}(M)/\sim$ is surjective for each $M^{d-1}$.  For this it is sufficient to show that for each germ $(\omega_0,J)\in Symp^\mathcal{J}(M)$ there exists $(\omega,J)\in Symp^\mathcal{J}(I\times M)$ restricting to $(\omega_0,J)$ on $\{0\}\times M$ and $S(\sigma,\Omega,\lambda, J)$ on $\{1\}\times M$ for some $(\sigma,\Omega,\lambda,J)\in Ham^{fr}(M)$.

Conveniently, there is a classical h-principle due to Gromov (see~\cite{eliashberg:h-principle} for a general reference on the subject) which implies that  the inclusion of the space of symplectic structures into the space of \textit{almost} symplectic structures, $Symp(M) \hookrightarrow Sp(M)$, is a weak homotopy equivalence.  Because $\mathcal{J}$ is a fiberwise structure on $M$, the same is true for the inclusion $Symp^{\mathcal{J}}(M)\xrightarrow{\simeq} Sp^\mathcal{J}(M)$.

From the lemma above, there is a path $\gamma(t) =(\omega_t,J_t)$ in $Symp^\mathcal{J}(M)$ from $(\omega,J)$ to $S([{\omega}^{d/2}_{\mid_M}],\omega_{\mid_M}, \lambda , J_{ker(\omega_{\mid_M})}) \in Symp^\mathcal{J}(M)$ for some characteristic $1$-form $\lambda$.  Because $Sp^\mathcal{J}$ is a fiberwise equivariant sheaf, $\gamma$ determines its `graph', $(\omega_I, J_I)\in Sp^\mathcal{J}(I\times M)$ (the notation here is explained by declaring $(\omega_I,J_I)_{(x,t)}= (\omega_t,J_t)_x)$).  The h-principle above informs us that we can integrate this (almost-complex) \textit{almost}-symplectic structure, relative to its boundary, to an (almost-complex) symplectic structure $(\omega,J)\in Symp^\mathcal{J}(I\times M)$.  This is the desired morphism from $(M,\omega,J)$ to $(M,S([{\omega}^{d/2}_{\mid_M}],\omega_{\mid_M}, \lambda , J_{ker(\omega_{\mid_M})}))$.
\end{proof}

The rest of this subsection sketches a tentative proof of Conjecture~\ref{thm: contact-symplectic} from the previous section.

\begin{lem}\label{lem: contact}
Symplectization
\[
Cont\xrightarrow{S} Symp
\]
is level-wise a homotopy equivalence.
\end{lem}

\begin{proof}[Sketch]
Similar to the proof of Lemma~\ref{lem: symplectization}, there is a map 
\[
r:Symp(M)\rightarrow Ham(M)
\]
given by $\omega\mapsto \omega_{\mid_M}$.  Similar to Lemma~\ref{lem: symplectization}, this map is a homotopy equivalence.

Let $Ham_0(M)$ be the space of \textit{exact} Hamiltonian structures on $M^{d-1}$.  For each such $M$ there is a fibration sequence
\[
Ham_0(M)\hookrightarrow Ham(M) \xrightarrow{[-]} H^2(M;\mathbb{R}).
\]
The right term being a vector space over $\mathbb{R}$ implies $Ham_0(M)\xrightarrow{\simeq}Ham(M)$ is a homotopy equivalence for each $M^{d-1}$.

Let $Ham^{char}_0(M)$ be the space of pairs $(\Omega,\lambda)$ such that $\lambda$ is a characteristic $1$-form for the \textit{exact} Hamiltonian structure $\Omega$ on $M$.  There is an obvious forgetful map $Ham^{char}_0(M)\rightarrow Ham_0(M)$ which is a fibration.  The fiber is the affine space of characteristic $1$-forms for a fixed Hamiltonian structure and is thus contractible.  There is a canonical map $Cont(M)\rightarrow Ham^{char}_0(M)$ given by $\alpha\mapsto (d\alpha, \alpha)$.  This map is a homotopy equivalence.  The result follows.
\end{proof}

\begin{conj}\label{prop: contact}
The universal map
\[
B\mathsf{Cob}^{cSymp}_d\rightarrow B\mathsf{Cob}^{Symp}
\]
is level-wise a homotopy equivalence.
\end{conj}

\begin{proof}[Sketch]
We follow the proof of Theorem~\ref{thm: SFT} above.  Namely, we show that for each germ $(\omega_0)\in Symp(M)$ there exists $(\omega)\in Symp(I\times M)$ restricting to $(\omega_0)$ on $\{0\}\times M$ and $S(\alpha)$ on $\{1\}\times M$ for some $\alpha\in Cont(M)$.  The map
\[
Symp_0(M)\hookrightarrow Symp(M)\xrightarrow{[-]} H^2(M;\mathbb{R})
\]
is a fibration sequence.  Note that $Symp_0(M)\hookrightarrow Symp(M)$ is then a homotopy equivalence.  Clearly, symplectization factors
\[
Cont\xrightarrow{S} Symp_0 \hookrightarrow Symp.
\]
It follows from Lemma~\ref{lem: contact} that $Cont\rightarrow Symp$ is level-wise a homotopy equivalence.  
Thus, there is a path $\gamma(t) =(\omega_t)$ in $Symp(M)$ from $(\omega,J)$ to $S(\alpha) \in Symp(M)$.  Because $Sp^\mathcal{J}$ is a fiberwise equivariant sheaf, $\gamma$ determines its `graph', $(\omega_I)\in Sp(I\times M)$ where $(\omega_I)_{(x,t)}= (\omega_t)_x)$).  There is an h-principle due to Gromov (see~\cite{eliashberg:h-principle} for an survey and~\cite{gromov:foliations} for a first account) stating that the inclusion $Symp_0\rightarrow Symp$ is a weak homotopy equivalence. We can therefore integrate this \textit{almost}-symplectic structure, relative to its boundary, to a symplectic structure $(\omega,J)\in Symp^\mathcal{J}(I\times M)$.  This is the desired morphism from $(M,\omega,J)$ to $(M,S(\alpha))$.
\end{proof}

\section{Gauge theoretic sheaves and smooth $4$-manifolds}\label{subsec: gauge theory}

In this section we identify the homotopy type of various gauge theoretic cobordism categories.  The current techniques for defining invariants of smooth structures, particularly on $4$-manifolds, uses gauge theory.  Namely, these invariants are characteristic classes of gauge structures on smooth manifolds, that is cohomology classes of moduli spaces $\mathcal{M}_{gauge}$ of gauge structures. Taking the point of view that we are interested primarily in characteristic classes of \textit{smooth} structures on smooth manifolds, that is cohomology classes of $BDiff$, we could ask how close the forgetful map $\mathcal{M}_{gauge} \rightarrow BDiff$ is to being an equivalence.  Appealing to \textsection\ref{subsubsec: group completion}, one interpretation of some of the results in this section is that this map is an equivalence after `stabilization'.  In this sense, the `stable' characteristic classes of gauge structures \textit{are} the `stable' characteristic classes of smooth structures.

\subsection{Flat connections}

Fix a Lie group $G$ and choose as a model for $EG$ the space of embeddings $Emb(G,\mathbb{R}^\infty)$ with the action of $G$ by pre-multiplication.  In what follows, all principal $G$-bundles are taken to be smooth.  Consider the sheaf $G\text{-}bun:= map^{sm}_{BG}$ of maps to $BG$ which classify \textit{smooth}  principal $G$-bundles in the sense of \textsection\ref{subsec: smooth families}.  Let $(P\rightarrow W)\in G\text{-}bun$ be a principal $G$-bundle.  There is the adjoint action $ad$ of $G$ on the Lie algebra $\mathcal{G}:=T_e G$ resulting in the vector bundle $ad(\mathcal{G})$ over $W$ with total space $P\times_G \mathcal{G}$.  Recall that a \textit{connection} on such a principal $G$-bundle $P\xrightarrow{\pi} W$ is a (smooth) $G$-equivariant splitting 
\[
\tau_P\xleftarrow{A} \pi^* \tau_W
\]
of the short exact sequence of vector bundles over $P$
\[
0\rightarrow\pi^*ad\mathcal{G}\rightarrow\tau_P\rightarrow \pi^* \tau_W\rightarrow 0;
\]
see~\cite{atiyah-bott:yang-mills} for a general reference on the subject.  
For a fixed principal $G$-bundle $P\xrightarrow{\pi} W$, topologize the set of (smooth) connections, denoted $Conn(P\xrightarrow{\pi} W)$, with the $C^\infty$ Whitney topology.

The data of the splitting $A$ above is equivalent to the data of a splitting $\pi^*ad\mathcal{G}\xleftarrow{\omega_A} \tau_P$.  Define the \textit{curvature $2$-form on $W$ with values in $\mathcal{G}$}, 
\[
F_A\in\Gamma(\Lambda^2\tau^*_W\otimes \mathcal{G}),
\]
by $F(v,w) = \omega_A([A(v),A(w)])$.  The connection $A$ is said to be \textit{flat} if $F_A \equiv 0$.  Notice that to say $A$ is flat is to say that the $d$-dimensional distribution $A(\pi^*\tau_W)$ on $P$ is integrable.  Thus a flat connection determines and is determined by a foliation of $P$ whose tangent distribution projects isomorphically onto $\pi^*\tau_W$.

The functor 
\[
G\text{-}conn: {\mathsf{Emb}_d}^{op}\rightarrow\mathsf{Top},
\]
given by $W\mapsto \{(p,A)\mid p\in G\text{-}bunn(W)\text{ and }A\in Conn(p)\}$, is an equivariant sheaf.  It is not obvious how to topologize the set $G$-$conn(W)$; this is done as follows.  The mapping space $Map^{sm}(W,BG)$ is locally contractible.  Consider a contractible neighborhood $f\in U\subset Map^{sm}(W,BG)$.  The choice of a contraction $c:U\simeq *$, canonically determines a morphism of principal $G$-bundles,
\[
\xymatrix{
\tilde{P} \ar[r]^{c}  \ar[d]^{\tilde{\pi}}
&
P \ar[d]^{\pi}
\\
U\times W  \ar[r]
&
W,
}
\]
where the lower horizontal map is projection, and the total space $\tilde{P}:= ev^*EG$ with 
\[
ev: U\times W\rightarrow BG
\]
the evaluation map.  Write $\tau^{vert}_{\tilde{P}}$ for the vector bundle over $\tilde{P}$ whose restriction to $\{f\}\times W$ is the tangent bundle $\tau_{f^*EG}$ of the smooth principal bundle $f^*EG$ over $W$.  Let $V\subset Map(\tilde{\pi}^*(U\times TW), \tau^{vert}_{\tilde{P}})$ be an open set.  With such a choice of contraction of $U$, there is a canonical embedding $U\times V\hookrightarrow G\text{-}conn(W)$ given by 
$(f,A)\mapsto (f,A_{\mid_f})$ where $A_{\mid_f}$ is the map 
\[
\tilde{\pi}^*(\{f\}\times TW)\rightarrow Tf^*EG
\]
whose composition with the inclusion $\tau_{f^*EG}\hookrightarrow\tau^{vert}_{\tilde{P}}$ is $A$.  Topologize $G\text{-}conn(W)$ as generated by the images of these embeddings $U\times V$ for each choice of contraction $c$ of $U$.

The obvious forgetful morphism $G\text{-}conn\xrightarrow{\simeq} G\text{-}bun$ is level-wise a homotopy equivalence  because the space $Conn(P\xrightarrow{\pi} W)$ is affine on the space of $1$-forms on $W$ and therefore contractible.  As a consequence,
\[
B\mathsf{Cob}^{G\text{-}conn}_d\simeq B\mathsf{Cob}^{G\text{-}bun}_d\simeq\Omega^{\infty-1}MTO(d)\wedge BG_+.
\]
Of more interest is the equivariant sheaf 
\[
G\text{-}flat:{\mathsf{Emb}_d}^{op}\rightarrow \mathsf{Top}
\]
of flat connections.  The theorem to follow is theorem 2 of~\cite{galatius-cohen-kitchloo:flat-connections}; presented below is an altogether different proof.

\begin{theorem}
The morphism of equivariant sheaves $G\text{-}flat\rightarrow G\text{-}bun$ induces a weak homotopy equivalence
\[
B\mathsf{Cob}^{G\text{-}flat}_d\rightarrow B\mathsf{Cob}^{G\text{-}bun}_d\simeq \Omega^{\infty-1}MTO(d)\wedge BG_+.
\]
\end{theorem}

\begin{proof}
We are interested in the stalk $Stalk_0({G\text{-}flat}_{\mid_{\mathbb{R}^d}})$.  Any principal $G$-bundle over $\mathbb{R}^d$ (induced by a map $\mathbb{R}^d\xrightarrow{p} BG$) is canonically isomorphic to the trivial bundle $p(0)\times\mathbb{R}^d\xrightarrow{pr_2}\mathbb{R}^d$ where $G\cong p(0)\in BG$.  Any foliation of $p(0)\times\mathbb{R}^d$  coming from a flat connection, when restricted to a small enough neighborhood of $0\in\mathbb{R}^d$, is standard; amounting to the data of a map of vector spaces $T_0\mathbb{R}^d\rightarrow  T_e p(0)\times T_0\mathbb{R}^d$ such that composition with the projection onto the second factor is the identity map on $T_0\mathbb{R}^d$.  Such data is contractible, contracting onto the standard inclusion $\{0\}\times T_0\mathbb{R}^d \hookrightarrow T_e p(0)\times T_0\mathbb{R}^d$.  It follows that 
\[
Stalk_0({G\text{-}flat}_{\mid_{\mathbb{R}^d}})\xrightarrow{\simeq} Stalk_0({G\text{-}bun}_{\mid_{\mathbb{R}^d}})
\]
is a weak homotopy equivalence and the result follows.
\end{proof}

\begin{remark}\label{rmk:flat-moduli}
Here, we relate the moduli space $G-flat(W)//{Diff}(W)$ of $G-flat$-structures on $W$ to similar moduli spaces found in the literature.
Let $\alpha= (E\rightarrow W)$ be a principal $G$-bundle.  Write $Aut_1(\alpha)$ for the \textit{gauge group} of $\alpha$ defined as the bundle isomorphisms
\[
\xymatrix{
E  \ar[d]  \ar[r]^{\tilde{\phi}}
&
E \ar[d]
\\
W \ar[r]^{id_W}
&
W
}
\]
which live over the identity on $W$.
The \textit{moduli space of flat connections on a principal $G$-bundle $\alpha=(E\rightarrow W^d)$} is the homotopy orbit space $G\text{-}flat(W)_\alpha//Aut_1(\alpha)$ where $G\text{-}flat(W)_\alpha \subset G\text{-}flat(W)$ is the subspace of pairs $(f,A)$ such that $f\in Map^{sm}(W,BG)_\alpha$ is in the component which classifies $\alpha$.

Write $Aut(\alpha)$ for the topological group of automorphisms
\[
\xymatrix{
E  \ar[d]  \ar[r]^{\tilde{\phi}}
&
E \ar[d]
\\
W \ar[r]^{\phi}
&
W
}
\]
which do not necessarily live over a diffeomorphism $\phi$ of $W$.  There is a semi-direct product sequence
\[
1\rightarrow Aut_1(\alpha)\rightarrow Aut(\alpha)\rightarrow {Diff}(W)\rightarrow 1.
\]
Namely, there is a section ${Diff}(W)\rightarrow Aut(\alpha)$ given by $\phi \mapsto (\phi^*,\phi)$.  This map is well-defined because our pull backs are taken inside ambient Euclidean space and are thus determined on the nose, not just up to isomorphism.
The \textit{universal} moduli space of flat connections on $\alpha$ is the homotopy quotient $G\text{-}flat(W)_\alpha//Aut(\alpha)$.  The disjoint union $\mathcal{M}^{flat}(W):= \amalg_{\alpha}( G\text{-}flat(W)_\alpha//Aut(\alpha))$ is the \textit{universal moduli space of flat connections on $W$}.

Because of the section of semi-direct product sequence above, there is a projection map $G\text{-}flat(W)//{Diff}(W) \rightarrow \mathcal{M}^{flat}(W)$ from the moduli space of $G\text{-}flat$-structures on $W$ to the universal moduli space of flat connections on $W$.   In \cite{cohen-galatius-kitchloo:flat-connections} it is shown that a natural map $\mathcal{M}^{flat}(W)\rightarrow Map^{sm}(W,BG)$ is highly connected.  
\end{remark}

\subsection{Solutions to the Yang-Mills equations for $d=4$}\label{subsubsec: yang-mills}

Let $met$ denote the equivariant sheaf of Riemannian metrics.  A  metric $<,>\in met(W)$ induces a metric on the vector bundle $\Lambda^k \tau^*_W$ whose sections are $k$-forms on the $d$-manifold $W$; to see this declare the basis $\{e^*_{i_1<...< i_k}\}$ to be orthonormal for $\{e_i\}$ orthonormal with respect to $<,>$.  Choose a metric $<,>^\prime$ on the vector bundle $ad(\pi)$ over $W$.  Such choices of metrics determine a \text{Hodge star} operator on the space of forms on $W$ with values in $ad(\pi)$ as follows.  For $\alpha\otimes s$ such a $k$-form, $*(\alpha\otimes s)$ is the $(d-k)$-form on $W$ defined implicitly by 
\[
-\wedge (*\alpha\otimes s) =< -, s>^\prime  <-,\alpha>dvol_{<,>}.  
\]
For $W$ closed, this $*$-operator demonstrates PoincarŽ duality on the level of forms.

Take $G = SU(2)$ from the section above.
Briefly, the Yang-Mills functional 
\[
\mathcal{YM}:Conn(P\xrightarrow{\pi} W)\rightarrow \mathbb{R}
\]
is defined as the integral $A\mapsto \int_W {\lVert F_A\rVert}^2 dvol_{<,>}$, referred to as the energy of a connection.  A connection $A$ is said to be a solution to the Yang-Mills equations if $A$ is a local minimum of $\mathcal{YM}$.  The Euler-Lagrange equations inform us that solutions to the Yang-Mills functional in dimension four satisfy $*F_A = \pm F_A$, the self-dual and anti-self-dual solutions depending on the sign of the Chern class $-c_2(\pi)$.  See~\cite{donaldson-kronheimer:four-manifolds} for a full account of the material.  See also~\cite{scorpan:wild-world} for an overview.

Consider the equivariant sheaf 
\[
\mathcal{YM}:\mathsf{Emb}^{op}_4\rightarrow \mathsf{Top}
\]
of (anti-)self-dual connections given by 
\[
w\mapsto \{(p,A,<,>,<,>^\prime)\}
\]
where $<,>\in met(W)$,  $p\in G\text{-}bun(W)$, $<,>^\prime\in Met(ad(p))$,  and $A\in Conn(p)$   with $*F_A=\pm F_A$.
The set $\mathcal{YM}(W)$ is topologized in an obvious way.

\begin{theorem}
The forgetful morphism of equivariant sheaves $\mathcal{YM}\rightarrow G\text{-}bun$ induces a weak homotopy equivalence
\[
B\mathsf{Cob}^\mathcal{YM}_4\rightarrow B\mathsf{Cob}^{G\text{-}bun}_4 \simeq \Omega^{\infty-1}MTO(4)\wedge BG_+.
\]
\end{theorem}

\begin{proof}
We identify the homotopy type of the stalk $Stalk_0(\mathcal{YM}_{\mid_{\mathbb{R}^4}})$.  
Begin as before by noting the canonical homotopy equivalence $G\text{-}bun(\mathbb{R}^4)\simeq * = \{G\times\mathbb{R}^4\xrightarrow{pr_2}\mathbb{R}^4\}$.  Such a contraction pulls-back the remaining data $\{A,<,>,<,>^\prime\}$.  It can thus be assumed that $p\in G\text{-}bun(\mathbb{R}^4)$ is constant.  There is a similar contraction of $met(\mathbb{R}^d)$ onto $\delta_{ij}\in met(\mathbb{R}^4)$.  The splitting of $p = (G\times\mathbb{R}^4\rightarrow\mathbb{R}^4)$ results in a preferred connection 
\[
\tau_G\times\tau_{\mathbb{R}^4}\xleftarrow{A_0} \pi^*\tau_{\mathbb{R}^4}
\]
given by $A_0(g,(m,v))= ((g,0),(m,v))$.  There is a straight-line homotopy $A_t:= tA + (1-t) A_0$ of connections from $A$ to $A_0$.  A calculation using the triviality of $A_0$ shows  $F_{A_t} = t^2 F_A$ which is clearly (anti-)self-dual exactly when $F_A$ is (anti-)self-dual.  Lastly, because the space of metrics $<,>^\prime$ is affine and thus contractible, there is a homotopy equivalence
\[
Stalk_0(\mathcal{YM}_{\mid_{\mathbb{R}^4}})\xrightarrow{\simeq} Stalk_0({G\text{-}bun}_{\mid_{\mathbb{R}^4}})
\]
given by $(p,A,<,>,<,>^\prime)\mapsto p$.  The result follows from the main theorem.
\end{proof}

\begin{remark}
We relate the moduli space $\mathcal{YM}(W)$ of Yang-Mills structures on $W$ to those found in the literature, namely the moduli space of (anti-) self-dual connections on $W$ as in Donaldson's theory.  The procedure is identical to that of the above Remark~\ref{rmk:flat-moduli} with $\mathcal{YM}$ in place of $G\text{-}flat$.  
\end{remark}

\subsection{Solutions to the Seiberg-Witten equations}

This section begins with a refresher on Seiberg-Witten equations.  The reader is referred to~\cite{morgan:seiberg-witten} for a full introductory account and to~\cite{scorpan:wild-world} for an overview.

The fundamental group $\pi_1 SO(d)\cong \mathbb{Z}/2$ for $d>2$ with universal double-cover group written $Spin(d)$.  The group $\mathbb{Z}/2$ being normal in both $U(1)$ and $Spin(4)$, consider the Lie group 
\[
Spin^\mathbb{C}(d):= U(1)\times_{\mathbb{Z}/2} Spin(d).
\]
For $d=4$, via quaternionic trickery, there is a splitting $Spin(4)\cong SU(2)\times SU(2)$.  There are then three representations 
\[
Spin^\mathbb{C}(4)\xrightarrow{det} U(1)\hspace{.5cm}\text{ and } \hspace{.5cm}  Spin^\mathbb{C}(4)\xrightarrow{\rho_\pm} U(2)
\]
classifying a complex line bundle, called the determinate bundle, and two complex plane bundles.  Think of these complex plane bundles as quaternionic line bundles tensored with the ``square root" of this determinate bundle.

Let $W^4\subset\mathbb{R}^\infty$ be a smooth embedded $4$-manifold with induced Riemannian metric $<,>$.  Suppose $W$ has a $spin^\mathbb{C}$-structure $l$ as in the diagram
\[
\xymatrix{
&
BSpin^\mathbb{C}(4)  \ar[d]
\\
W  \ar[ur]^l  \ar[r]^{\tau_W}
&
BO(4).
}
\]
There are then the determinant complex line bundle $det$ and two complex plane bundles $\omega^\pm$ over $W$.  Together with the Levi-Civita connection on $(W,<,>)$, a connection $A$ on the principal $U(1)$-bundle $W\xrightarrow{det}BU(1)$ induces a connection on the bundles $(P_{l}\rightarrow W)$ and $(P_{\omega^\pm}\rightarrow W)$ still referred to as $A$.  Because the Lie algebra of $U(1)$ is $\mathbb{R}$, the curvature $2$-form $F_A$ on $W$ from $det$ has values in $\mathbb{R}$.

There is a bundle morphism
\begin{equation}\label{eq: clifford}
\omega^\pm\otimes  \tau_W \xrightarrow{\bullet} \omega^\mp
\end{equation}
given locally by quaternionic multiplication $(w,v)\mapsto w*v$.  This is called the \textit{Clifford} action.  Iterating this Clifford action, there is an obvious extension of the Clifford action to all tensor powers of $\tau_W$.  The same is true for tensor powers of $\tau^*_W$ because of the metric $<,>$ on $W$ exhibiting an isomorphism $\tau_W\cong\tau^*_W$.  In particular, it is possible to regard the curvature $2$-form $F_A$ of $det$ as a section of the endomorphism bundle
\[
F_A \in \Gamma(End_\mathbb{C}(\omega^+)).
\]
There is a morphism of bundles
\[
\sigma: \omega^+\rightarrow End_\mathbb{C}(\omega^+),
\]
called the \textit{squaring map}, given when $w$ has unit length by $w\mapsto Proj_{Span_\mathbb{C} \{w\}}-(1/2) id $, the projection onto the span of $w$ minus its trace.

For $\alpha=(E\rightarrow B)$ a general rank $n$ vector bundle with structure group $G$, write $P_\alpha$ for the total space of its associated principal bundle.  So $(P_\alpha\times_G \mathbb{R}^n\rightarrow B) \cong \alpha$.  Using this equivalence, there is an assignment 
\[
\Gamma(\alpha)\rightarrow Map_G(P_\alpha, \mathbb{R}^n)
\]
from the space of sections to the space of $G$-equivariant maps given by $(b\mapsto(p,v))\mapsto(p\mapsto v)$.  This assignment is a homeomorphism.  In this way it is possible to view a section of $\alpha$ as $n$ real valued functions on $P_\alpha$.

In view of the paragraph above, there is a map
\[
\Gamma(\omega^+ )\xrightarrow{\nabla^A} \Gamma(\omega^+ \otimes  \tau^*_W)
\]
given by $s\mapsto (X\mapsto A(X)(s))$, the latter being the the vector $A(X)$ on $P_{\omega^+}$ acting on the coordinate functions of $s$.  This is the familiar notion of covariant derivative induced by a connection.  Using the isomorphism $\tau^*_W\cong \tau_W$ induced from the metric $<,>$, and composing with the action~(\ref{eq: clifford}), we obtain the map
\[
D^A:\Gamma(\omega^+)\rightarrow \Gamma(\omega^+\otimes \tau_W)\rightarrow \Gamma(\omega^-)
\]
called the \textit{Dirac operator}.  Locally, for $\{e_i\}$ a local orthonormal frame, the Dirac operator is expressed as 
\[
\Sigma_i e_i\bullet \nabla^A_{e_i} 
\]
where $\bullet$ denotes the Clifford action.

A \textit{Seiberg-Witten structure} on $W^4$ is the data $(l,<,>, A,s)$ where $l$ is a $Spin^\mathbb{C}(4)$-structure on $W$, $<,>$ is a metric on $W$, $A$ is a connection on the associated principal $U(1)$-bundle, and $s\in \Gamma(\omega^+)$ is a section of the induced complex plane bundle from $l$ such that
\begin{equation}\label{eq: seiberg}
D^A s = 0\hspace{.5cm}\text{ and }\hspace{.5cm} F_A = \sigma(s)\in End(\omega^+).
\end{equation}
There is an obvious topology on the set of Seiberg-Witten structures on $W$ whose resulting space is  denoted 
\[
\mathcal{SW}(W) =\{(l,A,s)\mid~(\ref{eq: seiberg})\}.
\]
The local character of a Seiberg-Witten structure results in $\mathcal{SW}$ being an equivariant sheaf.

\begin{theorem}
The forgetful morphism of equivariant sheaves $\mathcal{SW}\rightarrow Spin^\mathbb{C}(4)$ induces a weak homotopy equivalence
\[
B\mathsf{Cob}^{\mathcal{SW}}_4 \simeq \Omega^{\infty-1}MTSpin^\mathbb{C}(4).
\]
\end{theorem}

\begin{proof}
We analyze $\mathcal{SW}(\mathbb{R}^4)$ in a neighborhood of the origin.  There is a canonical isomorphism (and thus path) from any $Spin^\mathbb{C}(4)$-structure $l:\mathbb{R}^4\rightarrow BSpin^\mathbb{C}(4)$ to the constant structure $l(0)$ for which all of the bundles $\omega^\pm$ and $det$ are trivial.  Independently, there is a canonical path of metrics on $\mathbb{R}^4$ from $<,>$ to the standard metric $\delta_{ij}$ on $\mathbb{R}^4$.  Under these simplifications, because the connection on $\omega^+$ is induced from that on $det$, the Dirac operator $D^A$ is independent of $A$.  The splitting of $P_{det} = U(1)\times\mathbb{R}^4$ over $\mathbb{R}^4$ results in a preferred connection 
\[
\tau_{U(1)}\times\tau_{\mathbb{R}^4}\xleftarrow{A_0} \pi^*\tau_{\mathbb{R}^4}
\]
given by $A_0(\lambda,(m,v))= ((\lambda,0),(m,v))$.  There is a straight-line homotopy $A_t:= tA + (1-t) A_0$ of connections from $A$ to $A_0$.  A calculation using the triviality of $A_0$ shows  $F_{A_t} = t^2 F_A$.  Similarly, the squaring map scales as $\sigma(ts) = t^2 s$.  Thus, for $(A,s)\in \mathcal{SW}(\mathbb{R}^4)$  there is a canonical path from $(A,s)$ to $(A_0,0)$ demonstrated by $(A_t, t s)$, see~(\ref{eq: seiberg}).    This proves that the assignment 
\[
\mathcal{SW}(\mathbb{R}^4)\rightarrow Spin^\mathbb{C}(\mathbb{R}^4),
\]
given by $(l,<,>,A,s)\mapsto l$, is a homotopy equivalence on stalks.  The result follows from the main theorem.
\end{proof}

\begin{remark}
We relate the moduli space $\mathcal{SW}(W)$ of Seiberg-Witten structures on $W$ to those found in the literature.  Once again, the procedure is identical to that of the above Remark~\ref{rmk:flat-moduli} with $\mathcal{SW}$ in place of $G\text{-}flat$ and the gauge group $Gauge(det) \cong Map^{sm}(W,U(1))$ in place $Aut(\alpha)$.  
\end{remark}

\section{Configurations, $2$-dimensional oriented orbifolds, marked holomorphic curves, and submanifolds}\label{subsec: configurations}

In what follows we will make use of the following well-known proposition.
\begin{prop}\label{prop: cofibration sequence}
Let $\alpha=(E\rightarrow B)$ be an $O(n)$-vector bundle with associated unit sphere bundle $S(E)\xrightarrow{p} B$.  Let $\alpha^\prime=(E^\prime\rightarrow B)$ be an $O(k)$-vector bundle.  There is a natural cofibration sequence of spectra
\[
Th(p^*\alpha^\prime)\xrightarrow{\overline{p}} Th(\alpha^\prime) \rightarrow Th(\alpha\oplus \alpha^\prime).
\]
\end{prop}
\begin{proof}
The argument is simple.  For $u\in S(E_b)$, $v\in E^\prime_b$, and $t\in [0,\infty]$, the assignment $(u,v,t)\mapsto (tu,v)$ describes a map of spectra from the mapping cone of $\overline{p}$ to $Th(\alpha\oplus \alpha^\prime)$.  This map is a homotopy equivalence.
\end{proof}

\subsection{Configurations}

Consider the equivariant sheaf
\[
con: {\mathsf{Emb}_d}^{op}\rightarrow \mathsf{Top}
\]
described by $con(U)=\Psi_0(U)$.  So a point in $con(U)$ is a discrete subset of $U$.  The work has already been done in~\S\ref{subsec: homotopy type of psi} (there when $d=0$) to establish $Stalk_0(con_{\mid_{\mathbb{R}^d}})\simeq S^d=\mathbb{R}^d\cup\{\infty\}$.  The homotopy orbit $S^d//O(d)$ is the  sphere bundle of the vector bundle $\epsilon\oplus \gamma_{d-1}$ over $BO(d)$.  
As such, for $p:S^d//O(d)\rightarrow BO(d)$ the projection, the Thom spectrum $(S^d//O(d))^{-\gamma_d}$ fits into the cofibration sequence of spectra
\[
(S^d//O(d))^{-\gamma_d}\rightarrow BO(d)^{-\gamma_d}\rightarrow BO(d)^{\epsilon\oplus \gamma_d\oplus (-\gamma_d)}.
\]
The right term is the suspension spectrum $\Sigma^{\infty+1} BO(d)_+$ while the middle term is by definition $MTO(d)$.  There is a canonical null homotopy of the right map due to the factor $\epsilon$.  It follows that 
\[
(S^d//O(d))^{-\gamma_d} \simeq MTO(d)\vee \Sigma^{\infty} BO(d)_+.
\]
We conclude from the main theorem and this discussion that 
\[
B\mathsf{Cob}^{con}_d\simeq   \Omega^{\infty-1}  MTO(d) \times Q \Sigma BO(d)_+.
\]

It can be checked that the forgetful morphism of equivariant sheaves $con\rightarrow *$ induces the map of spectra which simply collapses the term $\Sigma^{\infty} BO(d)_+$.  The morphism of equivariant sheaves $*\rightarrow con$, with image the empty configuration, describes the splitting.

Similarly, there is a functor $\mathsf{Cob}^{con}_d\rightarrow \mathsf{Cob}^{map_{BO(d)}}_0$ given by $(W,K)\mapsto (K,{\tau_W}_{\mid_K})$ where $K\subset W$ is a configuration.  This functor induces the map of spectra which collapses the $MTO(d)$ term.  The splitting of this term is related to Euler characteristic in the following way.

Write $S(\gamma_{d,N}) = S(U_{d,N})\xrightarrow{p} Gr_{d,N})$ for the unit sphere bundle of $\gamma_{d,N}=(U_{d,N}\rightarrow Gr_{d,N})$.  There is a map $Gr_{d-1,N-1}\xrightarrow{\iota} S(U_{d,N})$ given by $V\mapsto (\mathbb{R}\times V, e_1)$.  The connectivity of this map grows with $N$ and realizes a homotopy equivalence in the colimit over $N$.  Moreover, the canonical morphism
\[
\gamma^\perp_{d-1,N-1}\rightarrow \iota^*\gamma^\perp_{d,N}
\]
of rank $(N-d)$ vector bundles over $Gr_{d-1,N-1}$ is an isomorphism.  Thus, in the limit $\Sigma^{-1} MTO(d-1) = \Sigma^{-1}Th(-\gamma_{d-1})\xrightarrow{\simeq} Th(-p^*\gamma_d)$ is a homotopy equivalence.
Upon identifying $Th(\gamma_d\oplus -\gamma_d) =  \Sigma^\infty BO(d)_+$, Proposition~\ref{prop: cofibration sequence} gives the cofibration sequence of Thom spectra
\[
\Sigma^{-1} MTO(d-1)\rightarrow MTO(d)\rightarrow \Sigma^{\infty} BO(d)_+.
\]
Presented below is a geometric interpretation of the right map.

Recall the space $\Psi^\mathcal{F}_d(\mathbb{R}^N)$ of pairs $(W,g)$ where $W\subset\mathbb{R}^N$ is an embedded $d$-manifold and $g\in \mathcal{F}(W)$.  Although not induced by a map of equivariant sheaves, there is a homotopy class of a map
\[
\Psi_d \rightarrow \Psi^{con}_d
\]
described as follows.  Over $\Psi_d(\mathbb{R}^N)$ is the tautological space 
\[
E\Psi_d(\mathbb{R}^N):=\{(W,w)\mid w\in W\}\subset\Psi_d(\mathbb{R}^N)\times \mathbb{R}^N.  
\]
Choose a continuous section $s\in\Gamma(E\Psi_k(\mathbb{R}^N)\rightarrow \Psi_k(\mathbb{R}^N))$ so that $s(W)$ represents the Euler characteristic of the $d$-manifold $W$.  Such a continuous section exists on compact smooth families $X\xrightarrow{f} \Psi_d(\mathbb{R}^N)$  as the (isolated) zeros of a vertical vector field on the pull back $f^* E\Psi_d(\mathbb{R}^N)$ .  Such a section is a lift 
\[
\xymatrix{
&
\Psi^{con}_d(\mathbb{R}^N)  \ar[d]
\\
X  \ar[r]^f  \ar[ur]^s
&
\Psi_d(\mathbb{R}^N).
}
\]
Composing $s$ with the projection 
$\Psi^{con}_d(\mathbb{R}^N)\rightarrow \Psi^{map_{BO(d)}}_0(\mathbb{R}^N)$
described above \textit{realizes} the quotient map $MTO(d)\rightarrow \Sigma^{\infty}BO(d)_+$.

\subsection{$2$-dimensional oriented orbifolds}

A $2$-dimensional oriented orbifold is the data of an oriented surface with unordered marked discrete points each labeled by a cyclic group.  Consider then the equivariant sheaf
\[
\mathbb{N}^{con}: \mathsf{Emb}^{op}_2\rightarrow   \mathsf{Top}
\]
given by 
\[
\mathbb{N}^{con}(U) = \{(\sigma,K,q)\mid\sigma\in \mathcal{O}r(U),\text{ } K\in con(U),\text{ and }q:K\rightarrow \mathbb{N}\}.
\]  
It is simple to check that $Stalk_0({\mathbb{N}^{con}}_{\mid_{\mathbb{R}^2}}) \simeq \mathbb{N}\times S^2$ and thus 
\[
B\mathsf{Cob}^{\mathbb{N}^{con}}_2\simeq   \Omega^{\infty-1}  (\mathbb{N}\times  S^2//SO(2))^{-\gamma_2}.
\]
The observation that there is a bijection between $mor\text{ }\mathsf{Cob}^{\mathbb{N}^{con}}_2$ and the space of $2$-dimensional embedded oriented orbifolds in $\mathbb{R}^\infty$ elects $\mathsf{Cob}^{\mathbb{N}^{con}}_2$ as a fine example of a $2$-dimensional oriented cobordism category.  One might request that the objects of such an orbifold cobordism category be either smooth manifolds themselves or at least orbifolds in their own right.  In either case, as all objects in the stated orbifold cobordism category are \textit{homotopy isomorphic} to such objects, from~\textsection\ref{subsec: Restricting objects} the weak homotopy of the classifying spaces of these various candidate orbifold cobordism categories will be identical.  In this ad hoc way, we have thus identified its weak homotopy type.

\subsection{Marked maps}

Fix a smooth closed manifold $Y$ and a compact submanifold $Z\subset Y$.  Choose a complete metric on $Y$.  Choose $\delta>0$ such that a $\delta$-neighborhood of $Z$ in $Y$ is a tubular neighborhood.  Refer to such a tubular neighborhood as $\nu_Z$.

Consider the equivariant sheaf
$con_{Y,Z}$ given by 
\[
con_{Y,Z}(U)=\{(K,f)\mid K\in con(U)\text{ and } f\in map_Y(U)\text{ such that } f(K)\subset Z\}.
\]
Consider the diagram 
\[
Y\leftarrow \nu_Z\times S^{d-1}\rightarrow Z
\]
where the left arrow is projection onto $\nu_Z$ followed by inclusion into $Y$; the right arrow is projection onto $\nu_Z$ followed by projection onto $Z$.  With trivial $O(d)$-actions on the left and right terms, the maps in the diagram are $O(d)$-equivariant.  Write $H$ for the homotopy colimit of the diagram with the induced $O(d)$-action.

\begin{theorem}\label{thm:configurations-maps}
There is a weak homotopy equivalence
\[
B\mathsf{Cob}^{con_{Y,Z}}_d \simeq \Omega^{\infty-1} ( H//O(d))^{-\gamma_d}.
\]
\end{theorem}

\begin{remark}
Let $Y=(Y,i)$ be a complex (or symplectic) manifold.  Though no details will be provided, the proof below holds for the equivariant sheaf $(con\times\mathbb{C})_{Y,Z}$ described by 
\[
(con\times\mathbb{C})_{Y,Z}(W) =\{(K,j,f)\mid K\in con(W) \text{  and } f\in \mathbb{C}_Y(W) \text{ with }f(K)\subset Z\}.
\]
(Recall that $\mathbb{C}_Y$ is the sheaf of complex structures along with a holomorphic map to $Y$).  
In the case $d=2$, the cobordism category $\mathsf{Cob}^{(con\times\mathbb{C})_{Y,Z}}_2$ is the category of \textit{marked} holomorphic curves with markings landing in a prescribed submanifold $Z$ of $Y$.  The space of morphisms of this category consisting of closed surfaces is homotopy equivalent to the moduli spaces in enumerative geometry and Gromov-Witten theory.
\end{remark}

Before presenting the proof, we make $H$ more explicit.   There is a morphism of horizontal $O(d)$-diagrams
\[
\xymatrix{
(\nu_Z  \ar[d]
&
\nu_Z\times S^{d-1} \ar[l]   \ar[r]  \ar[d] 
&
Y)  \ar[d]
\\
(\{0\}
&
S^{d-1}  \ar[l]  \ar[r]
&
\{\infty\})
}
\]
inducing an $O(d)$-equivariant map on homotopy colimits as $H\rightarrow S^d\cong \mathbb{R}^d\cup\{\infty\}$.  Away from $\infty$, this map is a fibration.  The fiber over $x\in \mathbb{R}^d$ is the space 
\[
\{p\in \nu_Z\mid dist(p,pr_Z(p))\leq\phi( \lVert x \rVert)  \}\simeq\nu_Z
\]
where $\phi:[0,\infty) \rightarrow [0,\delta)$ is some homeomorphism.  The fiber over $\infty$ is of course $Y$.

\begin{proof}[Proof of Theorem~\ref{thm:configurations-maps}]

We will build the diagram
\[
\xymatrix{
(V_0  \ar[d]^h
&
V_0\cap V_{min}  \ar[l]  \ar[r]  \ar[d]^h
&
V_{min})  \ar[d]^h
\\
(Y 
&
\nu_Z \times S^{d-1}    \ar[l]    \ar[r]
&
Z).
}
\]
Consider the subspace $con^\delta_{Y,Z}(\mathbb{R}^d)$ consisting of those pairs $(K,f)$ for which the image $f(\mathbb{R}^d)\subset B_\delta(f(0))$ is contained in a $\delta$-neighborhood of $f(0)\in Y$.  Recognize $con^\delta_{Y,Z}(\mathbb{R}^d) = V_0 \cap V_{min}$ as in the proof from~\S\ref{subsec: homotopy type of psi} where $V_0$ consists of those pairs $(K,f)$ for which $0\notin K\subset\mathbb{R}^d$ and $V_{min}$ consists of those pairs for which $K$ has a unique closest point to $0\in\mathbb{R}^d$.  The top horizontal maps are inclusions of subspaces.
The left vertical map is given by $h(K,f) = f(0)\in Y$.  The right vertical map is given by $h(K,f) = f(k)$ where $k\in K$ is the unique closest point in $\mathbb{R}^d$ to $0$.  Finally, the middle vertical map is given by $h(K,f) = (f(0) , k/\lVert k \rVert)$.  There are obvious actions of $O(d)$ on the spaces in the diagram (some of them trivial actions) with respect to which the diagram is $O(d)$-equivariant.

The diagram commutes up to a specific homotopy.  Namely, the left square commutes up to the homotopy $((K,f),t)\mapsto f(tk)$.  But more, each vertical map is a weak homotopy equivalence; the argument being similar to those in~\S\ref{subsec: homotopy type of psi}.  Regarding the diagram as a morphism of horizontal diagrams, it follows that the map of homotopy colimits of the horizontal diagrams is a weak homotopy equivalence.  Thus, denoting the homotopy colimit of the bottom horizontal diagram by $H$, the obvious zig-zag map $con^\delta_{Y,Z}(\mathbb{R}^d)\xrightarrow{\simeq} H$ is a weak homotopy equivalence.

Taking germs of such maps as $\mathbb{R}^d$ is replaced by $\epsilon$-neighborhoods of the origin, we arrive at the $O(d)$-equivariant zig-zag of weak homotopy equivalences
\[
Stalk_0({con_{Y,Z}}_{\mid_{\mathbb{R}^d}}) \xrightarrow{\simeq} H.  
\]
Indeed, the inclusion 
\[
colim_\epsilon con^\delta_{Y,Z}(D_\epsilon) \rightarrow colim_\epsilon con_{Y,Z}(\mathbb{R}^d) = Stalk_0({con_{Y,Z}}_{\mid_{\mathbb{R}^d}})
\]
is a weak homotopy equivalence for the following reason.  For $K$ a compact CW-complex, and $K\rightarrow con_{Y,Z}(\mathbb{R}^d) $ a family of such configurations, there is an $\epsilon>0$ such that the composition $K\rightarrow con_{Y,Z}(\mathbb{R}^d) \xrightarrow{restrict} con_{Y,Z}(D_\epsilon)$ canonically factors as $K\rightarrow con^\delta_{Y,Z}(D_\epsilon)\hookrightarrow con_{Y,Z}(D_\epsilon)$.
The result follows from the main theorem.  
\end{proof}

\subsection{Higher dimensional configurations: $\Psi_k$}

Consider the equivariant sheaf $\Psi_k:{\mathsf{Emb}_d}^{op}\rightarrow \mathsf{Top}$ from \textsection\ref{subsec: the sheaf} given as 
\[
\Psi_k(W)=\{K^k\subset W^d\}
\]
where $K\subset W$ is a $k$-dimensional submanifold of $W^d$ without boundary and closed as a subset.  From the main theorem there is a weak homotopy equivalence
\[
B\mathsf{Cob}^{\Psi_k}_d \simeq \Omega^{\infty-1}(Th(\gamma^\perp_{k,d})//O(d))^{-\gamma_d}
\]
since from \textsection\ref{subsec: homotopy type of psi}, $\Psi_k(\mathbb{R}^d)\simeq Th(\gamma^\perp_{k,d})$ as $O(d)$-spaces.  Here, the action of $O(d)$ on $Th(\gamma^\perp_{k,d})$ comes from the standard action of $O(d)$ on $\mathbb{R}^d$.

Consider the fiberwise equivariant sheaf $map^\tau_{BO(d)}:\mathsf{Emb}^{op}_k \rightarrow \mathsf{Top}$ associated to the following fibration over $BO(k)$.  Define 
\[
B_\tau = \{(U,V)\mid U\in BO(k)\text{ and }U\subset V\in BO(d)\}.
\]
The obvious projection of $B_\tau$ onto $BO(k)$ is a fibration.  The fiber over $U\in BO(k)$ is the space $BO(d-k)$.  In fact, there is a splitting $B_\tau \xrightarrow{\simeq} BO(k)\times BO(d-k)$ with homotopy inverse given by direct sum.  
There are two forgetful functors
\[
\mathsf{Cob}_d \leftarrow \mathsf{Cob}^{\Psi_k}_d\rightarrow \mathsf{Cob}^{map^\tau_{BO(d)}}_k.
\]
The left arrow is given by $(K^k\subset W^d)\mapsto W^d$ while the right arrow is given by $(K\subset W)\mapsto (K,{\tau_W}_{\mid_K})$.  Upon taking classifying spaces there results a map of spectra
\[
(\Psi_k(\mathbb{R}^d)//O(d))^{-\gamma_d}\rightarrow BO(d)^{-\gamma_d}\vee (BO(k)^{-\gamma_k} \wedge BO(d-k)_+).
\]

Let $X^r$ be a compact $r$-manifold and let $X\xrightarrow{f} \Psi_d(\mathbb{R}^N)$ be a smooth family. After possibly increasing $N$, a smooth lift $X\xrightarrow{\tilde{f}} \Psi^{\Psi_k}_d(\mathbb{R}^N)$ of $f$ yields a smooth family $X\rightarrow \Psi^{map_{BO(d-k)}}_k(\mathbb{R}^N)$.  For $k=1$, there is always such a lift given by `continuously' choosing fiberwise representatives of the dual of $w_{d-1}$, the $(d-1)^{st}$ Stiefel-Whitney class.  The resulting composition $X\xrightarrow{w_{d-1} f}  \Psi^{map_{BO(d-1)}}_1(\mathbb{R}^N)$ being null homotopic implies the family $X\xrightarrow{f} \Psi_d(\mathbb{R}^N)$ is concordant to a composition 
\[
X\rightarrow \Psi_{d-2}(\mathbb{R}^{N-2})\xrightarrow{\mathbb{R}^2\times -} \Psi_d(\mathbb{R}^N).
\]  
The same is true for general $k\leq d$ provided it is possible to choose continuous fiberwise \textit{manifold} representatives of the dual of $w_{d-k}$.  In this way, $MTO(k)\wedge BO(d-k)_+$ is a close description of the space of obstructions for a map
\[
X\rightarrow MTO(d)
\]
to lift to $\Sigma^{-k-1}MTO(d-k-1)$.

We now describe the actual cofiber of the map $\Sigma^{-k}MTO(d-k)\rightarrow MTO(d)$.  
Write $Fr_k(\gamma_{d,N})=(Fr_k(U_{d,N})\xrightarrow{p} Gr_{d,N})$ for the \textit{unit} $k$-frame bundle associated to $\gamma_{d,N}=(U_{d,N}\rightarrow Gr_{d,N})$.  
There is a map $Gr_{d-k,N-k}\xrightarrow{\iota} Fr_k(U_{d,N})$ given by $V^{d-k}\mapsto (\mathbb{R}^k\times V, \{e_1,...,e_k\})$.  The connectivity of this map grows with $N$ and realizes a homotopy equivalence in the colimit over $N$.  Moreover, the canonical morphism
\[
\gamma^\perp_{d-k,N-k}\rightarrow \iota^*\gamma^\perp_{d,N}
\]
of rank $(N-d)$ vector bundles over $Gr_{d-k,N-k}$ is an isomorphism.  Thus, in the limit $\Sigma^{-k}MTO(d-k) = \Sigma^{-k}Th(-\gamma_{d-k})\xrightarrow{\simeq} Th(-p^*\gamma_d)$ is a homotopy equivalence.

Let $\alpha:=(P\rightarrow B)$ be a principal $O(n)$-bundle.  Denote also by $\alpha$ its asociated rank $n$ vector bundle with metric.   For $V$ a vector space with metric $<,>$, a collection $\{v_1,...,v_k\}\subset V$ is said to be a \textit{scaling} $k$-frame of $V$ if there is an $r\geq 0$ such that $<v_i,v_j> = r^2\delta_{ij}$.  Let $Fr^+_k(V)$ be the space of \textit{scaling} $k$-frames in $V$ topologized in the obvious way.  Define $Fr^+_k(\alpha) =(P\times_{O(n)} Fr^+_k(\mathbb{R}^n)\rightarrow B) $ as the bundle of \textit{scaling} $k$-frames of $\alpha$.  In this way there is the fiber bundle
\[
Fr^+_k(\gamma_{d,N}) = (Fr^+_k(U_{d,N}) \rightarrow Gr_{d,N}).
\]
The projection 
\[
Fr^+_k(U_{d,N})\times_{Gr_{d,N}}U^\perp_{d,N} \rightarrow Gr_{d,N}
\]
is a fiber bundle over $Gr_{d,N}$ denoted by $Fr^+_k(\gamma_{d,N})\oplus \gamma^\perp_{d,N}$.  Extend notation and denote by $Th(Fr^+_k(\gamma_{d,N})\oplus \gamma^\perp_{d,N})$ the one-point compactification of $Fr^+_k(U_{d,N})\times_{Gr_{d,N}}U^\perp_{d,N}$.  Write $TF_kO(d):= Th(Fr^+_k(\gamma_d)\oplus -\gamma_d)$ for the spectrum whose $N^{th}$ space is $Th(Fr^+_k(\gamma_{d,N})\oplus \gamma^\perp_{d,N})$ with structures maps induced from those of the Thom spectrum $Th(-\gamma_d)$.  A nearly identical proof to Proposition~\ref{prop: cofibration sequence} proves

\begin{prop}\label{prop:scaling-frames}
There is a cofibration sequence of spectra
\[
\Sigma^{-k} MTO(d-k)\rightarrow MTO(d)\xrightarrow{\chi_k} TF_kO(d).
\]
\end{prop}

It follows that $\Omega^\infty TF_kO(d)$ measures obstructions for a family $X\rightarrow \Omega^\infty MTO(d)$ to be concordant to a family $X\rightarrow \Omega^{\infty+k} MTO(d-k)\rightarrow \Omega^\infty MTO(d)$.  Said another way, $\Omega^\infty TF_kO(d)$ measures the obstruction for a stable family of $d$-manifolds to be concordant to a stable family of $(d-k)$-manifolds fiberwise crossed with $\mathbb{R}^k$.

As should be expected, the spectra $TF_kO(d)$ are built inductively out of the $TF_1O(l) = Th(\gamma_l\oplus-\gamma_l)\cong\Sigma^\infty BO(l)_+$ with $l\leq d$. 

\begin{prop}\label{prop:scaling-frames2}
There is a cofibration sequence of spectra
\[
\Sigma^{-1} TF_{k-1}O(d-1)\rightarrow TF_kO(d) \rightarrow \Sigma^\infty BO(d)_+.
\]
\end{prop}

\begin{proof}
It is routine to verify that this is indeed a cofibration sequence.  For concreteness, the left map is given by assigning to a scaling $(k-1)$-frame $\{v_1,...,v_{k-1}\}\subset V^{d-1}$ the scaling $k$-frame $\{e_1, v_1,...,v_{k-1}\}\subset \mathbb{R}\oplus V $.  The right map is given by assigning to a scaling $k$-frame $\{v_1,...,v_k\}\subset V^d$ the scaling $1$-frame $\{v_1\}\subset V$.  
\end{proof}

\begin{remark}
There are obvious analogues to the previous two propositions \ref{prop:scaling-frames} and \ref{prop:scaling-frames2} for vector spaces endowed with tangential structure.  Notable examples of such extra structure are that of an orientation, replacing the groups $O(-)$  with $SO(-)$ in all places; and framings, replacing the symbol $O(-)$ with the symbol $1(-)$.
\end{remark}

\begin{prop}
Suppose $W^d$ is a parallelizable $d$-manifold.  Then the map induced by inclusion of $1$-simplicies
\[
B{Diff}(W)\xrightarrow{\iota} B\mathsf{Cob}_d \simeq \Omega^\infty MTO(d)
\]
is the zero-map in rational homotopy.  Choosing an orientation of $W$, the same is true for
\[
BDiff^{\mathcal{O}r}(W) \xrightarrow{\iota} B\mathsf{Cob}^{\mathcal{O}r}_d \simeq \Omega^\infty MTSO(d).
\]
\end{prop}

\begin{proof}
Only the non-oriented case will be presented, the oriented case being nearly identical.  The case $d=0$ is obviously true.  Assume $d\geq 1$.  First note from the fibration sequence
\[
\Omega^{\infty+d}MTO(0) \rightarrow \Omega^\infty MTO(d) \xrightarrow{\chi_d} \Omega^\infty TF_dO(d)
\]
that the right map induces an isomorphism on rational homotopy groups.  Next, notice by induction on $d$ that Proposition~\ref{prop:scaling-frames2} implies 
\[
\mathbb{Q}\cong \mathbb{Q}\otimes \pi_q \Omega^\infty TF_d1(d).
\]
Let $W\subset\mathbb{R}^\infty$ be the base point in $B{Diff}(W)$.  So to prove the corollary we need only show that any smooth based map 
\[
S^q\xrightarrow{f} B{Diff}(W)\xrightarrow{\iota} \Omega^\infty MTO(d) \xrightarrow{\chi_d} \Omega^\infty TF_dO(d)
\]
factors through $\Omega^\infty TF_d1(d) \rightarrow \Omega^\infty TF_dO(d)$ up to homotopy.  For this we describe an explicit construction for the composite map $S^q \xrightarrow{\chi_d \circ \iota \circ f}  \Omega^\infty TF_dO(d)$.

A smooth based map $S^q\xrightarrow{f} B{Diff}(W)$ is the data of a smooth fiber bundle $E^{q+d}\xrightarrow{p} S^q$ with fiber over the base point $*\in S^q$ \textit{equal} to $W$.  Over $E$ there is the vertical tangent bundle $\tau^v_p$ whose fiber over $e\in E$ is the vertical tangent space $T_e {p^{-1}(p(e))}$.  Associated to $\tau^v_p$ is the fiber bundle $Fr^+_d(\tau^v_p)$ over $E$ of scaling $d$-frames of $\tau^v_p$.  The zero-section of this fiber bundle gives the composite map $S^q \rightarrow \Omega^\infty TF_dO(d)$.

Unfortunately, the map $B{Diff}(W)\xrightarrow{\iota} \Omega^\infty MTO(d)$ is not based.  In fact, this map does not necessarily land in the component of the base point.  To resolve this issue, we choose the image of $W\in B{Diff}(W)$ as the base point in the appropriate component of $\Omega^\infty MTO(d)$.  To ensure that the zero-section map $B{Diff}(W) \rightarrow \Omega^\infty MTO(d)\rightarrow \Omega^\infty TF_dO(d)$ is a based map, we describe once and for all a path from the image of $W$ to the base point $\emptyset\in \Omega^\infty TF_dO(d)$ of the empty frame.  Fix a framing of $W$.  Scalar-multiplying this frame by $t$ with $t\in[1,\infty]$ gives a path to from the image of $W$ to the empty frame.

Now choose a smooth bump function $\phi: S^q\rightarrow [0,\infty]$ which is identically $\infty$ on a ball $*\in B^\prime \subset  S^q$ and has support in a larger ball $B^\prime \subset B\subset S^q$.  Choose a trivialization of $E_{\mid_B} \rightarrow B$.  This choice along with the framing of $W$ gives a map $B \xrightarrow{s^\prime} \Omega^\infty TF_dO(d)$.  The map $\phi s^\prime $, scalar-multiplying each scaling-frame, extends as the zero-frame to a map $s: S^q\rightarrow\Omega^\infty TF_dO(d)$.  This map $s$ is homotopic to the zero-section map from two paragraphs above.  Note that $s$ sends ${B^\prime\ni *}$ to the empty frame, that is, to the base point in $\Omega^\infty TF_dO(d)$.

Choose a trivialization of $E_{\mid_{{S^q\setminus *}} }\rightarrow {S^q\setminus *}$.  This choice along with the framing of $W$ gives a lift $\tilde{s}:{S^q\setminus *} \rightarrow \Omega^\infty TF_d1(d)$ of $s_{\mid_{{S^q\setminus *}}}$.  But this lift can be extended to all of $S^q$ by declaring $*\mapsto \emptyset$, the empty scaling-frame.  This shows that the map $S^q\xrightarrow{\chi_d\circ \iota\circ f} \Omega^\infty TF_dO(d)$ factors through $\Omega^\infty TF_d1(d)$ up to homotopy as desired.  

\end{proof}

\begin{cor}
For $W^3$ a closed oriented $3$-manifold, the map
\[
BDiff^{\mathcal{O}r}(W)\rightarrow \Omega^\infty MTSO(3)
\]
is the zero-map in rational homotopy.  
\end{cor}

\begin{remark}
Recently, in~\cite{ebert:vanishing} Ebert proved using index theory that the above map is the zero-map in rational homology.
\end{remark}

\bibliographystyle{alpha}
\bibliography{references}

\def\cprime{$'$}
\begin{thebibliography}{GMTW06}

\bibitem[AB83]{atiyah-bott:yang-mills}
M.~F. Atiyah and R.~Bott.
\newblock The {Y}ang-{M}ills equations over {R}iemann surfaces.
\newblock {\em Philos. Trans. Roy. Soc. London Ser. A}, 308(1505):523--615,
  1983.

\bibitem[Ada74]{adams:stable-homotopy}
J.~F. Adams.
\newblock {\em Stable homotopy and generalised homology}.
\newblock University of Chicago Press, Chicago, Ill., 1974.
\newblock Chicago Lectures in Mathematics.

\bibitem[Aya08]{ayala:param-segal}
David Ayala.
\newblock Homological stability for the moduli space of holomorphic curves in
  complex projective space.
\newblock Preprint available on math.stanford.edu/~ayala, on the ArXiv
  (November), 2008.

\bibitem[BD95]{baez-dolan:cobordism-hypothesis}
John~C. Baez and James Dolan.
\newblock Higher-dimensional algebra and topological quantum field theory.
\newblock {\em J. Math. Phys.}, 36(11):6073--6105, 1995.

\bibitem[BP72]{barratt-priddy:symmetric-group}
Michael Barratt and Stewart Priddy.
\newblock On the homology of non-connected monoids and their associated groups.
\newblock {\em Comment. Math. Helv.}, 47:1--14, 1972.

\bibitem[CJS95]{cohen-jones-segal:morse}
Ralph~L. Cohen, John~D.S. Jones, and Graeme~B. Segal.
\newblock Morse theory and classifying spaces.
\newblock Available on online at http://math.stanford.edu/~ralph/morse.ps,
  1995.

\bibitem[DK90]{donaldson-kronheimer:four-manifolds}
S.~K. Donaldson and P.~B. Kronheimer.
\newblock {\em The geometry of four-manifolds}.
\newblock Oxford Mathematical Monographs. The Clarendon Press Oxford University
  Press, New York, 1990.
\newblock Oxford Science Publications.

\bibitem[Dwy96]{dwyer:centralizer}
W.~G. Dwyer.
\newblock The centralizer decomposition of {$BG$}.
\newblock In {\em Algebraic topology: new trends in localization and
  periodicity ({S}ant {F}eliu de {G}u\'\i xols, 1994)}, volume 136 of {\em
  Progr. Math.}, pages 167--184. Birkh\"auser, Basel, 1996.

\bibitem[Ebe09]{ebert:vanishing}
Johannes Ebert.
\newblock A vanishing theorem for characteristic classes of odd-dimensional
  manifold bundles.
\newblock On the ArXiv, 2009.

\bibitem[EM02]{eliashberg:h-principle}
Y.~Eliashberg and N.~Mishachev.
\newblock {\em Introduction to the {$h$}-principle}, volume~48 of {\em Graduate
  Studies in Mathematics}.
\newblock American Mathematical Society, Providence, RI, 2002.

\bibitem[Gal06]{galatius:graphs}
S{\o}ren Galatius.
\newblock Stable homology of automorphism groups of free groups.
\newblock On the ArXiv, 2006.

\bibitem[GCK08a]{cohen-galatius-kitchloo:flat-connections}
S{\o}ren Galatius, Ralph Cohen, and Nitu Kitchloo.
\newblock Universal moduli spaces of surfaces with flat connections and
  cobordism theory.
\newblock On the ArXiv, 2008.

\bibitem[GCK08b]{galatius-cohen-kitchloo:flat-connections}
S{\o}ren Galatius, Ralph Cohen, and Nitu Kitchloo.
\newblock Universal moduli spaces of surfaces with flat connections and
  cobordism theory.
\newblock On the ArXiv, 2008.

\bibitem[GMTW06]{galatius-madsen-tillmann-weiss:cobordism}
S{\o}ren Galatius, Ib~Madsen, Ulrike Tillmann, and Michael Weiss.
\newblock The homotopy type of the cobordism category.
\newblock On the ArXiv, 2006.

\bibitem[Gro69]{gromov:foliations}
M.~L. Gromov.
\newblock Stable mappings of foliations into manifolds.
\newblock {\em Izv. Akad. Nauk SSSR Ser. Mat.}, 33:707--734, 1969.

\bibitem[Gro86]{gromov:h-principle}
Mikhael Gromov.
\newblock {\em Partial differential relations}, volume~9 of {\em Ergebnisse der
  Mathematik und ihrer Grenzgebiete (3) [Results in Mathematics and Related
  Areas (3)]}.
\newblock Springer-Verlag, Berlin, 1986.

\bibitem[Har85]{harer:stability}
John~L. Harer.
\newblock Stability of the homology of the mapping class groups of orientable
  surfaces.
\newblock {\em Ann. of Math. (2)}, 121(2):215--249, 1985.

\bibitem[HV98]{hatcher-vogtmann:splitting-conjecture}
Allen Hatcher and Karen Vogtmann.
\newblock Rational homology of {${\rm Aut}(F\sb n)$}.
\newblock {\em Math. Res. Lett.}, 5(6):759--780, 1998.

\bibitem[HW07]{hatcher-wahl:stability}
Allen Hatcher and Nathalie Wahl.
\newblock Stabilization for mapping class groups of 3-manifolds.
\newblock On the ArXiv, 2007.

\bibitem[KS77]{kirby-siebenmann:topological-manifolds}
Robion~C. Kirby and Laurence~C. Siebenmann.
\newblock {\em Foundational essays on topological manifolds, smoothings, and
  triangulations}.
\newblock Princeton University Press, Princeton, N.J., 1977.
\newblock With notes by John Milnor and Michael Atiyah, Annals of Mathematics
  Studies, No. 88.

\bibitem[Lur09]{hopkins-lurie:tfts}
Jacob Lurie.
\newblock On the classification of topological field theories.
\newblock online at http://math.mit.edu/~lurie/, 2009.

\bibitem[McD75]{mcduff:configurations-scanning}
Dusa McDuff.
\newblock Configuration spaces of positive and negative particles.
\newblock {\em Topology}, 14:91--107, 1975.

\bibitem[Mor96]{morgan:seiberg-witten}
John~W. Morgan.
\newblock {\em The {S}eiberg-{W}itten equations and applications to the
  topology of smooth four-manifolds}, volume~44 of {\em Mathematical Notes}.
\newblock Princeton University Press, Princeton, NJ, 1996.

\bibitem[MS76]{mcduff-segal:group-completion}
D.~McDuff and G.~Segal.
\newblock Homology fibrations and the ``group-completion'' theorem.
\newblock {\em Invent. Math.}, 31(3):279--284, 1975/76.

\bibitem[Mum83]{mumford:conjecture}
David Mumford.
\newblock Towards an enumerative geometry of the moduli space of curves.
\newblock In {\em Arithmetic and geometry, {V}ol. {II}}, volume~36 of {\em
  Progr. Math.}, pages 271--328. Birkh\"auser Boston, Boston, MA, 1983.

\bibitem[MW07]{madsen-weiss:mumford}
Ib~Madsen and Michael Weiss.
\newblock The stable moduli space of {R}iemann surfaces: {M}umford's
  conjecture.
\newblock {\em Ann. of Math. (2)}, 165(3):843--941, 2007.

\bibitem[Phi67]{phillips:submersions}
Anthony Phillips.
\newblock Submersions of open manifolds.
\newblock {\em Topology}, 6:171--206, 1967.

\bibitem[Rez01]{rezk:segal-spaces}
Charles Rezk.
\newblock A model for the homotopy theory of homotopy theory.
\newblock {\em Trans. Amer. Math. Soc.}, 353(3):973--1007 (electronic), 2001.

\bibitem[RW08]{randal-williams:embedding}
Oscar Randal-Williams.
\newblock Embedded cobordism.
\newblock Oxford Transfer Thesis - to appear, 2008.

\bibitem[Sco05]{scorpan:wild-world}
Alexandru Scorpan.
\newblock {\em The wild world of 4-manifolds}.
\newblock American Mathematical Society, Providence, RI, 2005.

\bibitem[Seg74]{segal:gamma-spaces}
Graeme Segal.
\newblock Categories and cohomology theories.
\newblock {\em Topology}, 13:293--312, 1974.

\bibitem[Seg78]{segal:classifying-spaces}
Graeme Segal.
\newblock Classifying spaces related to foliations.
\newblock {\em Topology}, 17(4):367--382, 1978.

\bibitem[Seg79]{segal:scanning}
Graeme Segal.
\newblock The topology of spaces of rational functions.
\newblock {\em Acta Math.}, 143(1-2):39--72, 1979.

\bibitem[Til99]{tillmann:mapping-class-group}
Ulrike Tillmann.
\newblock A splitting for the stable mapping class group.
\newblock {\em Math. Proc. Cambridge Philos. Soc.}, 127(1):55--65, 1999.

\bibitem[Wei05]{weiss:classify}
Michael Weiss.
\newblock What does the classifying space of a category classify?
\newblock {\em Homology Homotopy Appl.}, 7(1):185--195 (electronic), 2005.

\end{thebibliography}

\end{document}